\documentclass[12pt]{article}

\usepackage[margin=2.5cm, top=2.5cm, bottom=2.5cm]{geometry}

\usepackage[utf8]{inputenc} 
\usepackage[T1]{fontenc}    
\usepackage{url}            
\usepackage{booktabs}       
\usepackage{amsfonts}       
\usepackage{nicefrac}       
\usepackage{microtype}      

\usepackage[sort&compress,numbers]{natbib}

\usepackage{amsmath,amssymb,amsthm,xcolor}
\usepackage{aliascnt} 
\usepackage{array}
\usepackage{float}
\usepackage[hidelinks,hypertexnames=false]{hyperref}
\usepackage[capitalize]{cleveref}
\theoremstyle{definition}
\newtheorem{example}{Example}[section]

\newaliascnt{definition}{example}
\newtheorem{definition}[definition]{Definition}
\aliascntresetthe{definition}
\crefname{definition}{Definition}{Definitions}

\newaliascnt{remark}{example}
\newtheorem{remark}[remark]{Remark}
\aliascntresetthe{remark}
\crefname{remark}{Remark}{Remarks}

\theoremstyle{plain}

\newaliascnt{theorem}{example}
\newtheorem{theorem}[theorem]{Theorem}
\aliascntresetthe{theorem}
\crefname{theorem}{Theorem}{Theorems}

\newaliascnt{assumption}{example}
\newtheorem{assumption}[assumption]{Assumption}
\aliascntresetthe{assumption}
\crefname{assumption}{Assumption}{Assumptions}

\newaliascnt{lemma}{example}
\newtheorem{lemma}[lemma]{Lemma}
\aliascntresetthe{lemma}
\crefname{lemma}{Lemma}{Lemmas}

\newaliascnt{proposition}{example}
\newtheorem{proposition}[proposition]{Proposition}
\aliascntresetthe{proposition}
\crefname{proposition}{Proposition}{Propositions}

\newaliascnt{corollary}{example}
\newtheorem{corollary}[corollary]{Corollary}
\aliascntresetthe{corollary}
\crefname{corollary}{Corollary}{Corollaries}

\usepackage{graphicx}
\usepackage{subcaption}
\usepackage{mdframed}
\usepackage{lipsum}
\usepackage{wrapfig}
\usepackage{tikz}
\usetikzlibrary{arrows.meta}
\usepackage{enumitem}
\allowdisplaybreaks
\usepackage{epstopdf}
\usepackage{booktabs}
\usepackage{changepage}
\usepackage{multirow}
\usepackage{algorithm}
\usepackage{algorithmic}
\usepackage{mathtools}
\crefname{assumption}{Assumption}{Assumptions}

\numberwithin{equation}{section}
\usepackage[nottoc]{tocbibind}
\usepackage{mathrsfs}
\usepackage{stmaryrd}

\newcommand{\prox}{{\rm prox}}
\newcommand{\dom}{\mathop{\rm dom}} 
 
\newcommand{\epi}{\mathop{\rm epi}} 
\newcommand{\co}{\mathop{\rm co}}

\newcommand{\gra}{{\rm graph}}

\newcommand{\argmin}{\mathop{\mathrm{argmin}}}

\title{\LARGE Proximal random reshuffling under local Lipschitz continuity}

\begin{document}

\author{\large C\'edric Josz\thanks{\url{cj2638@columbia.edu}, IEOR, Columbia University, New York. Research supported by NSF EPCN grant 2023032 and ONR grant N00014-21-1-2282.} \and Lexiao Lai\thanks{\url{lai.lexiao@hku.hk}, Department of Mathematics, The University of Hong Kong, Hong Kong. Research supported by Hong Kong Research Grants Council (ECS project 27301425) and a start-up fund from The University of Hong Kong.} \and Xiaopeng Li\thanks{\url{xiao pengli1@cuhk.edu.cn}, School of Artificial Intelligence, The Chinese University of Hong Kong, Shenzhen, Shenzhen. Research supported by the Guangdong Provincial Key Laboratory of Big Data Computing at The Chinese University of Hong Kong, Shenzhen.}} 
\date{}

\maketitle
\vspace*{-5mm}
\begin{center}
    \textbf{Abstract}
    \end{center}
    \vspace*{-4mm}
 \begin{adjustwidth}{0.2in}{0.2in}
~~~~ We study proximal random reshuffling (PRR) for minimizing the sum of locally Lipschitz or locally smooth functions and a proper lower semicontinuous convex function without assuming coercivity or the existence of limit points. The algorithmic guarantees pertaining to near approximate stationarity rely on a new tracking lemma linking the iterates to trajectories of conservative fields. One of the novelties in the analysis consists in handling set-valued mappings with unbounded values. In the locally smooth case, it improves the known convergence rate from nearly $O(k^{-1/4})$ to nearly $o(k^{-1/2})$.
\end{adjustwidth} 
\vspace*{3mm}
\noindent{\bf Keywords:} Nonsmooth nonconvex optimization, random reshuffling, conservative fields, differential inclusions. 
\vspace*{-3mm}

\tableofcontents

\section{Introduction}\label{sec:Introduction}

There is a growing interest in solving optimization problems where either the objective function or its gradient is merely locally Lipschitz continuous, rather than globally Lipschitz. One approach, popular in the last five years, has been to modify classical algorithms using line search methods (\citet{kanzow2022convergence,jia2023convergence,lu2023accelerated,zhang2024first}), normalized search directions (\citet{grimmer2019convergence,li2025revisiting}), non-Euclidean metrics (\citet{bauschke2017descent,bolte2018first,bauschke2019linear,dragomir2022optimal}), or radial projections (\citet{renegar2016efficient,grimmer2018radial}). Another approach, and the one that we will follow in this paper, is to resort to continuous-time dynamics (\citet{davis2020stochastic,bolte2021conservative,bianchi2022convergence,josz2024global}) to study classical algorithms without making any modifications. The latter approach seems to offer the only avenue so far that succeeds without assuming coercivity, boundedness of iterates, or the existence of limit points. Without these assumptions, it is however limited to deterministic algorithms for unconstrained optimization, namely the gradient method (\citet{josz2023global}) and the momentum method (\citet{josz2023convergence}).

The goal of this paper is to take a first step towards stochastic and constrained optimization by considering the simplest implementation of the stochastic subgradient method, along with constraints of a particularly simple kind. Specifically, we consider proximal random reshuffling (i.e., \cref{alg:prr}) for solving the composite model
\begin{equation}
\label{eq:obj}
	\inf\limits_{x\in \mathbb{R}^n} \Phi(x) := \sum\limits_{i= 1}^N f_i(x) + g(x)
\end{equation}
where $f_1,\hdots,f_N:\mathbb{R}^n \rightarrow \mathbb{R}$ are either locally Lipschitz or differentiable with locally Lipschitz gradients, and $g:\mathbb{R}^n \rightarrow \overline{\mathbb{R}}:= \mathbb{R} \cup \{\pm\infty\}$ is proper, convex, and locally Lipschitz on its closed domain. We allow the user to choose any conservative field $D_i:\mathbb{R}^n \rightrightarrows \mathbb{R}^n$ for each $f_i$, for example the gradient $\nabla f_i$ if $f_i$ is continuously differentiable, or the Clarke subdifferential $\partial f_i$ if $f_i$ is semialgebraic (\citet[Definition 1.4]{pham2016genericity}), but it could also be the output of an automatic differentiation procedure (see, e.g., \citet[Theorem 8]{bolte2021conservative}). 

\begin{algorithm}[ht]
\caption{Proximal random reshuffling (PRR)}\label{alg:prr}
\begin{algorithmic}
\STATE{\textbf{choose} $\alpha_0,\alpha_1,\alpha_2,\hdots>0$, $x_0 \in \dom \Phi$.}
\FOR{$k = 0,1,2, \ldots$}
    \STATE{ $x_{k,0} = x_k$}
    \STATE{choose a permutation $\sigma^k$ of $\llbracket 1 , N \rrbracket$}
    \FOR{$i = 1, \ldots, N$}
        \STATE{ $x_{k,i} \in x_{k,i-1}  - \alpha_k D_{\sigma^k_i} (x_{k,i-1})$}
    \ENDFOR
    \STATE{ $x_{k+1} = \prox_{\alpha_k g}(x_{k,N})$}
\ENDFOR
\end{algorithmic}
\end{algorithm}

Inspired by the pioneering works of  \citet{ljung1977analysis,kushner1977general,kushner1977convergence,benaim2005stochastic,benaim2006stochastic,benaim1999dynamics}, we consider the differential inclusion
\begin{equation}
    \label{eq:di}
    x'(t) \in -(D_1+\cdots+D_N+\partial g)(x(t)), ~~~\mathrm{for~a.e.}~ t>0,
\end{equation}
where $x:[0,\infty)\rightarrow \mathbb{R}^n$ is an absolutely continuous function. Their insight, extended to the proximal setting by \citet{davis2020stochastic} (i.e., $N=1$, $D_1 = \partial f_1$), is that trajectories of \eqref{eq:di} track the iterates $(x_k)_{k\in\mathbb{N}}$ generated by \cref{alg:prr} in the following sense. For any accuracy $\epsilon>0$ and any time $T>0$, after sufficiently many iterations, say $l$, we have $\|x_k - x(t_k-t_l)\| \leqslant \epsilon$ for some solution $x(\cdot)$ to \eqref{eq:di}, where the times $t_k := \alpha_0+\hdots+\alpha_{k-1}$ lie in $ [t_l , t_l + T]$. This requires bounded iterates, a regularity condition that is slightly stronger than local Lipschitz continuity, and a specific choice of step sizes.

However, such a bounded-iterate assumption can be artificial and unnecessary for understanding the fundamental behavior of the algorithm. In the absence of a global gradient Lipschitz constant, even if iterates remain bounded in practice, gradient descent (a special case of \cref{alg:prr}) may fail to converge (see \cref{fig:unbdd_const} below). This occurs because the iterates are not uniformly bounded for all sufficiently small step sizes, and the sufficient decrease condition (see \citet[H1]{attouch2013convergence}) does not hold. A similar issue arises in machine learning problems such as $\ell_1$ matrix completion (\cref{ex:nrpca} without constraints), where constant step sizes lead to non-uniformly bounded iterates (see \cref{fig:constant}). Moreover, while many well-tuned algorithms exhibit bounded iterates in practice, there exist natural scenarios where unbounded iterates arise. For instance, in the same $\ell_1$ matrix completion problem, using diminishing step sizes can lead to unbounded iterates (see \cref{fig:unbdd_iter}). Finally, verifying boundedness empirically is inherently difficult. As shown in \cref{fig:unbdd_summable}, even with summable step sizes, after running $10^9$ iterations, the iterate norm appears to increase, but it does so at such a slow rate that it is unclear whether it is truly unbounded or simply growing at a diminishing rate and will eventually stabilize. This highlights the need for theoretical results that do not rely on bounded iterate assumptions. 
\begin{figure}[!ht]
\centering 
\includegraphics[width=0.33\textwidth]{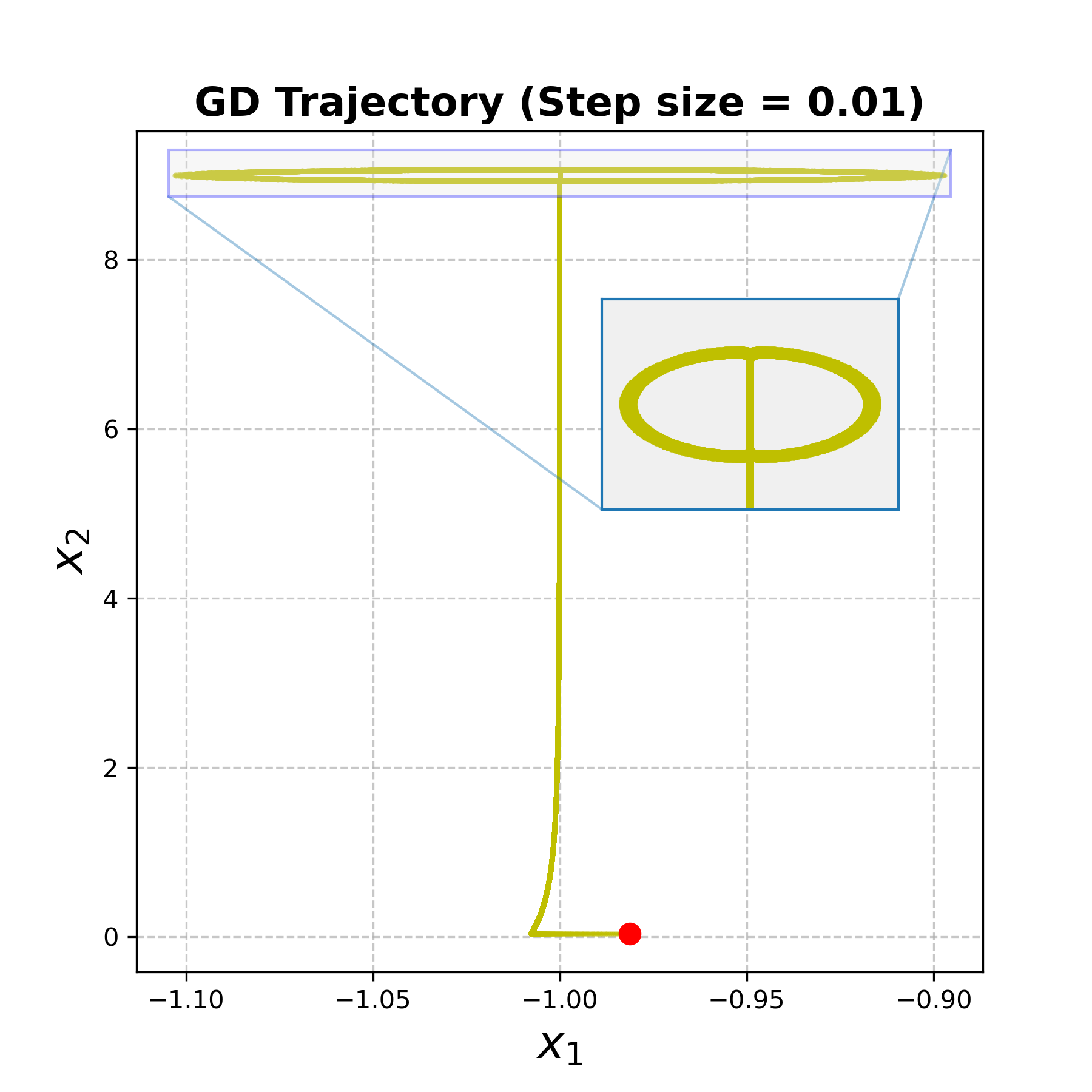}\hspace{-1ex}
\includegraphics[width=0.33\textwidth]{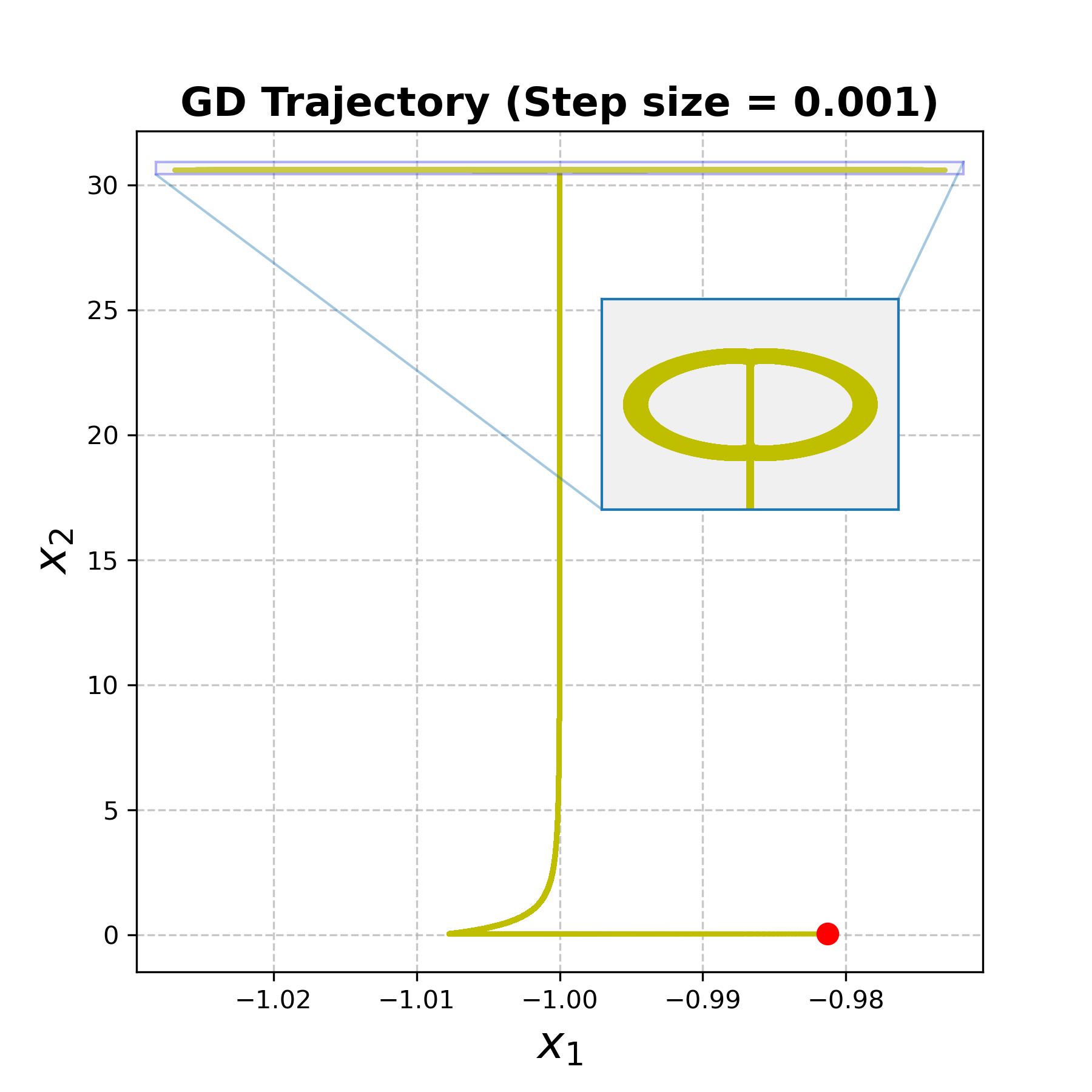}\hspace{-1ex}
\includegraphics[width=0.33\textwidth]{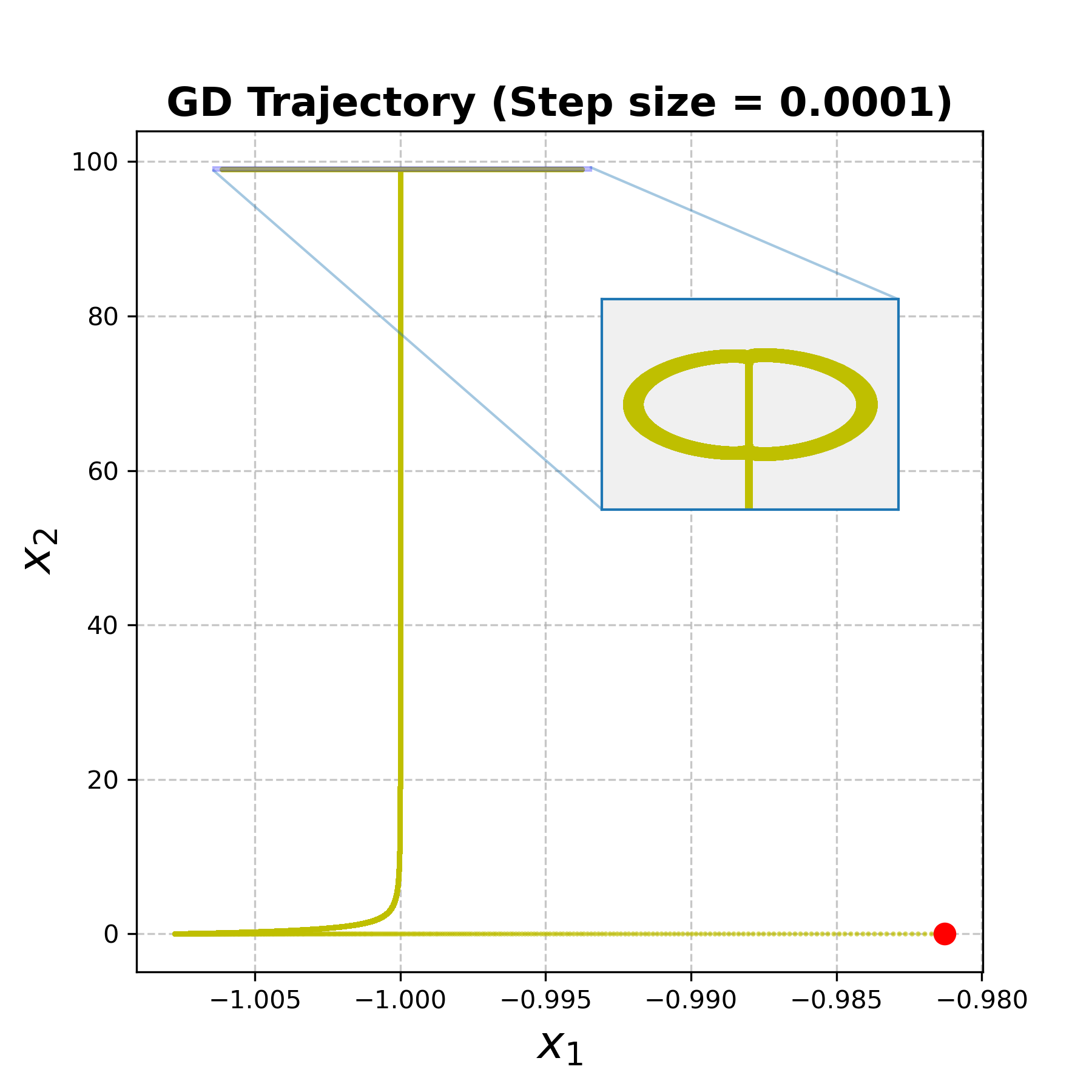}
\caption{Plots of gradient descent iterates for $f(x_1,x_2)=(x_1+1)^2(x_2+1)^2 + \exp\left(-(1+x_1^2)^{1/8}-(1+x_2^2)^{1/8}\right)$ with a common initial point and different constant step sizes. A similar behavior can be observed for a simpler but not lower bounded function $f(x_1,x_2)=x_1^2x_2^2-x_1-x_2$.} 
\label{fig:unbdd_const}
\end{figure}

In this work, we propose a tracking lemma devoid of such assumptions, which is at the same time more intuitive. For any initial iterate in some bounded set, we show that $\|x_k - x(t_k)\| \leqslant \epsilon$ for all times $t_k$ in $[0,T]$, provided that the step sizes are sufficiently small. Notably, this result extends \citet{josz2024global} in a nontrivial way, even if we ignore the constraints, as it applies over any finite time horizon without requiring boundedness of iterates. In the general case, the dynamics is governed by the set-valued mapping $D_1+\cdots+D_N+\partial g$ with possibly unbounded values. This leads us to extend the notion due to \citet{bolte2021conservative} initially conceived with compact values. Some set-valued analysis reveals that the local boundedness of $x \in \mathbb{R}^n \mapsto \inf \{ \|y\| : y \in (D_1+\cdots+D_N+\partial g)(x)\}$ plays a role instead. The idea of selecting the minimal norm element from unbounded set-valued mappings is well-established in convex subgradient dynamics (\citet{brezis1973operateurs}). In nonconvex settings, some authors (e.g., \citet{schechtman2023stochastic}, \citet{cornet1983existence}) have instead truncated unbounded set-valued mappings, a concept closely linked to the local boundedness of their distance to the origin.

The tracking lemma enables us to ensure that \cref{alg:prr} can recover stationary points of the objective function $\Phi$ from any initial point, with varying degrees of precision depending on the regularity of $f_1,\hdots,f_N$ and $g$. If $f_1,\ldots,f_N$ are locally Lipschitz, by a direct application of the tracking lemma, we establish that \cref{alg:prr} reaches $(\epsilon,\delta)$-near approximate stationary points (\cref{def:nas}) under mild conditions (\cref{thm:exists_small}). If the objective is weakly convex, then the algorithm converges to $(\epsilon,\delta)$-near approximate stationary points when certain summable step sizes are used (\cref{thm:weakly_convex}). If additionally requiring the uniqueness and convergence of the solution to \eqref{eq:di}, then the algorithm converges to $(\epsilon,0)$-near approximate stationary points (\cref{thm:nonsmooth_bounded_flow}). The power of the tracking lemma becomes especially apparent when $f_1, \dots, f_N $ are locally smooth. Inspired by \citet{josz2023global} and \citet{li2023convergence}, we establish that if $\Phi$ is semialgebraic and has bounded subgradient trajectories, then \cref{alg:prr} converges to stationary points at a rate approaching $o(1/\sqrt{k})$ (\cref{thm:sgd_finite_len}). This rate improves to $o(1/k)$ in the deterministic setting ($N=1$), where \cref{alg:prr} reduces to the proximal gradient method. Applied to nonnegative matrix factorization (NMF) (\citet{paatero1994positive,lee1999learning,wang2012nonnegative}), it provides the first convergence guarantee to stationary points of NMF from any initial point via the projected gradient method without assuming the existence of limit points. Our sole assumption is to use sufficiently small constant step sizes -- how small depends on a numerically intractable quantity. To the best of our knowledge, other iterative methods all require the existence of limit points to ensure convergence of the iterates, but some benefit from an explicit and tractable choice of constant step sizes (see \cref{ex:ngmf} and its subsequent discussion for details).

This paper is organized as follows. \cref{sec:Literature review} provides a detailed review of related work, comparing our results and techniques with existing approaches to highlight our contributions. \cref{sec:main results} contains the main results and is divided into several subsections. They respectively contain the definitions, assumptions, theorems, and examples. \cref{sec:proofs} contains the proofs. It begins with main technical lemma regarding tracking by trajectories of conservative fields. The four main results are then proved successively using this lemma. Finally, we conclude in \cref{sec:conclusion} with a summary of our findings, a discussion of limitations, and potential directions for future research. Omitted proofs and supplementary figures are provided in \cref{sec:omittedproof,sec:suppfigs}.

\section{Literature review}
\label{sec:Literature review}

This section reviews the literature closely related to \cref{alg:prr}, comparing their results and assumptions to ours. We divide the discussion into two subsections based on the regularity of the component functions $f_i$'s: nonsmooth and locally smooth. At the end of each subsection, we also compare the proof techniques used in our analysis with those in the related works. 

\subsection{Nonsmooth component functions}\label{subsec:nonsmoothliterature}

We start by reviewing the literature that studies algorithms which tackle \eqref{eq:obj} with locally Lipschitz component functions $f_i$'s. To the best of our knowledge, the proximal random reshuffling algorithm (\cref{alg:prr}) has not been studied at this level of generality, but there are other algorithms which apply to settings similar to \eqref{eq:obj} with various assumptions. Many of their techniques serve as inspirations for our own approach. A recent work of \citet{solodov2025convergence} studies a class of proximal algorithms, which includes a special case of \cref{alg:prr} without permutations. If the component functions are Clarke regular (\citet{clarke1990}) and the iterates are bounded, then the sequence converges subsequentially to a critical point. The proximal stochastic subgradient method (with replacement) (\citet{duchi2009efficient,davis2019stochastic,davis2020stochastic}) is closely related to \cref{alg:prr}, whose updates consist of a step of stochastic subgradient method with respect to $f:= f_1+\cdots+f_N$ and a proximal step with respect to $g$. Convergence of the method was shown in \citet{duchi2009efficient} if 
$f$ is convex and certain Lipschitz conditions (\citet[equation (6)]{duchi2009efficient}) are satisfied. If one relaxes the convexity and assumes that $f_1,\ldots, f_N$ are Lipschitz continuous and $\Phi$ is weakly convex instead, then the method generates a point with small gradient norm of a Moreau envelope of $\Phi$ in expectation (\citet[Theorem 3.4]{davis2019stochastic}). If one further relaxes the regularity assumptions and only assumes that $\Phi$ is Whitney stratifiable, then the iterates subsequentially converge to composite critical points, given that they are bounded almost surely, among other assumptions (\citet[Theorem 6.2]{davis2020stochastic} and \citet{schechtman2023stochastic}). 

We next turn to the other subgradient-based algorithms that apply to general nonsmooth and nonconvex optimization problems, which might not admit the same composite structure as in \eqref{eq:obj}. This includes variants of the subgradient method that allow stochasticity (\citet{davis2020stochastic,josz2024global}), inexact subgradient oracles (\citet{bolte2025inexact}), and/or momentum terms in the updates (\citet{josz2024global,ding2025adamfamily,xiao2026stochastic}). These algorithms are analyzed in settings where \cref{assumption_standing} always holds. Moreover, it is worth noticing that all of these works assume either the iterates are bounded or the objective function is coercive. Under such assumptions, the algorithms are shown to converge subsequentially to critical points with diminishing step sizes (\citet{davis2020stochastic,ding2025adamfamily,bolte2025inexact}), or eventually stay close to the set of critical points with non-diminishing step sizes (\citet{josz2024global,xiao2026stochastic,bolte2025inexact}). In terms of the stationary measure in this work (see \cref{def:nas}), for any $\epsilon>0$, these algorithms return $(\epsilon,0)$-near approximate stationary (NAS) points after a sufficiently large number of iterations, given that the step sizes are sufficiently small. The analyses in these works rely on the continuous-time limit of the iterates as we have discussed in the introduction. 

We discuss the relevance of \cref{thm:exists_small,thm:weakly_convex,thm:nonsmooth_bounded_flow} in light of the previous results. First, we avoid assuming coercivity of the objective function or boundedness of the iterates, which allows us to apply these theorems to applications in data science (see \cref{subsec:example}). Specifically, \cref{thm:exists_small} establishes that an approximate stationary point can be reached under mild assumptions (\cref{assumption_standing}). While this result may seem limited, it cannot be strengthened without additional hypotheses. To analyze the asymptotic behavior of \cref{alg:prr}, we rely on summable step sizes in \cref{thm:weakly_convex,thm:nonsmooth_bounded_flow} under stricter assumptions. These theorems follow directly from the tracking lemma (\cref{lemma:stoc_tracking}), proven in \cref{subsec:proof_of_stoc_tracking}. Notably, \cref{thm:weakly_convex} provides new asymptotic guarantees for locally Lipschitz weakly convex functions, advancing beyond prior works (\citet{davis2019stochastic,solodov2025convergence}). Meanwhile, \cref{thm:nonsmooth_bounded_flow} appears to be the first result leveraging bounded trajectory assumptions for subgradient-based methods, whereas prior studies applied this condition only to gradient and momentum methods (\citet{josz2023global,josz2023convergence}) for smooth functions.

\subsection{Locally smooth component functions}\label{subsec:smoothliterature}

We now discuss the contributions of our results when the component functions $f_i$ are locally smooth, and relate them to existing work. We divide our discussion into the deterministic case (\cref{alg:prr} with $N=1$) and the stochastic case (\cref{alg:prr} with $N>1$), with the latter further divided based on the presence or absence of the nonsmooth term $g$. We conclude by highlighting the technical contributions of \cref{thm:sgd_finite_len} in comparison to the literature. 

In the deterministic case (\cref{alg:prr} with $N=1$), the algorithm reduces to the proximal gradient method. \cref{thm:sgd_finite_len} establishes that, if $\Phi$ is locally smooth semialgebraic and has bounded subgradient trajectories, the iterates converge to stationary points at a rate of $o(1/k)$ under sufficiently small constant step sizes. This rate is faster than the general $O(1/\sqrt{k})$ convergence for nonconvex globally smooth functions when $\Phi$ attains its infimum (\citet[Theorem 10.15]{beck2017first}, \citet[Theorem 3]{nesterov2013gradient}), and is consistent with the known $O(1/k)$ rate in the convex case (\citet[Theorem 10.26]{beck2017first}, \citet[Theorem 4]{nesterov2013gradient}). The global smoothness assumption can be relaxed to local smoothness if a sufficient decrease condition is satisfied (\citet[H1]{attouch2013convergence}). Under this condition, convergence to critical points holds (\citet[Theorem 3.1]{kanzow2022convergence}), and the iterates have finite length when a limit point exists and $\Phi$ satisfies the Kurdyka-\L{}ojasiewicz inequality (\citet[Theorem 5.1]{attouch2013convergence}, \citet[Theorem 3.1]{frankel2015splitting}, \citet[Theorem 4.5]{jia2023convergence}). Furthermore, explicit convergence rates can be obtained if the desingularizing function in the K\L{} inequality is a power function (\citet[Theorem 3.4]{frankel2015splitting}, \citet[Theorem 4.6]{jia2023convergence}).

If $N > 1$, \cref{thm:sgd_finite_len} shows that if $\Phi$ is locally smooth semialgebraic and has bounded subgradient trajectories, then the PRR algorithm (\cref{alg:prr}) with diminishing step sizes $\alpha_k = \alpha/(k+1)^\beta$ converges to stationary points at a rate $o(1/k^{1-\beta})$ for any $\beta \in (1/2,1)$ and sufficiently small $\alpha$. When $g = 0$, \cref{alg:prr} corresponds to the random reshuffling (RR) algorithm, also known as stochastic gradient descent (SGD) without replacement (\citet{bertsekas2011incremental, mishchenko2020random, gurbuzbalaban2021random}). If a fixed permutation $\sigma^k$ is used for all $k \in \mathbb{N}$, the method reduces to the classical incremental gradient (IG) algorithm, whose asymptotic behavior has been extensively studied (\citet{luo1991convergence, Mangasarian01011994, bertsekas2000gradient, gurbuzbalaban2019convergence}). RR has been shown to converge faster than IG or SGD with replacement when the $f_i$ are convex with globally Lipschitz gradients and Hessians (\citet{gurbuzbalaban2021random, haochen2019random}). In the nonconvex setting with globally smooth $f_i$, RR achieves $O(1/\epsilon^3)$ complexity using $O(\epsilon)$ constant step sizes under certain expected smoothness conditions (\citet{khaled2022better}), as shown in \citet[Corollary 3]{mishchenko2020random} and \citet[Corollary 1]{nguyen2021unified}. The same complexity is obtained using optimized diminishing step sizes (\citet[Theorem 6]{nguyen2021unified}). When $f_i$ are semialgebraic and iterates are bounded, RR with the same step sizes as in \cref{thm:sgd_finite_len} converges to stationary points (\citet[Corollary 3.8]{li2023convergence}) at the rate $O(k^{-(3\beta - 1)/2})$ for $\beta$ close to $1/2$ (\citet[Theorem 5.3]{qiu2025new}), approaching $O(k^{-1/4})$, which can be improved with knowledge of a \L{}ojasiewicz exponent.

When $g \ne 0$, related results for PRR and proximal SGD are more limited. For convex, globally smooth $f_i$, \citet{rosasco2020convergence} show proximal SGD converges at nearly $O(1/\sqrt{k})$, matching our rate for nonconvex $f_i$. \citet{majewski2018analysis} study projected proximal SGD with $C^1$ $f$, and $g$ locally Lipschitz, regular, and lower bounded, assuming a compact constraint set. Their framework handles broader function classes without requiring definability or convexity of $g$, but critically relies on boundedness of the domain, making it inapplicable to problems like NMF. In contrast, PRR achieves better iteration complexity, $O(1/\epsilon^3)$ versus $O(1/\epsilon^4)$ for proximal SGD, when $f_i$ are globally smooth and additional regularity conditions hold (\citet[Section 4]{mishchenko2022proximal}, \citet{davis2019stochastic}). A normal-map-based PRR algorithm was proposed in \citet{qiu2025new}, showing similar improvements over the proximal stochastic subgradient method under analogous assumptions. By further assuming bounded iterates and a Kurdyka-\L{}ojasiewicz property, the asymptotic behavior of the iterates is also established (\citet[Sections 4 and 5]{qiu2025new}).

We conclude by highlighting the contributions of \cref{thm:sgd_finite_len} in the context of related work. Global convergence of deterministic first-order methods for smooth semialgebraic functions is well established (\citet{absil2005convergence, attouch2013convergence, ochs2014ipiano, josz2023global, josz2023convergence}), while convergence of SGD with replacement has been studied in \citet{tadic2015convergence} and \citet{dereich2024convergence}. The analysis of RR is more recent in \citet{li2023convergence}. Although \citet{josz2023global} and \citet{josz2023convergence} also leverage bounded subgradient trajectories, their results are restricted to deterministic methods with a local descent property. In contrast, the PRR algorithm analyzed in \cref{thm:sgd_finite_len} does not satisfy such a property and instead relies on an approximate descent condition (\cref{lem:approx_descent}), which complicates the proof of a length formula (\cref{prop:discrete_length_formula,cor:length}). The inclusion of a convex extended-value term $g$ further introduces difficulties due to the potential unboundedness of $\partial g$. Finally, unlike \citet{li2023convergence}, we remove the assumptions of global smoothness and bounded iterates, and extend the analysis from semialgebraic to general definable functions. This generalization is nontrivial in the RR setting, where previous work relies on a quasi-additivity condition for the desingularizing function, an assumption that fails for general definable functions. Our approach circumvents this via a new uniform Kurdyka–\L{}ojasiewicz inequality (\cref{lemma:ukl}).

\section{Main results}\label{sec:main results}

Given two integers $a \leqslant b$, we use the notation $\llbracket a,b \rrbracket := \{a,a+1,\hdots,b\}$. Let $\|\cdot\|$ be the induced norm of an inner product $\langle \cdot, \cdot\rangle$ on $\mathbb{R}^n$. Given $\Phi:\mathbb{R}^n\rightarrow \overline{\mathbb{R}}$, the domain, graph, and epigraph are respectively given by $\dom \Phi := \{x\in \mathbb{R}^n:\Phi(x)<\infty\}$, $\gra \Phi := \{(x,t) \in \mathbb{R}^n \times \mathbb{R}: \Phi(x)=t \}$, and $\epi \Phi := \{(x,t)\in \mathbb{R}^{n+1}:\Phi(x)\leqslant t\}$. A function $\Phi$ is convex (respectively lower semicontinuous) if $\epi \Phi$ is convex (respectively closed). A function $\Phi:\mathbb{R}^n \rightarrow \overline{\mathbb{R}}$ is weakly convex (\citet{nurminskii1973quasigradient,vial1983strong}) if there exists $\rho\geqslant 0$ such that $\Phi + \rho\|\cdot\|^2$ is convex.

Given $g:\mathbb{R}^n\rightarrow \overline{\mathbb{R}}$, the proximal mapping $\prox_g:\mathbb{R}^n \rightrightarrows \mathbb{R}^n$ (\citet[3.b]{moreau1965proximite}) is defined by $\prox_g(x):= \argmin\{y\in \mathbb{R}^n:g(y) + \|y-x\|^2/2\}$ for all $x\in \mathbb{R}^n$. If $g$ is proper, lower semicontinuous and convex, then $\prox_g(x)$ is a singleton for any $x\in \mathbb{R}^n$ (\citet[Theorem 6.3]{beck2017first}) and it lies in the domain of the convex subdifferential of $g$ (\citet[Theorem 6.39]{beck2017first}).

\subsection{Definitions}
\label{subsec:def}

We begin by stating the definition of a solution to a differential inclusion. Let $\mathbb{R}_+ :=[0,\infty)$ and $\mathbb{R}_{++} :=(0,\infty)$.

\begin{definition}
    Let $I$ be an interval of $\mathbb{R}_+$ and $D:\mathbb{R}^n\rightrightarrows \mathbb{R}^n$ be a set-valued mapping. We say that $x:I\rightarrow \mathbb{R}^n$ is a solution to $x' \in - D(x)$ if $x(\cdot)$ is absolutely continuous and $x'(t) \in - D(x(t))$ for almost every $t\in I$.
\end{definition}
We say that $x(\cdot)$ is a $D$-trajectory if it is maximal (see, e.g., \citet[Definition 5]{josz2023global}). Also, we say that $x(\cdot)$ is a subgradient trajectory of $\Phi:\mathbb{R}^n\rightarrow \overline{\mathbb{R}}$ if it is a $\partial \Phi$-trajectory, where $\partial \Phi:\mathbb{R}^n\rightrightarrows\mathbb{R}^n$ denotes the Clarke subdifferential (see, e.g., \citet[p. 336]{rockafellar2009variational}, \citet[Chapter 2]{clarke1990}).

In order to measure the optimality of the iterates generated by \cref{alg:prr}, we will borrow the notion of $(\epsilon,\delta)$-near approximate stationarity (see, e.g., \citet[Definition 4]{tian2024no}, \citet[Definition 2.7]{tian2021hardness}, \citet{davis2019proximally,davis2019stochastic}). It was introduced due to the intractability of finding near stationary points for Lipschitz functions, in the sense that the Clarke subdifferential admits small elements (\citet[Theorem 5]{zhang2020complexity}). It is thus well suited for this work. We remark that this is a stronger stationary notion compared to the one based on the Goldstein subdifferential (\citet[Definition 4]{zhang2020complexity}, \citet{goldstein1977optimization}). Given $S\subset \mathbb{R}^n$, $x\in \mathbb{R}^n$, and $r\geqslant 0$, let $B(x,r):=\{y \in \mathbb{R}^n:\|x-y\|\leqslant r\}$, $d(x,S) := \inf\{\|x-y\|:y\in S\}$, and $P_S(x) := \argmin\{\|x-y\|:y\in S\}$.

\begin{definition}\label{def:nas}
Given $\epsilon,\delta\geqslant 0$ and a set-valued mapping $D:\mathbb{R}^n\rightrightarrows \mathbb{R}^n$, a point $x \in \mathbb{R}^n$ is $(\epsilon,\delta)$-near approximate $D$-stationary if $d(0,D(B(x,\epsilon))) \leqslant \delta$. 
\end{definition}

\begin{remark}
    In particular, a $(0,0)$-near approximate $D$-stationary point is also called a $D$-stationary point or a $D$-critical point. 
\end{remark}

\subsection{Assumptions}

We first introduce the standing assumption for all the main results.
\begin{assumption}
    \label{assumption_standing} \text{ } 
    \begin{enumerate} 
        \item $f_1,\hdots,f_N:\mathbb{R}^n \rightarrow \mathbb{R}$ are locally Lipschitz.
        \item $D_i:\mathbb{R}^n \rightrightarrows \mathbb{R}^n$ is a locally bounded conservative field with nonempty convex values for $f_i$ for all $i \in \llbracket 1, N \rrbracket$.
        \item $g:\mathbb{R}^n \rightarrow \overline{\mathbb{R}}$ is proper, convex, and locally Lipschitz in its closed domain.
        \item $\Phi:= f+g$ is lower bounded, where $f:=f_1+\cdots+f_N$, and $D_\Phi := \sum\limits_{i = 1}^N D_i + \partial g$.
    \end{enumerate}
\end{assumption}
Most of the conditions above are mild and standard in the literature, as discussed in the literature review. The convex value requirement for $D_i$ can be relaxed, by substituting $D_i$ with its convex hull $\co D_i$. Since $g$ is locally Lipschitz in its closed domain, it holds that $\dom \Phi = \dom g = \dom \partial g$ and $\partial g$ has closed graph (\citet[Theorem 24.4]{rockafellar1970convex}). Thus, $D_\Phi$ has closed graph (\cref{lemma:closed}). Also, thanks to the outer sum rule of conservative fields (\citet[Corollary 4]{bolte2021conservative}), $D_\Phi$ satisfies a chain rule along $D_\Phi$-trajectories (Lemma \ref{lem:chain}). If both $\Phi$ and $D_\Phi$ are additionally semialgebraic, then $D_\Phi$ possesses properties analogous to conservative fields of real-valued functions \citet{bolte2021conservative,josz2023global}; see the last author's PhD thesis \cite[Chapter 5.4]{li2025bounded}. 

\begin{remark}
A common strategy in the literature (see, e.g., \citet{majewski2018analysis}, \citet{davis2020stochastic}, \citet{schechtman2023stochastic}) for handling unbounded set-valued mappings is to explicitly decompose \(g\) as \(g = \tilde{g} + \delta_C\), where \(\tilde{g}\) is a locally Lipschitz convex function and \(\delta_C\) is the indicator function of a closed convex set. However, this approach implicitly assumes that \(g\) admits a globally defined convex extension, which is not always guaranteed. For instance, consider  
\begin{equation*}
    g(x,y):= \begin{cases}
        -2\sqrt{xy} & \text{if } xy \geqslant 1 \text{ and } x > 0, \\
        +\infty & \text{elsewhere}. 
    \end{cases}
\end{equation*}
It is clear that \(\dom(g)\) is a closed convex set and that \(g\) is smooth on its domain. Nevertheless, as shown in \citet{schulz1979finite}, \(g\) does not possess a convex extension over \(\mathbb{R}^2\). This example illustrates the limitations of the decomposition approach and motivates our choice to adopt a more general framework that directly handles convex functions that are locally Lipschitz on their domains. 
\end{remark}

The following assumption is needed for \cref{thm:weakly_convex}, where we assume access to Clarke subdifferentials of $f_i$. Recall that a locally Lipschitz function $f:\mathbb{R}^n\rightarrow\mathbb{R}$ is subdifferentially regular (\citet[2.3.4 Definition]{clarke1990}) if the classical directional derivative exists and agrees with the generalized directional derivative, that is to say, we have
\begin{equation*}
    \lim_{
    t \searrow 0
    } \frac{f(x+th)-f(x)}{t}
     ~ = ~ 
     \limsup_{\text{\scriptsize$\begin{array}{c} y\rightarrow x \\
    t \searrow 0
    \end{array}$}
    } \frac{f(y+th)-f(y)}{t} \vspace*{-1mm}
\end{equation*}
for all $x,h\in \mathbb{R}^n$.
\begin{assumption} \text{ } 
    \label{assumption_weakly_convex}
    \begin{enumerate}
        \item $f_1,\ldots, f_N$ are subdifferentially regular.
        \item $D_i := \partial f_i$ for all $i \in \llbracket 1, N \rrbracket$.
        \item $\Phi$ is weakly convex.
    \end{enumerate}
\end{assumption}
When both \cref{assumption_standing,assumption_weakly_convex} hold, we have that $D_\Phi = \sum_{i = 1}^N \partial f_i + \partial g = \partial \Phi$ from \citet[10.9 Corollary and p. 337]{rockafellar2009variational}. We remark that many functions arising in practice are weakly convex (see, e.g., \citet{chandrasekaran2011rank,abbe2014decoding,davis2020nonsmooth}). The next assumption will be used to establish convergence to stationary points of $\Phi$ (\cref{thm:sgd_finite_len}).
\begin{assumption} \text{ } 
    \label{assumption_smooth} 
    \begin{enumerate}
        \item $f_1,\hdots,f_N$ 
        are differentiable with locally Lipschitz continuous gradients.
        \item $D_i:=\{\nabla f_i\}$ for all $i \in \llbracket 1, N \rrbracket$.
    \end{enumerate}
\end{assumption}

\subsection{Theorems}\label{subsec:thm}

We are now ready to state our main theorems concerning iterates generated by proximal random reshuffling (\cref{alg:prr}). The four theorems in this subsection treat objective functions with varying degrees of regularity, where different step size strategies are needed. None of the theorems assume global Lipschitz continuity or require prior knowledge on the iterates. The results are new even when $g=0$, in which case \cref{alg:prr} reduces to random reshuffling (\citet{bertsekas2011incremental,mishchenko2020random,gurbuzbalaban2021random,pauwels2021incremental}). We discuss some applications of these results in \cref{subsec:example}. The proofs of the theorems are deferred to \cref{sec:proofs}. 

We begin with a theorem that applies to the most general setting in this work, where \cref{alg:prr} is guaranteed to reach a near approximate stationary point.
\begin{theorem}
\label{thm:exists_small}
    Let \cref{assumption_standing} hold and $\delta,\epsilon>0$. For any bounded set $X_0 \subset \dom \Phi$, there exists $\bar{\alpha}>0$ such that for any sequence generated by \cref{alg:prr} initialized in $X_0$ such that \vspace*{-2mm}
    \begin{equation*}
         \alpha_k \in (0,\bar{\alpha}] ~~~ \text{and} ~~~ \sum_{k = 0}^\infty \alpha_k = \infty, \vspace*{-1mm}
    \end{equation*}
    at least one iterate is $(\epsilon,\delta)$-near approximate $D_\Phi$-stationary.
\end{theorem}
In the above setting, we require the step sizes to be nonsummable, which is standard for first-order methods (\citet{robbins1951stochastic,polyak1969minimization}). We also need the step sizes to be small, otherwise the iterates may diverge and never reach a near approximate stationary point (consider $\Phi(x):= x^4$). Examples of such step sizes include those that are eventually constant or diminishing at certain rates (e.g., $\alpha_k = \alpha/(k+1)^p$, $\alpha\in(0,\bar{\alpha}]$, $p\in (0,1]$). 

The conclusion of \cref{thm:exists_small} may appear weak at first glance, but as discussed in our literature review, the PRR algorithm has not been studied under such general assumptions. The closest related result we are aware of is \citet[Theorem 6.2]{davis2020stochastic}, which applies to the stochastic proximal subgradient method with replacement. Unlike the unbiased subgradient estimates assumed in \citet[p. 141]{davis2020stochastic}, PRR relies on biased estimates. This key complication is resolved by \cref{lemma:stoc_tracking}. Moreover, \cref{thm:exists_small} only requires nonsummable step sizes, thus accommodating constant step sizes. In this regime, the result offers new insights even when iterates are bounded. The closest related work we know of considers the random reshuffling algorithm without constraints and requires uniform boundedness of iterates for all sufficiently small step sizes (\citet[Corollary 1, Remark 2]{josz2024global}). However, this condition is difficult to verify and is known to hold only in special cases such as coercive objectives. As illustrated by the example in \cref{fig:constant}, for robust matrix completion, iterates remain bounded but over increasingly large regions as step sizes decrease. Nonetheless, \cref{thm:exists_small} remains applicable in such cases.

In the remaining theorems of this subsection, we address the convergence of \cref{alg:prr}. By additionally assuming that $\Phi$ is weakly convex, we establish convergence of the iterates to an $(\epsilon,\delta)$-near approximate stationary point in the following theorem.

\begin{theorem}\label{thm:weakly_convex}
   Let \cref{assumption_standing,assumption_weakly_convex} hold and $\delta,\epsilon>0$. For all $x_0 \in \dom \Phi$, there exists $T_0\geqslant 0$ such that for all $T\geqslant T_0$, there exists $\bar{\alpha}>0$ such that any sequence generated by \cref{alg:prr} such that \vspace*{-2mm}
   \begin{equation*}
        \alpha_k \in (0,\bar{\alpha}] ~~~ \text{and}~~~ \sum_{k = 0}^\infty \alpha_k = T \vspace*{-1mm}
   \end{equation*}
   converges to an $(\epsilon,\delta)$-near approximate $\partial \Phi$-stationary point.
\end{theorem}
In this theorem, weak convexity plays a crucial role by ensuring the uniqueness of the subgradient trajectory and that its velocity eventually becomes arbitrarily small (see \cref{lemma:essential_conv}). These two properties, together with the tracking lemma and the proposed step sizes, guarantee that the iterates converge to a $(\epsilon,\delta)$-NAS. 

We next elaborate on the choice of step sizes. Given $\bar{\alpha},T>0$, one option is to run the algorithm with constant step size and then stop after finitely many iterations. More precisely, if we let $K\in \mathbb{N}$ be such that $T/K \in (0,\bar{\alpha}]$, then we may implement \cref{alg:prr} with $\alpha_k = T/K$ for all $k \in \llbracket 0, K-1 \rrbracket$ and then terminate. In this case, the $K$th iterate is guaranteed to be an $(\epsilon,\delta)$-near approximate $\partial \Phi$-stationary point. Another option is to use a sequence of geometrically decaying step sizes, i.e., $\alpha_k:= \alpha \rho^k$ for some $\alpha>0$ and $\rho \in (0,1)$, which are known to speed up the local convergence of the (projected) subgradient method (\citet{goffin1977convergence,davis2018subgradient,li2020nonconvex,davis2024stochastic}), if a certain sharpness condition holds around the set of global minima. In order for such step sizes to satisfy the criterion in the above theorem, it suffices to take $\alpha \in (0,\bar{\alpha}]$ and $\rho = 1- \alpha/T$. In order to manage the overall magnitude of the sum in applications such as the nonnegative $\ell_1$ matrix completion problem in \Cref{ex:nrpca}, our theorem suggests increasing $T$ and potentially reducing $\bar{\alpha}$. The theorem ensures that by doing so one reaches regions of interest in the state space.

In \cref{thm:weakly_convex}, we show the convergence to an $(\epsilon,\delta)$-near approximate stationary point, which is not necessarily close to a stationary point of $\Phi$ (i.e., a $(0,0)$-near approximate $\partial \Phi$-stationary point). In fact, simple examples (e.g., $\Phi:= \exp$) can satisfy all assumptions in \cref{thm:weakly_convex} yet fail to admit any stationary points. In the following theorem, we prove convergence to somewhere near a stationary point, if an assumption on $D_\Phi$-trajectories holds.

\begin{theorem}
    \label{thm:nonsmooth_bounded_flow}
   Let \cref{assumption_standing} hold, $\epsilon>0$, and $x_0 \in \dom \Phi$. Assume that there exists a unique $D_\Phi$-trajectory initialized at $x_0$ and that it converges. There exists $T_0\geqslant 0$ such that for all $T\geqslant T_0$, there exists $\bar{\alpha}>0$ such that any sequence generated by \cref{alg:prr} initialized at $x_0$ with \vspace*{-2mm}
   \begin{equation*}
      \alpha_k \in (0,\bar{\alpha}] ~~~\text{and} ~~~ \sum_{k = 0}^\infty \alpha_k = T  \vspace*{-1mm} 
   \end{equation*}
   converges to an $(\epsilon,0)$-near approximate $D_\Phi$-stationary point.
\end{theorem}
We remark that $D_\Phi$-stationary points (i.e., $(0,0)$-near approximate $D_\Phi$-stationary points) exist in the setting of the above theorem, as bounded $D_\Phi$-trajectories converge to $D_\Phi$-stationary points. We use summable step sizes in both \cref{thm:weakly_convex,thm:nonsmooth_bounded_flow} because \cref{alg:prr} with nonsummable step sizes might not converge in either setting. Nonconvergence can easily be verified for $\Phi(x,y):= |x|/y + \delta_C(x,y)$ where $\delta_C$ is the indicator function of the set $C:= \{(x,y) \in \mathbb{R}^2:y\geqslant 1\}$. 

In our final convergence theorem, we provide sufficient conditions for \cref{alg:prr} to converge to stationary points of $\Phi$ with certain nonsummable step sizes. While the regularity assumption is stronger than in the previous theorems, this theorem is connected to them due to the use of the same tracking lemma. Besides differentiability of $f_1,\ldots,f_N$ (\cref{assumption_smooth}), we also require the objective function to be semialgebraic and have bounded subgradient trajectories.

A subset $A$ of $\mathbb{R}^n$ is semialgebraic \cite{bochnak2013real} if it is a finite union of basic semialgebraic sets, which are of the form 
\begin{equation*}
    \{ x \in \mathbb{R}^n : f_1(x) > 0 , \hdots , f_p(x) > 0, f_{p+1}(x) = 0, \hdots , f_q(x) = 0\}
\end{equation*}
where $f_1,\hdots,f_q$ are polynomials with real coefficients. A function $\Phi:\mathbb{R}^n\rightarrow\overline{\mathbb{R}}$ is semialgebraic if its epigraph (or equivalently its graph) is semialgebraic. Similarly, a set-valued mapping $D:\mathbb{R}^n\rightrightarrows \mathbb{R}^m$ is semialgebraic if its graph is semialgebraic. All the results in this manuscript hold more generally in an o-minimal structure on the real field \cite{van1998tame}.

We are now ready to state the last convergence theorem in this paper. Let $\mathbb{N}^* := \{1,2,\hdots\}$. 

\begin{theorem}
    \label{thm:sgd_finite_len}
    Let \cref{assumption_standing,assumption_smooth} hold and $\beta = 0$ if $N=1$, and  $\beta \in (1/2,1)$ otherwise. If $\Phi$ is semialgebraic and has bounded subgradient trajectories, then for any bounded set $X_0\subset\dom \Phi$, there exist $\bar{\alpha},c>0$ such that any sequence generated by \cref{alg:prr} initialized in $X_0$ with \vspace*{-1mm}
    \begin{equation}\label{eq:prr_smooth_steps}
        \alpha \in (0,\bar{\alpha}] ~~~\text{and}~~~  \alpha_k = \frac{\alpha}{(k+1)^\beta}, ~~~ \forall k \in \mathbb{N}, \vspace*{-2mm}
    \end{equation}
    converges to a $(0,0)$-near approximate $\partial \Phi$-stationary point, \vspace*{-2mm}
    \begin{equation*}
         \sum_{i=0}^\infty\|x_{i+1}-x_i\|\leqslant c, ~~~\text{and}~~~ \min_{i\in \llbracket 0,k \rrbracket} d(0,\partial \Phi(x_i)) = o(k^{\beta-1}), ~~~\forall k \in \mathbb{N}^*.
    \end{equation*}
\end{theorem}

One might expect a two-step route: combine the tracking lemma (\cref{lemma:stoc_tracking}) with bounded subgradient trajectories to first deduce bounded iterates, and then invoke standard arguments to obtain convergence. In our setting, however, we are not aware of a clean way to prove boundedness of a single discrete trajectory as a separate preliminary step. Instead, the proof in \cref{app:thm4proof} proceeds by establishing the stronger statement $\Gamma(X_0,\bar{\alpha})<\infty$, namely a uniform finite-length bound that holds for all initial points $x_0\in X_0$ and all step sizes $\alpha\le \bar{\alpha}$. Importantly, this uniform bound is obtained by using all ingredients together---the tracking lemma in \cref{lemma:stoc_tracking}, the length formula in \cref{prop:discrete_length_formula}, and the recursive estimate \cref{eq:main_recur}. Once $\Gamma(X_0,\bar{\alpha})<\infty$ is established, boundedness of iterates follows immediately as a direct corollary. Consequently, restating \cref{thm:sgd_finite_len} under an a priori bounded-iterates assumption would obscure the actual structure of the proof.

\begin{remark}\label{rem:cstEstGen}
    In practice, for \cref{thm:exists_small,thm:sgd_finite_len}, the constant $\bar{\alpha}$ needs to be estimated by tuning step sizes, typically by decreasing them or using various line search methods. For \cref{thm:weakly_convex,thm:nonsmooth_bounded_flow}, both $\bar{\alpha}$ and $T_0$ require estimation, implying adjustments of step sizes along with an increase in the number of iterations. In general, these constants depend on how accurately the iterates approximate the subgradient trajectory described in \cref{lemma:stoc_tracking}, making their precise determination challenging. Nevertheless, the practical tuning of step sizes and iteration counts remains feasible, and the theorems ensure convergence or desired outcomes when appropriate parameters are chosen. 
\end{remark}

\begin{remark}
    Inspired by the practical stopping criteria proposed in \citet{yu2023high}, we suggest using the following criterion:
    \begin{equation}
        s_k:=\frac{\|x_{k+1}-x_k\|}{\alpha_k} \leqslant \eta \delta,
    \end{equation}
    where $\eta>0$ is a tolerance parameter chosen by the user. This criterion provides a practical measure of stationarity, ensuring algorithm termination when close to a stationary point under the assumptions of \cref{thm:sgd_finite_len}. Defining the error term
    \begin{equation}
        e_k:=\nabla f(x_k) - \sum_{i=1}^N\nabla f_{\sigma_i^k}(x_{k,i-1}), 
    \end{equation}
    it follows from \cref{cor:Oalpha} that $s_k \leqslant  2d(0,\partial\Phi(x_k))+2\|e_k\|$. Given that $\|e_k\|=O(\alpha_k)$ by \eqref{eq:approxdes_inner2}, the stopping criterion is guaranteed to be triggered once a $\delta$-stationary point ($(0,\delta)$-NAS) is reached, provided sufficiently small step sizes and sufficiently large $\eta$ are used. Moreover, we also obtain $d(0,\partial\Phi(x_{k+1})) \leqslant 2s_k+\|e_k\|$, ensuring that the triggered stopping criterion indeed identifies a $\delta$-stationary point. Unfortunately, this criterion has no guarantee under the settings of \cref{thm:exists_small,thm:weakly_convex,thm:nonsmooth_bounded_flow}.
\end{remark}

\subsection{Examples}\label{subsec:example}
The four theorems can respectively be applied to four examples with increasing degrees of regularity. In the first example, we show that \cref{thm:exists_small} can be applied in deep learning (see, e.g., \citet{lecun2015deep}).
\begin{example}[training of neural networks]\label{ex:nn}
In a typical setting of supervised learning using neural networks, we are given a training set $\{(x_1,y_1),\ldots, (x_N,y_N)\}$ where $x_i \in \mathbb{R}^p$ and $y_i\in \mathbb{R}^q$ are the feature and label for the $i$th sample, respectively. The goal is to find a neural network $h(\cdot,\theta):\mathbb{R}^p \rightarrow \mathbb{R}^q$ by solving
\begin{equation*}
    \inf_{\theta \in \mathbb{R}^n} \ell(h(x_1,\theta),y_1) + \cdots + \ell(h(x_N,\theta),y_N). 
\end{equation*}
In most circumstances, $h$ is a composition of affine functions and nonlinear activation functions parametrized by the weights $\theta\in \mathbb{R}^n$. The loss function $\ell:\mathbb{R}^q\times \mathbb{R}^q \rightarrow \mathbb{R}$ (e.g., mean squared loss or cross entropy loss) is always lower bounded. Let $f_i := \ell(h(x_i,\cdot),y_i)$, which might not be differentiable everywhere. In practice, ``derivatives'' of $f_i$ are obtained via backpropagation, whose outputs constitute a conservative field $D_i$ for $f_i$ (\citet[Theorem 8]{bolte2021conservative}).
\end{example}

\cref{thm:exists_small} can readily be applied to provide a guarantee for the implementation of the stochastic subgradient method in practice (\citet{paszke2019pytorch}). Indeed, \cref{assumption_standing} is verified by letting $g:= 0$ and \cref{alg:prr} reduces to random reshuffling. We would like to highlight that we do not make any assumption on coercivity of the objective function or boundedness of the iterates, in contrast to \citet{josz2024global}. The theorem can also be applied to handle nonsmooth regularizers such as the group sparse/$\ell_1$ regularizers (\citet{scardapane2017group,yuan2006model}). 

In the next example, we see that if the objective function is weakly convex, then \cref{thm:weakly_convex} can be applied to guarantee the asymptotic behavior of \cref{alg:prr}.
\begin{example}[nonnegative $\ell_1$ matrix completion]\label{ex:nrpca}
Let $M\in \mathbb{R}^{m\times n}$ and $\Omega \subset \llbracket 1,m \rrbracket\times \llbracket 1,n \rrbracket$ be a collection of observed entries. We seek to
solve
\begin{equation*}
    \inf_{X,Y\geqslant 0} ~~~ \sum_{(i,j)\in \Omega} |(XY-M)_{ij}|
\end{equation*}
where $X \in \mathbb{R}^{m \times r}$ and $Y \in \mathbb{R}^{r \times n}$. This is a nonconvex formulation of nonnegative robust principal component analysis with partial observations (\citet{candes2011robust}), with the rank-one case studied in \citet{fattahi2020exact}. Let $f_1,\hdots,f_N$ denote the summands in the objective and $g$ denote the indicator of $\mathbb{R}_+^{m\times r} \times \mathbb{R}_+^{r\times n}$. Each function $f_i$ is weakly convex according to \citet[Lemma 4.2]{drusvyatskiy2019efficiency}, hence subdifferentially regular (\citet[2.5.6 Proposition]{clarke1990}). The sum $\Phi := f_1+\cdots+f_N+g$ is then weakly convex. We may thus apply \cref{thm:weakly_convex} to guarantee the convergence of \cref{alg:prr} to an $(\epsilon,\delta)$-near approximate $\partial \Phi$-stationary point.
\end{example}

Recall that the strengthened guarantees in \cref{thm:nonsmooth_bounded_flow,thm:sgd_finite_len} are due to the fact that $\Phi$ has bounded and convergent subgradient trajectories. We next show two examples where this assumption can be verified.

\begin{example}[$\ell_1$ matrix sensing]\label{ex:rms}
Given sensing matrices $A_1,\ldots,A_N \in \mathbb{R}^{m\times n}$ and measurements $b_1,\ldots, b_N\in \mathbb{R}$, we aim to recover a low-rank matrix (\citet{li2020nonconvex,ma2023global}) by solving
\begin{equation*}
    \inf_{X,Y} ~~ \sum_{i = 1}^N \left| \left\langle A_i, XY\right\rangle - b_i \right| 
\end{equation*}
where $X \in \mathbb{R}^{m \times r}$ and $Y \in \mathbb{R}^{r \times n}$. Assuming the lower bound in the $\ell_1/\ell_2$-restricted isometry property (\citet{chen2015exact,li2020nonconvex}), the objective function has convergent subgradient trajectories, following similar arguments as in the proof of \citet[Proposition 4.4]{josz2023certifying}. Uniqueness of the subgradient trajectories is due to weak convexity of the objective function (\citet{marcellin2006evolution}). This allows us to apply \cref{thm:nonsmooth_bounded_flow} with $D_i$'s being the Clarke subdifferentials of the summands and $g:=0$. Therefore, with the step sizes in \cref{thm:nonsmooth_bounded_flow}, \cref{alg:prr} is guaranteed to converge near a critical point.
\end{example}

\begin{example}[nonnegative $\ell_p$ matrix factorization]\label{ex:ngmf}
Given $M\in \mathbb{R}^{m\times n}$ and $p\geqslant 2$, we aim to solve
\begin{equation*}
    \inf_{X,Y\geqslant 0} ~~ \|XY-M\|_p^p
\end{equation*}
where $(X,Y) \in \mathbb{R}^{m\times r} \times \mathbb{R}^{r\times n}$ and $\|\cdot\|_p$ is the entrywise $\ell_p$-norm. One possible way is to set each $f_i$ to be one of $|\sum_{k=1}^rX_{\ell k}Y_{kj}-M_{\ell j}|^p$ for every $\ell=1,\ldots,m,\,j=1,\ldots,n$ and $g$ to be the indicator of $\mathbb{R}_+^{m\times r} \times \mathbb{R}_+^{r\times n}$. The summands $f_i$ and the function $g$ readily satisfy \cref{assumption_standing,assumption_smooth}. According to \citet[Proposition 2.5]{li2025bounded}, the objective function has bounded subgradient trajectories. Therefore, \cref{alg:prr} is guaranteed to converge to a stationary point by \cref{thm:sgd_finite_len}. When $N=1$, convergence happens at the rate $o(1/k)$. When $N>1$, convergence happens at the rate $o(1/k^{1/2-\epsilon})$ if $\beta = 1/2+\epsilon$ for any $\epsilon \in (0,1/2)$. We actually proved boundedness of subgradient trajectories for any $p\geqslant 1$, so that \cref{thm:nonsmooth_bounded_flow} can be applied when $p\in [1,2)$. 
\end{example}

Various algorithms have been proposed for solving nonnegative $\ell_2$ matrix factorization over the past three decades (\citet{gillis2020nonnegative}). We briefly summarize their convergence guarantees: 
\begin{enumerate}
    \item convergence in function value: multiplicative update (\citet{lee2000algorithms}), hierarchical alternating least squares (\citet{cichocki2007hierarchical,kim2014algorithms}) (if the columns of the factors remain nonzero);
    \item stationarity of limit points: active-set method (\citet{kim2008nonnegative}), proximal gradient method (\citet{lin2007projected}), alternating direction method of multipliers (\citet{hajinezhad2016nonnegative});
    \item convergence of bounded iterates to stationary points: proximal alternating linearized minimization (\citet{bolte2014proximal}), Bregman proximal gradient method (\citet{mukkamala2019beyond,teboulle2020novel}).
\end{enumerate}
By modifying the above algorithms, one can obtain convergence guarantees to modified versions of NMF. For example, $X,Y\geqslant 0$ can be replaced by $X,Y\geqslant \epsilon$ for some small $\epsilon>0$, or by $u \geqslant X,Y\geqslant 0$ for some large $u>0$. A regularizer can also be added to the objective. In this vein, modified multiplicative update (\citet{takahashi2014global}) and modified hierarchical alternating least squares (\citet{kimura2015global}) yield bounded iterates and subsequential convergence to modified problems. The projected gradient method (\citet{calamai1987projected,lin2007projected}) with line search and box constraints produces bounded iterates whose limit points are stationary points of the modified problem. A similar result holds with norm-based regularizers (\citet{rakotomamonjy2013direct}).

Despite the extensive convergence results available for various algorithms applied to NMF, to the best of our knowledge, no existing work establishes global convergence for finding stationary points of NMF without assuming the existence of limit points, even with the projected gradient method (adaptive step sizes or not). This is because global convergence of the forward-backward algorithm typically requires both bounded iterates and the sufficient decrease property (H1 in \citet{attouch2013convergence}). However, under only local gradient Lipschitz continuity and a constant step size, H1 may fail even when iterates are bounded, as illustrated in the examples from the fourth paragraph of \cref{sec:Introduction}. To address this, some authors (e.g., \citet{teboulle2020novel}) introduce Bregman divergences that ensure global relative smoothness, allowing the use of the H1-H3 framework with boundedness and constant step size $\alpha_k = 1/\max\{(1+\|M\|_F^2)/2,2\}$ to establish convergence.

\begin{remark}\label{rem:NMFstep}
    For \cref{ex:ngmf}, our theory establishes convergence guarantees for \cref{alg:prr} under a constant step size ($N=1$ and $\beta=0$) and a diminishing step-size schedule ($N>1$ and $\beta\in(1/2,1)$). While line-search rules are not covered by our analysis, in practice, backtracking line search is a standard and effective way to tune step sizes, as it can automatically enforce a sufficient decrease condition for NMF. 
\end{remark}

\section{Proofs}\label{sec:proofs}
In this section, we prove the results in \cref{sec:main results}.

\subsection{Tracking lemma}\label{subsec:proof_of_stoc_tracking}

The object of this subsection is to show that the iterates generated by \cref{alg:prr} can be tracked by $D_\Phi$-trajectories. This result, stated below, is one of the key technical contributions of this work. All of the theorems pertaining to \cref{alg:prr} rely on it.
\begin{lemma}
\label{lemma:stoc_tracking}
     Let \cref{assumption_standing} hold. For any compact set $X_0 \subset \dom \Phi$ and $\epsilon,T>0$, there exists $\bar{\alpha}>0$ such that for any sequence generated by \cref{alg:prr} initialized in $X_0$ with $\alpha_0,\alpha_1,\hdots \in (0, \bar{\alpha}]$, there exists a solution $x:[0,T]\rightarrow \dom \Phi$ to $x' \in -D_\Phi(x)$ such that $x(0) \in X_0$ and
\begin{equation*}
   \forall k\in \mathbb{N}^*,~~~ \alpha_0+\cdots+\alpha_{k-1} \leqslant T ~~~\Longrightarrow~~~ \|x_k - x(\alpha_0+\cdots+\alpha_{k-1})\| \leqslant \epsilon.
\end{equation*}
\end{lemma}
We clarify the contribution of \cref{lemma:stoc_tracking} in detail. In recent years, there has been significant progress in analyzing optimization algorithms through their continuous counterparts (see, e.g., \citet{duchi2018stochastic,davis2020stochastic,bolte2021conservative,salim2018random,pauwels2021incremental,bianchi2022convergence,josz2023lyapunov,josz2023global,josz2024global,josz2023convergence}), while this idea traces back to the stochastic approximation literature in the 1970s (\citet{ljung1977analysis,kushner1977general,kushner1977convergence}). To the best of our knowledge, \cref{lemma:stoc_tracking} is the first result to demonstrate that such approximations remain valid over any finite time period, even without the assumption of a smooth objective function (\citet[Proposition 4]{josz2023global}) or bounded iterates (\citet[Lemma 1]{josz2023lyapunov}). This nontrivial extension is made possible thanks to \cref{prop:bounded}, which asserts uniform boundedness of the iterates produced by \cref{alg:prr} when initialized in a bounded set. In addition, we have to overcome hurdles introduced by the variable step sizes, despite the use of similar techniques as in \citet[Lemma 1]{josz2023lyapunov} and \citet[Proposition 1]{josz2024global} at several places. 

We next provide a roadmap for the remainder of this subsection. We first study convergence towards solutions to differential inclusions with unbounded right-hand side (\cref{prop:conv}). It requires closedness of the graph of the sum of two set-valued mappings (\cref{lemma:closed}). Note that \cref{prop:conv} generalizes \citet[Theorem 1, p. 60]{aubin1984differential} by relaxing the local boundedness requirement, and is used later for the convergence to solutions of $x'\in -D_\Phi(x)$. After that, we turn our attention to the iterates generated by \cref{alg:prr}. We show that the displacement generated by the proximal operator is upper bounded by the distance from $0$ to the subdifferential of $g$ (\cref{lemma:Oalpha}), which is locally bounded (\cref{lemma:min_grad}). It then follows that the distances between consecutive iterates are well controlled by the step sizes (\cref{cor:Oalpha}). Combining these results, we show that iterates are uniformly bounded, given that the sum of step sizes is bounded (\cref{prop:bounded}). Finally, we prove the approximation of iterates by $D_\Phi$-trajectories (\cref{lemma:stoc_tracking}).

We first state some results regarding the sum of two set-valued mappings and the convergence to its solutions.
 
\begin{lemma}\label{lemma:closed}
    If $F:\mathbb{R}^n\rightrightarrows\mathbb{R}^n$ is locally bounded with closed graph and $G:\mathbb{R}^n\rightrightarrows\mathbb{R}^n$ has closed graph, then $F+G$ has closed graph. 
\end{lemma}
\begin{proof}
    Suppose $\gra (F+G) \ni(x_k,y_k) \rightarrow (x,y)$, i.e., $y_k = u_k + v_k$, for some $u_k \in F(x_k)$, $v_k \in  G(x_k)$. Since $F$ is locally bounded, $u_k$ admits a subsequence which converges to $u\in \mathbb{R}^n$. As $F$ has closed graph, $u\in F(x)$. Then $G(x_k) \ni v_k = y_k - u_k \rightarrow  y - u$. Since $G$ has closed graph, we have $G(x) \ni y-u$, i.e., $y \in u + G(x) \subset F(x)+G(x)$.
\end{proof}

\begin{proposition}
\label{prop:conv}
    Let $F,G:\mathbb{R}^n\rightrightarrows\mathbb{R}^n$ have closed graphs with convex values and $F$ be locally bounded. Assume that $F+G$ is proper. Let $I\subset\mathbb{R}$ be a bounded interval and $x_k,y_k:I\rightarrow \mathbb{R}^n$ be measurable functions such that for almost every $t\in I$ and for any neighborhood $U$ of $0$ in $\mathbb{R}^{2n}$, there exists $k_0\in\mathbb{N}$ such that 
    \begin{equation*}
        (x_k(t),y_k(t)) \in \gra(F+G)+U, \quad \forall k\geqslant k_0. 
    \end{equation*}
    Assume that 
    \begin{enumerate}
        \item $(x_k(\cdot))_{k\in\mathbb{N}}$ is uniformly bounded in $L^\infty(I,\mathbb{R}^n)$ and $x_k(\cdot)\to x(\cdot)$ a.e. on $I$; 
        \item $(y_k(\cdot))_{k\in\mathbb{N}}$ is uniformly bounded in $L^\infty(I,\mathbb{R}^n)$ and $y_k(\cdot)\to y(\cdot)$ weakly in $L^1(I,\mathbb{R}^n)$. 
    \end{enumerate}
    Then $(x(t),y(t))\in \gra(F+G)$ for almost all $t\in I$. 
\end{proposition}

For the proof of \cref{prop:conv}, refer to \cref{app:propconvproof}. The proof follows a similar approach to \citet[Lemma 2]{schechtman2023stochastic}, which truncates the unbounded-valued component of the system to a bounded set-valued mapping (see also \citet{brezis1973operateurs}). However, our result generalizes this idea, extending beyond the case where $F$ is the subdifferential of a locally Lipschitz function and $G$ is a normal cone. To the best of our knowledge, \cref{prop:conv} does not have an exact counterpart in the literature. Consequently, we retain the statement for convenience and ease of citation. 

We now move on to the study of iterates generated by \cref{alg:prr}. We begin with a simple upper bound on the length traveled after one iteration.
\begin{lemma}
\label{lemma:Oalpha}
    Let $g:\mathbb{R}^n\rightarrow\overline{\mathbb{R}}$ be proper and convex. For any $x\in \dom \partial g$, $\alpha>0$, $d\in \mathbb{R}^n$, it holds that
    \begin{equation*}
        \frac{1}{\alpha}\|\prox_{\alpha g}(x-\alpha d) - x\| \leqslant 2\|d\| + 2d(0,\partial g(x)).
    \end{equation*}
\end{lemma}
\begin{proof}
Let $y:= \prox_{\alpha g}(x-\alpha d)$. Applying the definition of the proximal operator, it holds that
    \begin{equation*}
        g(y) + \frac{1}{2\alpha} \|y - (x-\alpha d)\|^2 \leqslant g(x) + \frac{1}{2\alpha} \|x - (x - \alpha d)\|^2.
    \end{equation*}
    Rearranging the above inequality yields that
    \begin{equation}\label{eq:prox-1}
         \frac{1}{2\alpha} \|y - x\|^2 \leqslant g(x) - g(y) +  \langle x - y, d\rangle.
    \end{equation}
    Meanwhile, by convexity of $g$, for any $v\in\partial g(x)$, we have $g(y) \geqslant g(x) + \langle v, y-x\rangle$. Together with \eqref{eq:prox-1}, this yields that $\|y - x\|^2/(2\alpha) \leqslant \langle x - y, d + v\rangle$. By Cauchy-Schwarz and triangular inequalities, we have
     \begin{equation*}
         \frac{1}{2\alpha} \|y - x\|^2 \leqslant \langle x - y, d + v\rangle \leqslant \|y - x\| \|d + v\| \leqslant \|y - x\| (\|d\| + \|v\|).
     \end{equation*}
     We conclude by taking the infimum with respect to $v\in \partial g(x)$ on the right-hand side of the above inequality. 
\end{proof}

The bound on the length after one iteration previously obtained will prove useful once we can control the distance between the origin and $\partial g$. This is the object of the next result, which can be regarded as a simple corollary of the norm preserving extension result of convex Lipschitz functions in \citet{cobzas1978norm}. Let $B(S,r) := S+B(0,r)$. 

\begin{lemma}
\label{lemma:min_grad}
    If $g:\mathbb{R}^n\rightarrow\overline{\mathbb{R}}$ is proper, convex, locally Lipschitz on its domain, and has closed domain, then $d(0,\partial g)$ is locally bounded over $\dom \partial g$ and $\dom \partial g = \dom g$.
\end{lemma}
\begin{proof}
    Let $X\subset\mathbb{R}^n$ be any compact set. Notice that $Y:=\dom g\cap B(\co X,1)$ is a convex compact set, so the restriction of $g$ on $Y$, denoted as $g\big|_Y:Y\rightarrow\mathbb{R}$, is an $L$-Lipschitz convex function on $Y$. Applying \citet[Theorem 1]{cobzas1978norm}, we obtain an $L$-Lipschitz convex function  $\hat{g}:\mathbb{R}^n\rightarrow\mathbb{R}$ such that $\hat{g}(x)=g\big|_Y(x)$ for all $x\in Y$. Thus, $g(x)=\hat{g}(x)+\delta_Y(x)$ for all $x\in B(\co X,1)$ and $\partial g(x)=\partial  \hat{g}(x)+N_Y(x)$ for all $x\in X\cap\dom g$, where $N_Y(x)$ denotes the normal cone of $Y$ at $x$. Since $0\in N_Y(x)$ for all $x\in X\cap \dom g$, we have that $d(0,\partial g(x)) \leqslant d(0,\partial \hat{g}(x))\leqslant L$ for all $x\in X\cap \dom g$. Thus, $d(0,\partial g)$ is locally bounded over $\dom g$. 
\end{proof}

Combining \cref{lemma:Oalpha,lemma:min_grad}, one can control the length traveled after one iteration directly via the step sizes. To see why, it will be convenient to use the notation $\|S\| := \sup \{\|s\|: s \in S\}$ for any $S\subset \mathbb{R}^n$ and $r\geqslant 0$.

\begin{corollary}\label{cor:Oalpha}
Suppose \cref{assumption_standing} holds. Let $X \subset \mathbb{R}^n$ be bounded and $x_0,\ldots,x_{K+1}$ be generated by \cref{alg:prr}. Consider the constants
\begin{equation*}
    L_i:=\sup_{B(X,1)}\| D_i\|, ~~~ L:=\max\{L_1,\ldots,L_N\}, ~~~ L_g:=\sup_{B(X,1)\cap \dom g} ~ d(0,\partial g),
\end{equation*}
and assume that $x_0,\hdots,x_K\in X$. 
\begin{enumerate}
    \item If $\alpha_0,\hdots,\alpha_K \leqslant 1/(NL)$, then $\{x_{k,i}\}^{i\in \llbracket 1, N \rrbracket}_{k\in \llbracket 0, K \rrbracket}\subset B(X,1)$.
    \item If $\alpha_0,\hdots,\alpha_K \leqslant 1/(2NL+2L_g)$, then $\|x_{k+1}-x_k\|\leqslant 2(NL+L_g)\alpha_k$ for all $k \in \llbracket 0,K\rrbracket$ and $x_{K+1}\in B(X,1)$.
\end{enumerate}

\end{corollary}
\begin{proof}
We first note that $L_g<\infty$ by \cref{lemma:min_grad} and a standard compactness argument. The first part can be proved by induction, with hypothesis $H_i: ~ ``\forall k \in \llbracket 0,K \rrbracket, ~ \|x_{k,i-1}-x_k\|\leqslant (i-1)/N$''. By \cref{alg:prr}, $x_{k,0}=x_k\in X$ for all $k \in \llbracket 0,K \rrbracket$, hence $H_1$. If $H_i$ holds, then $x_{k,i-1}\in B(X,1)$ and 
    \begin{align*}
        \|x_{k,i}-x_k\| &\leqslant \|x_{k,i}-x_{k,i-1}\| + \|x_{k,i-1}-x_k\| \leqslant \alpha_k\left\|D_{\sigma_i^k}(x_{k,i-1})\right\| + \frac{i-1}{N} \\
        &\leqslant \frac{1}{NL}\cdot L_{\sigma_i^k} + \frac{i-1}{N} \leqslant \frac{i}{N},
    \end{align*}
    hence $H_{i+1}$ holds. As for the second part, by \cref{lemma:Oalpha} and part one, we have 
    \begin{equation*}
        \|x_{k+1}-x_k\| \leqslant 2\alpha_k\left(\left\|\sum_{i=1}^N D_{\sigma_i^k}(x_{k,i-1})\right\|+d(0,\partial g(x_k))\right) \leqslant 2\alpha_k(NL+L_g) \leqslant 1.
    \end{equation*}
    Finally, since $x_K\in X$, it holds that $x_{K+1}\in B(X,1)$. 
\end{proof}

We next record a property of $D_\Phi$-trajectories.
\begin{lemma}\label{lem:chain}
    Let \cref{assumption_standing} hold and let $I$ be an interval of $\mathbb{R}_+$. If $x:I\rightarrow \mathbb{R}^n$ is a solution to $x'\in -D_\Phi(x)$, then $\Phi\circ x$ is differentiable almost everywhere on $I$ with
        \begin{equation*}
            \left(\Phi\circ x\right)'(t) = -\|x'(t)\|^2 \quad\text{and}\quad \|x'(t)\| = d\left(0,D_\Phi(x(t))\right), \quad \text{for almost every $t\in I$}. 
        \end{equation*}
\end{lemma}

Proof of the lemma is omitted as it follows immediately from the definition of conservative fields and \citet[Proposition 4.10 \& Lemma 4.11]{drusvyatskiy2015curves} (see also \citet{brezis1973operateurs}). We next prove the final result needed for the proof of \cref{lemma:stoc_tracking}. We show that iterates generated by \cref{alg:prr} remain in a bounded set, whose radius only depends on the sum of step sizes.
\begin{proposition}
\label{prop:bounded}
    Let \cref{assumption_standing} hold. Then for any bounded set $X_0 \subset \dom \Phi$ and $T>0$, there exist $\bar{\alpha}, r>0$ such that for any sequence generated by \cref{alg:prr} initialized in $X_0$ with $\alpha_0,\alpha_1,\hdots \leqslant \bar{\alpha}$, we have
\begin{equation*}
   \forall k\in \mathbb{N}^*,~~~ \alpha_0+\cdots+\alpha_{k-1} \leqslant T ~~~\Longrightarrow~~~ x_0, \ldots, x_k \in B(0,r).
\end{equation*}
\end{proposition}
\begin{proof}
    We assume without loss of generality that $X_0 \subset \dom \Phi$ is nonempty and compact. Fix $T>0$. There exists $Q>0$ such that for any absolutely continuous $x:[0,T] \rightarrow \mathbb{R}^n$  that satisfies $x(0) \in X_0$ and is a solution to 
    \begin{equation}
    \label{eq:sum_di}
    x' \in -D_\Phi(x) =  -\sum_{i = 1}^N D_i (x) - \partial g(x),
\end{equation}
it holds that $d(x(\tilde{T}),X_0) \leqslant \|x(\tilde{T})-x(0)\| \leqslant \int_0^{\tilde{T}} \|x'(t)\|~dt \leqslant \int_0^{T} \|x'(t)\|~dt \leqslant Q$ for all $\tilde{T} \in [0,T]$. Indeed, let $x(\cdot)$ be any such function, we have
\begin{subequations}
\label{eq:cont_len}
    \begin{align}
        \int_{0}^T \|x'(t)\|~dt &\leqslant \sqrt{T} \sqrt{\int_{0}^T \|x'(t)\|^2~dt}\label{eq:cont_len_a} \\
        &= \sqrt{T} \sqrt{\int_{0}^T -(\Phi\circ x)'(t)~dt} \label{eq:cont_len_b}\\
        &= \sqrt{T} \sqrt{\Phi(x(0)) - \Phi(x(T))} \label{eq:cont_len_c}\\
        &\leqslant \sqrt{T} \sqrt{\sup_{x_0\in X_0}\Phi(x_0) - \inf_{y\in \mathbb{R}^n}\Phi(y)}=: Q. \label{eq:cont_len_d}
    \end{align}
\end{subequations}
Above, \eqref{eq:cont_len_a} follows from the Cauchy-Schwarz inequality and \eqref{eq:cont_len_b} follows from \cref{lem:chain}. Let $R>0$ such that $X_0 \subset B(0,R)$ and $L>0$ such that $D_i(B(0,R+Q+2)) \subset B(0,L)$ for all $i \in \llbracket 1,N \rrbracket$. By \cref{lemma:min_grad} and a standard compactness argument, there exists $L_g>0$ such that $d(0,\partial g(x))\leqslant L_g$ for all $x\in B(0,R+Q+2)\cap \dom \Phi$.

We next reason by contradiction and assume that for any $r>0$, there exist a positive sequence $\bar{\alpha}_m \rightarrow 0$, a sequence $(K_m)_{m\in \mathbb{N}}$ of natural numbers, and sequences of iterates  $(x_k^m)_{k\in \mathbb{N}}$ generated by \cref{alg:prr} with step sizes $\alpha_0^m, \alpha_1^m, \hdots \leqslant \bar{\alpha}_m$, $\sum_{k = 0}^{K_m - 1}\alpha_k^m \leqslant T$, and $x_0^m \in X_0$ such that $\max\{\|x_k^{m}\|:k = 0, \ldots,K_m\}>r$ for any $m\in \mathbb{N}$. Take $r = R+Q+2$ and assume that $\bar{\alpha}_m \leqslant 1/(2NL+2L_g)$ without loss of generality. For each $m\in \mathbb{N}$, let $k_m := \min\{k\in \mathbb{N}:x_{k}^m \in B(0,r), x_{k+1}^m \not\in B(0,r)\}$. We have that $k_m \leqslant K_m - 1$ following our assumption.  Thus $\sum_{k = 0}^{k_m - 1} \alpha_k^m \leqslant T$ for any $m \in \mathbb{N}$ and $\bar{T}:= \liminf_{m \rightarrow \infty}\sum_{k = 0}^{k_m - 1}\alpha_k^m \in [0, T]$. By taking a subsequence if necessary, assume that $\lim_{m \rightarrow \infty}\sum_{k = 0}^{k_m - 1}\alpha_k^m = \bar{T}$.

Let $T^m_0 := 0$ and $T^m_k := \sum_{i = 0}^{k-1}\alpha_i^m$ for all $k\in \mathbb{N}^*$ and $m \in \mathbb{N}$. For each sequence $x_0^m, x_1^m, \ldots, x_{k_m}^m$, consider the (extended) linear interpolation $\bar{x}^m:[0,\max\{T^m_{k_m},\bar{T}\}] \rightarrow \mathbb{R}^n$ defined by
	\begin{equation*}
		\bar{x}^m(t) = x_k^m + (t - T_k^m) \frac{x_{k+1}^m - x_{k}^m}{\alpha_k^m}
	\end{equation*}
	for any $t \in [T_k^m, T_{k+1}^m]$ and $k \in\{ 0,1, \ldots, k_m - 1\}$. Also, $\bar{x}^m(t) = x_{k_m}^m$ for $t \in [T^m_{k_m}, \bar{T}]$ if $T^m_{k_m} < \bar{T}$. As $x_k^m \in B(0,r)$ for $k = 0,\ldots, k_m$, we know that $\bar{x}^m(t) \in B(0,r)$ for all $t\in [0,\max\{T^m_{k_m},\bar{T}\}]$ by the convexity of $B(0,r)$. For any $t \in (T_k^m, T_{k+1}^m)$ and $k \in\{ 0,1, \ldots k_m - 1\}$, it holds that $(\bar{x}^m)'(t) = (x_{k+1}^m - x_{k}^m)/\alpha_k^m \in B(0,2(NL+L_g))$ by \cref{cor:Oalpha}. Also, $(\bar{x}^m)'(t) = 0$ for any $t \in (T_{k_m}^m, \bar{T})$ if $T_{k_m}^m < \bar{T}$. By successively applying the Arzel\`a-Ascoli and the Banach-Alaoglu theorems (see \citet[Theorem 4 p. 13]{aubin1984differential}), there exist a subsequence (again denoted $(\bar{x}^m(\cdot))_{m\in \mathbb{N}}$) and an absolutely continuous function $x:[0,\bar{T}]\rightarrow \mathbb{R}^n$ such that $\bar{x}^{m}\big|_{[0,\bar{T}]}(\cdot)$ converges uniformly to $x(\cdot)$ and $(\bar{x}^{m}\big|_{[0,\bar{T}]})'(\cdot)$ converges weakly to $x'(\cdot)$ in $L^1([0,\bar{T}],\mathbb{R}^n)$. In addition, for almost every $t \in (0,\bar{T})$, since $T_{k_m}^m \rightarrow \bar{T}$, for sufficiently large $m$, it holds that $t\in (T_{k}^m, T_{k+1}^m)$ for some $k \in\llbracket 0, k_m - 1\rrbracket$. Therefore, for any neighborhood $U$ of $0$ and for all sufficiently large $m$, it holds that
\begin{subequations}
\label{eq:nb}
    \begin{align}
        (\bar{x}^{m})'(t) &= \frac{x_{k+1}^m - x_{k}^m}{\alpha_k^m} \label{eq:nb-a}\\
        &= \sum_{i = 1}^N \frac{x_{k,i}^m - x_{k,i-1}^m}{\alpha_k^m} + \frac{x_{k+1}^m - x_{k,N}^m}{\alpha_k^m} \label{eq:nb-b}\\
        &\in -\sum_{i = 1}^N D_{\sigma_i^k}(x^m_{k,i-1}) - \left(\partial g(x_{k+1}^m)\cap B(0,L_2)\right) \label{eq:nb-c}\\
        &\subset -\sum_{i = 1}^N D_i(x(t)) - \partial g(x(t)) + U, \label{eq:nb-d}
    \end{align}
\end{subequations}
where $L_2:=3NL+2L_g$. Above, \eqref{eq:nb-c} follows from \citet[Theorem 6.39]{beck2017first} and \cref{lemma:Oalpha}. \eqref{eq:nb-d} is due to the fact that $x_{k,i}^m \rightarrow x(t)$ for all $i \in \llbracket 0,N-1 \rrbracket$, $x_{k+1}^m \rightarrow x(t)$ as $m\rightarrow \infty$, and that $D_i$ and $\partial g \cap B(0,L_2)$ are locally bounded with closed graphs by \cref{assumption_standing}.

By \cref{lemma:closed}, $D_1+\cdots+D_N$ is locally bounded with closed graph. According to \cref{prop:conv}, $x(\cdot)$ satisfies \eqref{eq:sum_di} for almost every $t\in (0,\bar{T})$. Also, $x(0) = \lim_{m\rightarrow \infty}\bar{x}^m(0) \in X_0$. Recall that $x(\bar{T}) \in B(X_0,Q) \subset B(0,R+Q)$. Notice that $\lim_{m\rightarrow \infty} \bar{x}^m(\bar{T}) = x(\bar{T}) \in B(0,R+Q)$ and $\|\bar{x}^m(\bar{T}) - x_{k_m}^m\| = \|\bar{x}^m(\bar{T}) - \bar{x}^m(T^m_{k_m})\| \leqslant 2(NL+L_g) |\bar{T} - T^m_{k_m}| \rightarrow 0$ as $m\rightarrow \infty$. Thus $x_{k_m}^m \in B(0,R+Q+1)$ for all sufficiently large $m$. By \cref{cor:Oalpha}, we have that $\|x_{k_m+1}^m  - x_{k_m}^m\| \leqslant 2\alpha_{k_m}^m (NL + L_g) \leqslant 1$, contradicting the assumption that $x_{k_m+1}^m \not\in B(0,R+Q+2)$.
\end{proof}

We are now ready to prove the desired result.
\begin{proof}[Proof of \cref{lemma:stoc_tracking}.]
      Let $X_0\subset \dom \Phi$ be nonempty and compact. Let $T>0$. By \cref{prop:bounded}, there exist $\hat{\alpha},r>0$ such that for any sequence generated by \cref{alg:prr} initialized in $X_0$ with $\alpha_0,\alpha_1,\hdots \leqslant \hat{\alpha}$, we have
\begin{equation*}
   \forall k\in \mathbb{N}^*,~~~ \alpha_0+\cdots+\alpha_{k-1} \leqslant T+1 ~~~\Longrightarrow~~~ x_0, \ldots, x_k \in B(0,r).
\end{equation*}
Let $(\bar{\alpha}_m)_{m\in \mathbb{N}}$ be a positive sequence that converges to zero. We assume that $\bar{\alpha}_m \leqslant \min\{1,\hat{\alpha}\}$. To each $\bar{\alpha}_m$, we attribute a sequence of iterates $(x_k^m)_{k\in \mathbb{N}}$ generated by \cref{alg:prr} with step sizes $\alpha_0^m, \alpha_1^m, \hdots \leqslant \bar{\alpha}_m$ and $x_0^m \in X_0$. 

Let $T^m_0 := 0$ and $T^m_k := \sum_{i = 0}^{k-1}\alpha_i^m$ for all $k\in \mathbb{N}^*$ and $m \in \mathbb{N}$. Let $T^m:= \min\{\sum_{k = 0}^\infty \alpha_k^m,T\}$. By taking a subsequence if necessary, we have that $T^m \to \bar{T} \in [0,T]$ as $m\to \infty$. For each sequence $(x_k^m)_{k\in \mathbb{N}}$, consider the (extended) linear interpolation $\bar{x}^m:[0,\max\{T^m,\bar{T}\}] \rightarrow \mathbb{R}^n$ defined by
	\begin{equation*}
		\bar{x}^m(t) = x_k^m + (t - T_k^m) \frac{x_{k+1}^m - x_{k}^m}{\alpha_k^m}
	\end{equation*}
	for any $t \in [T_k^m, \min\{T_{k+1}^m,T^m\}]$ and $k \in\mathbb{N}$ such that $T_k^m<T^m$. If $T^m<\bar{T}$, then $\bar{x}^m(t):= \lim_{s\nearrow T^m } \bar{x}^m(s)$ for all $t\in [T^m,\bar{T}]$. For any $m\in \mathbb{N}$, let $k_m := \sup\{k\in \mathbb{N}: T^m_k \leqslant T^m\}$. If $k_m = \infty$, then we have that $x_k^m \in B(0,r)$ for all $k\in \mathbb{N}$. Otherwise if $k_m<\infty$, then $T_{k_m+1}^m = T_{k_m}^m + \alpha_{k_m}^m \leqslant T^m + \bar{\alpha}_m \leqslant T+1$, and thus $x_0^m,\ldots,x_{k_m+1}^m \in B(0,r)$. In both cases, we know that $\bar{x}^m(t) \in B(0,r)$ for all $t\in [0,\max\{T^m,\bar{T}\}]$ by the convexity of $B(0,r)$. Using arguments similar to those in the proof of \cref{prop:bounded}, by passing to a subsequence if necessary, $(\bar{x}^{m}\big|_{[0,\bar{T}]}(\cdot))_{m\in \mathbb{N}}$ converges uniformly to a solution $x:[0,\bar{T}] \rightarrow \mathbb{R}^n$ to $x' \in - D_\Phi(x)$ with $x(0) = \lim_{m\rightarrow \infty}\bar{x}^m(0) = \lim_{m\rightarrow \infty}x_0^m \in X_0$. 
 
The conclusion of the lemma now follows. To see why, assume the contrary that there exists $\epsilon>0$ such that for any $\bar{\alpha}>0$, there exists $(x_k)_{k\in \mathbb{N}}$ generated by \cref{alg:prr} with step sizes $(\alpha_k)_{k\in \mathbb{N}} \subset (0,\bar{\alpha}]$ and $x_0 \in X_0$ such that for any absolutely continuous function $x:[0,T]\rightarrow \mathbb{R}^n$ that is a solution to $x' \in - D_\Phi(x)$ with $x(0)\in X_0$, there exists $K\in \mathbb{N}^*$ such that
    \begin{equation*}
   \alpha_0+\cdots+\alpha_{K-1} \leqslant T ~~~\text{and}~~~ \|x_K - x(\alpha_0+\cdots+\alpha_{K-1})\| > \epsilon.
\end{equation*}
Based on this, we may construct a sequence of iterates $(x_k^m)_{k\in \mathbb{N}}$ generated by \cref{alg:prr} with step sizes $(\alpha_k^m)_{k\in \mathbb{N}}$ and $x_0^m \in X_0$, where $\sup_k \alpha_k^m \to 0$ as $m \rightarrow \infty$. Moreover, for each $(x_k^m)_{k\in \mathbb{N}}$ and any solution $x:[0,T]\rightarrow \mathbb{R}^n$ to $x' \in - D_\Phi(x)$ with $x(0)\in X_0$, there exists $K_m\in \mathbb{N}^*$ such that
    \begin{equation*}
   \alpha^m_0+\cdots+\alpha^m_{K_m-1} \leqslant T ~~~\text{and}~~~ \|x^m_{K_m} - x(\alpha^m_0+\cdots+\alpha^m_{K_m-1})\| > \epsilon.
\end{equation*}
This is in contradiction with what we have shown above, namely the uniform convergence of the linear interpolations to a solution to $x' \in - D_\Phi(x)$.
\end{proof}

\subsection{Reachability of $(\epsilon,\delta)$-near approximate stationarity} \label{subsec:proof_exists_small}
\begin{proof}[Proof of \cref{thm:exists_small}.]
   We assume without loss of generality that $X_0 \subset \dom \Phi$ is nonempty and compact. Fix any $\delta,\epsilon>0$ and let $T:= (\sup_{y\in X_0}\Phi(y) - \inf_{x\in \mathbb{R}^n} \Phi(x))/\delta^2 \in \mathbb{R}_+$. Let $x:[0,T]\rightarrow \mathbb{R}^n$ be a solution to $x' \in - D_\Phi(x)$ with $x(0) \in X_0$. Following the same arguments as \eqref{eq:cont_len}, there exists $r>0$ such that $x([0,T]) \subset \dom \Phi \cap B(X_0,r)$ for any such $x(\cdot)$. According to  \cref{lemma:min_grad} and a standard compactness argument, $d(0,\partial g)$ is bounded over compact subsets of $\dom \Phi$. Thus, there exists $L>0$ such that $d(0,D_\Phi(y)) \leqslant L$ for all $y \in \dom \Phi \cap B(X_0,r)$.
    
    By \cref{lemma:stoc_tracking}, there exists $\bar{\alpha}\in (0,\epsilon/(2L)]$ such that for any sequence generated by \cref{alg:prr} with $\alpha_0,\alpha_1,\cdots \leqslant \bar{\alpha}$, $\sum_{k = 0}^\infty \alpha_k = \infty$, and $x_0 \in X_0$, there exists a solution $x:[0,T]\rightarrow \mathbb{R}^n$ to $x' \in - D_\Phi(x)$ such that $x(0) \in X_0$ and
\begin{equation*}
   \forall k\in \mathbb{N}^*,~~~ \alpha_0+\cdots+\alpha_{k-1} \leqslant T ~~~\Longrightarrow~~~ \|x_k - x(\alpha_0+\cdots+\alpha_{k-1})\| \leqslant \epsilon/2.
\end{equation*}
By \cref{lem:chain} it holds that
    \begin{equation*}
        \int_{0}^T d(0,D_\Phi(x(t)))^2~dt =  \Phi(x(0)) - \Phi(x(T)) \leqslant \sup_{y\in X_0}\Phi(y) - \inf_{x\in \mathbb{R}^n} \Phi(x).
    \end{equation*}
 Thus, there exists $t\in [0,T]$ such that $d(0,D_\Phi (x(t))) \leqslant \delta$. As $\alpha_k \leqslant \bar{\alpha}$, there exists $k\in \mathbb{N}$ such that $t_k := \sum_{i = 0}^{k-1}\alpha_i \in [\max\{0,t-\bar{\alpha}\},t]$.

 We next show that $\|x_k - x(t)\| \leqslant \epsilon$, then the conclusion of the theorem follows. Indeed,
 \begin{subequations}
     \begin{align}
        \|x_k - x(t)\| &\leqslant  \|x_k - x(t_k)\| + \|x(t_k) - x(t)\| \label{eq:dist_a}\\
        &\leqslant \epsilon/2 + \int_{t_k}^t \|x'(s)\|~ds \label{eq:dist_b}\\
        &= \epsilon/2 + \int_{t_k}^t d(0,D_\Phi(x(s)))~ds \label{eq:dist_c}\\
        &\leqslant \epsilon/2 + (t - t_k) L \label{eq:dist_d}\\
        &\leqslant \epsilon/2 + \bar{\alpha} L \leqslant \epsilon.\label{eq:dist_e}
     \end{align}
 \end{subequations}
Above, \eqref{eq:dist_a} and \eqref{eq:dist_b} follow from the triangular inequality. \eqref{eq:dist_c} follows again from \cref{lem:chain} and \eqref{eq:dist_d} follows from the fact that $d(0,D_\Phi)$ is locally bounded over $\dom \Phi$.
\end{proof}

\subsection{Convergence to $(\epsilon,\delta)$-near approximate stationarity} \label{subsec:proof_weakly_convex}
The proof of \cref{thm:weakly_convex} requires the following lemma, whose proof is omitted as it essentially follows from the proof of \citet[Theorem 6.2]{lewis2025identifiability}. We remark that weak convexity implies that the objective function is primal lower nice (\citet{poliquin1991integration}) everywhere in its domain, with the same constants.
\begin{lemma}
\label{lemma:essential_conv}
    Let $\Phi$ be proper, weakly convex, and lower bounded. For any subgradient trajectory $x:\mathbb{R}_+\rightarrow \mathbb{R}^n$ of $\Phi$ and for any $\epsilon>0$, there exists $T\geqslant 0$ such that $\|x'(t)\| \leqslant \epsilon$ for almost every $t\geqslant T$.
\end{lemma}
We proceed to prove the theorem.
\begin{proof}[Proof of \cref{thm:weakly_convex}.]
    Fix any $\delta,\epsilon>0$. As $\Phi$ is weakly convex and lower bounded, for any $x_0 \in \dom \Phi$, there exists a unique subgradient trajectory $x:\mathbb{R}_+\rightarrow \mathbb{R}^n$ of $\Phi$ initialized at $x_0$ (\citet{marcellin2006evolution}).  By \cref{lem:chain,lemma:essential_conv}, there exists $T_0>0$ such that $d(0,\partial \Phi(x(t))) = \|x'(t)\|\leqslant \delta$ for almost every $t \geqslant T_0$. Since $\gra \partial \Phi$ is closed, $\epi d(0,\partial \Phi)$ is closed. Fix $t \geqslant T_0$ and let $t_k\rightarrow t$ be such that $(x(t_k),\delta) \in \epi d(0,\partial \Phi)$. By continuity of $x(\cdot)$, $(x(t_k),\delta) \rightarrow (x(t),\delta) \in \epi d(0,\partial \Phi)$. In other words, $d(0,\partial \Phi(x(t))) \leqslant \delta$ actually holds for all $t \geqslant T_0$.

By the subdifferential regularity of $f_1,\ldots, f_N$ and the convexity of $g$, it holds that $\partial \Phi = \partial f_1 + \cdots + \partial f_N + \partial g$ (\citet[10.9 Corollary, p. 430]{rockafellar2009variational}). Thus by \cref{lemma:stoc_tracking}, for any $T\geqslant T_0$, there exists $\bar{\alpha}>0$ such that for any sequence generated by \cref{alg:prr} with $\alpha_0,\alpha_1,\hdots \leqslant \bar{\alpha}$, $\sum_{k = 0}^\infty \alpha_k = T$, and $x_0 \in \dom \Phi$, it holds that
\begin{equation*}
   \forall k\in \mathbb{N}^*,~~~ \|x_k - x(\alpha_0+\cdots+\alpha_{k-1})\| \leqslant \epsilon.
\end{equation*}
As a result, it holds that $\{x_k\}_{k\in \mathbb{N}} \subset B(x([0,T]),\epsilon)$, and that any limit point $x^*$ of the sequence must lie in $B(x(T),\epsilon)$. As $d(0, \partial \Phi(x(T))) \leqslant \delta$, it remains to show that the sequence $(x_k)_{k\in \mathbb{N}}$ is convergent. Indeed, by \cref{cor:Oalpha}, there exists $C>0$ such that $\|x_{k+1} - x_k\| \leqslant C\alpha_k$ for all $k\in \mathbb{N}$. Therefore, $\sum_{k = 0}^\infty \|x_{k+1} - x_k\| \leqslant C \sum_{k = 0}^\infty \alpha_k = CT$, thus the sequence $(x_k)_{k\in \mathbb{N}}$ is convergent. 
\end{proof}

\subsection{Convergence to $(\epsilon,0)$-near approximate stationarity}\label{subsec:proof_nonsmooth_bounded_flow}

The proof of \cref{thm:nonsmooth_bounded_flow} is a relatively direct consequence of the tracking lemma (\cref{lemma:stoc_tracking}). 

\begin{proof}[Proof of \cref{thm:nonsmooth_bounded_flow}.]

Fix any $\epsilon>0$ and $x_0 \in \dom \Phi$. According to the assumption, there exists a unique solution $x:\mathbb{R}_+\to \mathbb{R}^n$ to $x'(t) \in -D_\Phi(x(t))$ such that $x(0) = x_0$. By \cref{lem:chain} and the fact that $x(\cdot)$ is convergent, we know that $x^\sharp:= \lim_{t\to \infty} x(t)$ must be a $D_\Phi$-critical point. Also, there exists a $T_0>0$ such that $x([T_0,\infty))\subset B(x^\sharp,\epsilon/2)$.

For any $T\geqslant T_0$, by \cref{lemma:stoc_tracking}, there exists $\bar{\alpha}>0$ such that for any sequence generated by \cref{alg:prr} with $\alpha_0,\alpha_1,\hdots \leqslant \bar{\alpha}$ and $\sum_{k = 0}^\infty \alpha_k = T$, we have
\begin{equation}
\label{eq:close}
    \|x_k - x(\alpha_0 + \cdots + \alpha_{k-1})\| \leqslant \epsilon/2
\end{equation}
for any $k\in \mathbb{N}^*$. Using similar arguments as in the proof of \cref{thm:weakly_convex}, the sequence $(x_k)_{k\in \mathbb{N}}$ is convergent. Taking $k\rightarrow \infty$ in \eqref{eq:close}, we have that $\|x^* - x(T)\| \leqslant \epsilon/2$, where $x^*:= \lim_{k\rightarrow \infty} x_k$. Therefore, $\|x^* - x^\sharp\| \leqslant \|x^* - x(T)\| + \|x(T) - x^\sharp\| \leqslant \epsilon/2 + \epsilon/2 = \epsilon$.
\end{proof}

\subsection{Convergence to $(0,0)$-near approximate stationarity}\label{subsec:sgd_finite_len}

The proof of \cref{thm:sgd_finite_len} is derived from the following steps. We first establish an approximate descent lemma for the iterates generated by \cref{alg:prr} in \cref{lem:approx_descent}. A length formula is then proved in \cref{prop:discrete_length_formula} by using the approximate descent lemma (\cref{lem:approx_descent}), the uniform Kurdyka-\L{}ojasiewicz inequality (\cref{lemma:ukl}) and the upper bound on the distance between two successive iterations of \cref{alg:prr} (\cref{cor:Oalpha}). We then obtain the convergence of iterates in \cref{thm:sgd_finite_len}.

Different versions of approximate descent for \cref{alg:prr} have been established in the literature, even though it does not fit into the H1-H2-H3 framework in \citet{attouch2013convergence} or \citet{absil2005convergence}. Indeed, the sufficient decrease and relative error condition cannot be guaranteed due to random reshuffling. An approximate descent property of $\mathbb{E}[\Phi(x_k)]$ for \cref{alg:prr} was proved in \citet[E.2]{mishchenko2022proximal}. The approximate descent of $\Phi(x_k)$ for the special case of \cref{alg:prr} where $g=0$ was proved in \citet[Lemma 3.2]{li2023convergence}. In contrast to existing work, we do not require global Lipschitz continuity of $\nabla f_i$'s nor $g=0$. For our purposes, it is sufficient to restrict the iterates to a bounded subset, which incidentally makes the proof more direct.

\begin{lemma}\label{lem:approx_descent}
    Suppose \cref{assumption_standing,assumption_smooth} hold. Let $X$ be a bounded set. There exists $\bar{\alpha}>0$ such that if $(x_0,\ldots,x_{K+1}) \in X \times \cdots \times X \times \mathbb{R}^n$ is generated by \cref{alg:prr} with $\alpha_0,\hdots,\alpha_{K} \leqslant \bar{\alpha}$, then for all $k\in \llbracket 0,K\rrbracket$ we have
    \begin{equation*}
        \Phi(x_{k+1}) \leqslant \Phi(x_k)-\frac{\alpha_k}{4}d(0,\partial\Phi(x_{k+1}))^2-\frac{\|x_{k+1}-x_k\|^2}{8\alpha_k} + \frac{1}{12}(N-1)^2N(2N-1)L^2M^2\alpha_k^3 
    \end{equation*}
    where $L$ and $M$ are respectively Lipschitz constants of $f_1,\hdots,f_N$ on $B(X,2)$ and $\nabla f_1,\hdots,\nabla f_N$ on $\co B(X,1)$.
\end{lemma}
\begin{proof}
Let $L_g:=\sup_{x\in B(X,1)\cap \dom g}d(0,\partial g(x))$ and $\bar{\alpha}:=\min\{1/(2NL+2L_g),1/(2NM)\}$. Suppose $(x_0,\ldots,x_{K+1}) \in X \times \cdots \times X \times \mathbb{R}^n$ is generated by \cref{alg:prr} with $\alpha_0,\hdots,\alpha_{K} \leqslant \bar{\alpha}$. By \cref{cor:Oalpha}, $x_{K+1} \in B(X,1)$. Let $k \in \llbracket 0, K \rrbracket$. Since $NM$ is a Lipschitz constant of $\nabla f$ over $\co B(X,1)$, we have 
    \begin{equation*}
        f(x_{k+1}) \leqslant f(x_k) + \langle \nabla f(x_k), x_{k+1}-x_k \rangle + \frac{NM}{2}\|x_{k+1}-x_k\|^2. 
    \end{equation*}
    It holds also that $g(x_{k+1}) \leqslant g(x_k) + \langle g'(x_{k+1}), x_{k+1}-x_k\rangle$ for all $g'(x_{k+1})\in \partial g(x_{k+1})$ by convexity of $g$. Added to the previous inequality, this yields
    \begin{equation}\label{eq:approxdes_pre}
        \Phi(x_{k+1}) \leqslant \Phi(x_k) +\langle \nabla f(x_k)+ g'(x_{k+1}), x_{k+1}-x_k\rangle + \frac{NM}{2}\|x_{k+1}-x_k\|^2. 
    \end{equation}
    for all $g'(x_{k+1})\in \partial g(x_{k+1})$. We next handle the second term on the right-hand side. For all $g'(x_{k+1}) \in \partial g(x_{k+1})$, it holds that
    \begin{subequations}\label{eq:approxdes_inner}
        \begin{align}
                \langle \nabla f(x_k) + g'(x_{k+1}), x_{k+1}-x_k\rangle = & ~ \alpha_k \langle \nabla f(x_k) + g'(x_{k+1}), (x_{k+1}-x_k)/\alpha_k \rangle \\
                = & ~ \alpha_k \| \nabla f(x_k) + g'(x_{k+1})  + (x_{k+1}-x_k)/\alpha_k \|^2/2 \\
                & ~ - \alpha_k \| \nabla f(x_k) + g'(x_{k+1}) \|^2/2 - \|x_{k+1}-x_k\|^2/(2\alpha_k).
        \end{align}
    \end{subequations}
    On the one hand,
    \begin{align}
        d(0,\partial \Phi(x_{k+1}))^2 & = d(0,\nabla f(x_{k+1}) + \partial g(x_{k+1}))^2 \nonumber \\
        & \leqslant 2\| \nabla f(x_{k+1}) - \nabla f(x_k) \|^2 + 2d(0, \nabla f(x_k) + \partial g(x_{k+1}))^2 \nonumber \\
        & \leqslant 2\left(\sum_{i=1}^N \| \nabla f_i(x_{k+1}) - \nabla f_i(x_k) \|\right)^2 + 2\|\nabla f(x_k) + g'(x_{k+1})\|^2 \nonumber \\
        & \leqslant 2N^2M^2\| x_{k+1}-x_k \|^2 + 2\|\nabla f(x_k) + g'(x_{k+1})\|^2. \label{eq:approxdes_inner1}
    \end{align}
    On the other hand, since $x_{k+1}=\prox_{\alpha_kg}(x_{k,N})$, by Fermat's rule (\citet[Theorem 10.1]{rockafellar2009variational}) we have $0\in \alpha_k \partial g(x_{k+1}) + x_{k+1}-x_{k,N}$. Recall that $x_{k,N} = x_k - \alpha_k \sum_{i=1}^N\nabla f_{\sigma_i^k}(x_{k,i-1})$. Thus there exists $g'(x_{k+1}) \in \partial g(x_{k+1})$ such that 
    \begin{equation*}
        0 = g'(x_{k+1}) + \frac{x_{k+1}-x_k}{\alpha_k} + \sum_{i=1}^N\nabla f_{\sigma_i^k}(x_{k,i-1}).
    \end{equation*}
    By the QM-AM inequality (\citet{gwanyama2004hm}), we hence have
    \begin{align}
        \| \nabla f(x_k) + g'(x_{k+1})  + (x_{k+1}-x_k)/\alpha_k \|^2 & = \left\| \nabla f(x_k) - \sum_{i=1}^N\nabla f_{\sigma_i^k}(x_{k,i-1}) \right\|^2 \nonumber \\
         & = \left\| \sum_{i=1}^N \nabla f_{\sigma^k_i}(x_k) - \nabla f_{\sigma_i^k}(x_{k,i-1}) \right\|^2 \nonumber \\
        & \leqslant (N-1) \sum_{i=2}^N \| \nabla f_{\sigma^k_i}(x_k) - \nabla f_{\sigma_i^k}(x_{k,i-1}) \|^2 \nonumber \\ 
        & \leqslant (N-1) M^2 \sum_{i=2}^N \| x_k-x_{k,i-1} \|^2 \nonumber \\
        & = (N-1) M^2 \alpha_k^2 \sum_{i=2}^N \left\| \sum_{j=1}^{i-1} \nabla f_{\sigma^k_j}(x_{k,j-1})  \right\|^2 \nonumber \\
        & \leqslant (N-1) M^2 \alpha_k^2 \sum_{i=2}^N (i-1) \sum_{j=1}^{i-1} \|\nabla f_{\sigma^k_j}(x_{k,j-1})\|^2 \nonumber \\
        & \leqslant (N-1) L^2 M^2 \alpha_k^2 \sum_{i=1}^{N-1} i^2 \nonumber \\
        & \leqslant \frac{1}{6}(N-1)^2N(2N-1)L^2 M^2 \alpha_k^2. \label{eq:approxdes_inner2}
    \end{align}
    Above, we use the fact that $x_{k,i-1}\in B(X,1)$ for all $i \in \llbracket 1,N \rrbracket$ and $k\in \llbracket 0,K\rrbracket$ since $\alpha_k \leqslant 1/(NL)$ by \cref{cor:Oalpha}.

    Substituting \eqref{eq:approxdes_inner} into \eqref{eq:approxdes_pre}, and applying \eqref{eq:approxdes_inner1} and \eqref{eq:approxdes_inner2} to eliminate terms involving $g'(x_{k+1})$, we obtain 
    \begin{align*}
            \Phi(x_{k+1}) \leqslant & ~ \Phi(x_k)
          - \frac{\alpha_k}{4}d(0,\partial \Phi(x_{k+1}))^2 + \left(\frac{NM}{2} + \frac{N^2M^2}{2}\alpha_k - \frac{1}{2\alpha_k} \right) \|x_{k+1}-x_k\|^2 \\
          & ~ + \frac{1}{12}(N-1)^2N(2N-1)L^2 M^2 \alpha_k^3.
    \end{align*}
    This yields the desired inequality since $\alpha_k \leqslant 1/(2NM)$.
\end{proof}

With the approximate descent property, \citet{li2023convergence} first establish the convergence of $f(x_k)$ and $\|\nabla f(x_k)\|$ to obtain the convergence of $x_k$ provided that $g=0$ and the $f_i$'s have globally Lipschitz gradients. However, without either the global Lipschitz gradient continuity of $f_i$'s or boundedness of $x_k$, we cannot expect to first obtain the convergence of function value and gradient norm. Thus, a different approach inspired by \citet{josz2023global} is taken to directly show the convergence of $x_k$. In \citet[Proposition 8]{josz2023global}, the author proved a length formula for the gradient method, which is a special case of \cref{alg:prr} where $g=0$ and $N=1$. It is inspired by Kurdyka's original length formula (\citet[Theorem 2]{kurdyka1998gradients}). 

In order to establish a length formula for \cref{alg:prr}, we next state the uniform Kurdyka–\L ojasiewicz inequality (\citet[Proposition 5]{josz2023global}) for $D_\Phi$. The motivation behind the uniform Kurdyka–\L ojasiewicz inequality as opposed to the Kurdyka–\L ojasiewicz inequality is to extend the inequality to all points of a bounded set, without restricting the function value. We in fact propose a strengthened version of the uniform Kurdyka–\L ojasiewicz inequality for reasons discussed below. Also, we allow the inequality to hold at $D_\Phi$-critical points, which is new even for the Clarke subdifferential, and simplifies the subsequent analysis. Given $\Phi:\mathbb{R}^n\rightarrow \overline{\mathbb{R}}$ and $X\subset \mathbb{R}^n$, a scalar $v$ is called a $D_\Phi$-critical value of $\Phi$ in $X$ if there exists $x\in X$ such that $0\in D_\Phi(x)$ and $v = \Phi(x)$. In particular, if $D_\Phi=\partial \Phi$, we call this $v$ a critical value of $\Phi$ in $X$.

\begin{lemma}\label{lemma:ukl}
     Let \cref{assumption_standing} hold. Assume that $\Phi$ and $D_\Phi$ are semialgebraic. Let $X$ be a bounded subset of $\dom \Phi$ with $\overline{X}\subset\dom \Phi$. Define $V$ to be the set of $D_\Phi$-critical values of $\Phi$ in $\overline{X}$ if it is nonempty; otherwise $V:=\{0\}$. Let $\theta \in (0,1)$. There exist a concave semialgebraic diffeomorphism $\psi:\mathbb{R}_+\rightarrow\mathbb{R}_+$ and $\xi>0$ such that
    \begin{subequations}
        \begin{align}
         \forall\,x\in X,  \quad  &d(0,D_\Phi(x))\geqslant \frac{1}{\psi'(d(\Phi(x),V))}, \label{eq:ukl} \\[1mm]
         \forall s,t \geqslant 0, ~~~ &\frac{1}{\psi'(s+t)} \leqslant \frac{1}{\psi'(s)} + \max\left\{\frac{1}{\psi'(t)}, \xi t \right\}, \label{eq:bound_difference} \\[3mm]
         \forall t\geqslant 0, ~~~ &\psi'(t) \geqslant t^{-\theta}.
        \end{align}
    \end{subequations}
\end{lemma}
\begin{proof}
    Here, $\psi'(0)$ is the right derivative of $\psi$ at $0$, and we use the convention $1/\infty = 0$. Following the arguments in \citet[Theorem 6]{bolte2021conservative} (see also \citet[Corollary 15]{bolte2007clarke}) and the linear extension construction in \citet[Proposition 5]{josz2023global}, there exists a concave semialgebraic diffeomorphism $\psi:\mathbb{R}_+\rightarrow\mathbb{R}_+$ that satisfies \eqref{eq:ukl}. We may assume that $\psi'(t) \geqslant t^{-\theta}$ for all $t\geqslant 0$, after possibly replacing $\psi$ by $t\mapsto \int_0^t \max\{\psi'(s),s^{-\theta}\}ds$, which is concave semialgebraic. It is concave since the integrand is decreasing. To see why it is semialgebraic, note that $\{s >0 : \psi'(s)\geqslant s^{-\theta}\}$ is semialgebraic and hence a finite union of open intervals and points. Thus the integral is equal to $\psi$ up to a constant on finitely many intervals of $\mathbb{R}_+$, and equal to $t\mapsto t^{1-\theta}/(1-\theta)$ up to a constant elsewhere. The graph of such a function is hence semialgebraic. By the monotonicity theorem (\citet[(1.2) p. 43]{van1998tame}) and adapting the linear extension in \citet[Proposition 5]{josz2023global}, we may assume that $1/\psi'$ is continuously differentiable on $(0,T)$ and is constant on $(T,\infty)$ for some $T\in (0,\infty)$, with $(1/\psi')'$ monotone on $(0,T)$ and $\lim_{t\nearrow T}(1/\psi')'(t)\in (0,\infty)$. If $(1/\psi')'(t) \rightarrow \infty$ as $t\rightarrow 0$, then $(1/\psi')'$ is decreasing on $(0,T)$. By concavity of $\psi$, $(1/\psi')'(t) = -\psi''(t)/(\psi'(t))^2 \geqslant 0$ for all $t\in (0,T)$. As  $(1/\psi')'(t) = 0$ for all $t\in (T,\infty)$, $1/\psi'$ is concave on $\mathbb{R}_+$. Therefore, \eqref{eq:bound_difference} holds by \citet[Lemma 3.5]{josz2023convergence}. Otherwise if $\lim_{t\searrow 0}(1/\psi')'(t)<\infty$, then there exists $\xi>0$ such that $|(1/\psi')'(t)| \leqslant \xi$ for all $t\in (0,T)\cup (T,\infty)$. Therefore, $|1/\psi'(s+t) - 1/\psi'(s)| \leqslant \xi t$ for any $s,t\in \mathbb{R}_+$, and \eqref{eq:bound_difference} follows.
\end{proof}

We say that $\psi$ in \cref{lemma:ukl} is a desingularizing function of $\Phi$ over $X$ if it satisfies \eqref{eq:ukl} with $D_\Phi:= \partial \Phi$. The existing analysis of random reshuffling in \citet[Theorem 3.6]{li2023convergence} requires $1/\psi'$ to satisfy a quasi-additivity-type property, namely 
\begin{equation*}
   \forall s,t \in (0,\eta), ~~~ s+t < \eta ~~~ \Longrightarrow ~~~ \frac{1}{\psi'(s+t)} \leqslant C_\psi\left(\frac{1}{\psi'(s)} + \frac{1}{\psi'(t)}\right)
\end{equation*}
for some constants $\eta,C_\psi>0$. This is true for power functions with $C_\psi=1$ and any $\eta>0$. It is hence satisfied in polynomially bounded o-minimal structures  \cite{van1996geometric}. However, it is not clear why this should hold in general o-minimal structures. Thankfully, \cref{lemma:ukl} shows that one can actually get a somewhat weaker property for free, namely \eqref{eq:bound_difference}, which is sufficient for proving the length formula below. This is one of the key technical contributions of this work. The preceding argument is intended to ensure that ``semialgebraic'' can be replaced by ``definable'' throughout. 

\begin{proposition}\label{prop:discrete_length_formula}
Suppose \cref{assumption_standing,assumption_smooth} hold and $\Phi$ is semialgebraic. Let $X\subset\dom \Phi$ be bounded, $r\geqslant 1$, $\theta \in (0,1)$, and $m\in \mathbb{N}^*$ be an upper bound on the number of critical values of $\Phi$ in $\overline{X}$. Let $\psi$ be a desingularizing function of $\Phi$ over $X$ such that there exists $\xi\geqslant 2$ satisfying \eqref{eq:bound_difference} and $\psi'(t) \geqslant t^{-\theta}$ for all $t\geqslant 0$. There exist $c_1\geqslant 0$ and $\bar{\alpha},c_2>0$ such that for all $K\in \mathbb{N}$, if $(x_0,\ldots,x_{K+1}) \in X \times \cdots \times X \times \mathbb{R}^n$ is generated by \cref{alg:prr} with $\alpha_0,\hdots,\alpha_K \leqslant \bar{\alpha}$ and $\alpha_1/\alpha_0,\hdots,\alpha_{K+1}/\alpha_K \leqslant r$, then 
\begin{align}
\label{ga:length_formula}
    \frac{1}{2m} \sum_{k=0}^K\|x_{k+1}-x_k\| \leqslant & 2\sqrt{2}r\psi\left(\frac{1}{2m}\left[\Phi(x_0)-\Phi(x_K)+c_1\sum_{k=0}^{K-1}\alpha_k^3\right]\right) \\
    & ~ + \frac{c_1}{m}\sum_{k=0}^{K}\alpha_k\max\left\{\left(\sum_{i=k}^{K}\alpha_i^3\right)^\theta,\sum_{i=k}^{K}\alpha_i^3\right\} + c_2\max_{k\in\llbracket 0,K\rrbracket}\alpha_k. 
\end{align}
If $L,M \geqslant 1$ are respectively Lipschitz constants of $f_1,\hdots,f_N$ on $B(X,2)$ and $\nabla f_1,\hdots,\nabla f_N$ on $\co B(X,1)$, and $L_g:=\sup_{x\in B(X,1)\cap \dom g}d(0,\partial g(x))$, then one may choose $\bar{\alpha}:=1/(2\max\{NL+L_g,NM\})$,
\begin{equation}
    c_1 := \xi(N-1)^2N(2N-1)L^2M^2/12 ~~~ \text{and} ~~~ c_2 := 4(NL+L_g).
\end{equation}
\end{proposition}
\begin{proof}
    Let $K \in \mathbb{N}$ and $(x_0,\ldots,x_{K+1}) \in X \times \cdots \times X \times \mathbb{R}^n$ be generated by \cref{alg:prr} with $\alpha_0,\hdots,\alpha_K \leqslant \bar{\alpha}$ and $\alpha_1/\alpha_0,\hdots,\alpha_{K+1}/\alpha_K \leqslant r$. Define $y_k:=c_1\sum_{i=k}^K\alpha_i^3/\xi$ and $z_k:=\Phi(x_k)+y_k$ for all $k \in \llbracket 0,K\rrbracket$. By \cref{lem:approx_descent}, for all $k \in \llbracket 0,K-1\rrbracket$ we have
    \begin{equation}\label{eq:descent_property}
        z_{k+1} - z_k \leqslant -\frac{\alpha_k}{4}d(0,\partial\Phi(x_{k+1}))^2 -\frac{1}{8\alpha_k}\|x_{k+1}-x_k\|^2.
    \end{equation}
Let $V$ be the set of critical values of $\Phi$ in $\overline{X}$ if it is nonempty, and otherwise let $V:=\{0\}$. Since $\Phi$ is semialgebraic, $V$ has finitely many elements by the semialgebraic Morse-Sard theorem (\citet[Corollary 9]{bolte2007clarke}). Define $\tilde{z}_k := d(z_k,V)$ for all $k \in \llbracket 0,K\rrbracket$.

Assume that $[z_K,z_0)$ excludes the elements of $V$ and the averages of any two consecutive elements of $V$. Since $z_0\geqslant \cdots\geqslant z_K$, either $\tilde{z}_0 \leqslant \cdots \leqslant \tilde{z}_K$ or $\tilde{z}_0 \geqslant \cdots \geqslant \tilde{z}_K$. In the first case, for all $k \in \llbracket 0,K-1\rrbracket$, we have
    \begin{subequations}
    \label{eq:increasing}
        \begin{align}
        \psi(\tilde{z}_{k+1})-\psi(\tilde{z}_k) 
        \geqslant & ~ \psi'(\tilde{z}_{k+1})(\tilde{z}_{k+1}-\tilde{z}_k) \label{eq:increasing_b} \\
        = & ~ \psi'(\tilde{z}_{k+1})(z_k-z_{k+1}) \label{eq:increasing_c} \\
        \geqslant &~ \psi'(d(\Phi(x_{k+1})+y_{k+1},V))\left(\frac{\alpha_k}{4}d(0,\partial\Phi(x_{k+1}))^2+\frac{1}{8\alpha_k}\|x_{k+1}-x_k\|^2\right) \label{eq:increasing_d} \\
        \geqslant &~ \psi'(d(\Phi(x_{k+1}),V)+y_{k+1})\left(\frac{\sqrt{\alpha_k}}{2\sqrt{2}}d(0,\partial\Phi(x_{k+1}))+\frac{1}{4\sqrt{\alpha_k}}\|x_{k+1}-x_k\|\right)^2 \label{eq:increasing_e} \\
        \geqslant &~ \frac{\left(\frac{\sqrt{\alpha_k}}{2\sqrt{2}}d(0,\partial\Phi(x_{k+1}))+\frac{1}{4\sqrt{\alpha_k}}\|x_{k+1}-x_k\|\right)^2}{\frac{1}{\psi'(d(\Phi(x_{k+1}),V))}+\max\left\{\frac{1}{\psi'(y_{k+1})},\xi y_{k+1}\right\}} \label{eq:increasing_f} \\
        \geqslant &~ \frac{\left(\frac{\sqrt{\alpha_k}}{2\sqrt{2}}d(0,\partial\Phi(x_{k+1}))+\frac{1}{4\sqrt{\alpha_k}}\|x_{k+1}-x_k\|\right)^2}{d(0,\partial\Phi(x_{k+1}))+\max\left\{\frac{1}{\psi'(y_{k+1})},\xi y_{k+1}\right\}}. \label{eq:increasing_g}
        \end{align}
    \end{subequations}
    Indeed, \eqref{eq:increasing_b} holds because $\psi$ is concave. \eqref{eq:increasing_c} is due to the existence of $v \in V$ such that $\tilde{z}_k = v - z_k$ for all $k \in \llbracket 0,K \rrbracket$. \eqref{eq:increasing_d} follows by the descent property \eqref{eq:descent_property}. \eqref{eq:increasing_e} uses the fact that 
    $\psi'$ is decreasing and 
    \begin{align*}
        d(\Phi(x_{k+1})+y_{k+1},V) &= \min_{v\in V} |\Phi(x_{k+1})+y_{k+1}-v|  \\
        &\leqslant \min_{v\in V} |\Phi(x_{k+1})-v|+y_{k+1} = d(\Phi(x_{k+1}),V)+y_{k+1}
    \end{align*}
    for all $k\in \llbracket 0,K-1\rrbracket$. It also uses the QM-AM inequality (\citet{gwanyama2004hm}). \eqref{eq:increasing_f} is an application of \eqref{eq:bound_difference}. Finally, \eqref{eq:increasing_g} is a consequence of the uniform Kurdyka-\L{}ojasiewicz inequality \eqref{eq:ukl}. Using the AM-GM inequality, \eqref{eq:increasing} yields
    \begin{subequations}
    \label{eq:am-gm}
        \begin{align}
        &\frac{\sqrt{\alpha_k}}{2\sqrt{2}}d(0,\partial\Phi(x_{k+1})) +\frac{1}{4\sqrt{\alpha_k}}\|x_{k+1}-x_k\| \\
        &\quad \leqslant  \sqrt{(\psi(\tilde{z}_{k+1})-\psi(\tilde{z}_k))\left(d(0,\partial\Phi(x_{k+1}))+\max\left\{\frac{1}{\psi'(y_{k+1})},\xi y_{k+1}\right\}\right)} \\
        &\quad\leqslant  \frac{\psi(\tilde{z}_{k+1})-\psi(\tilde{z}_k)}{2(1-\iota)\sqrt{\alpha_k}} + \frac{(1-\iota)\sqrt{\alpha_k}}{2}\left(d(0,\partial\Phi(x_{k+1}))+\max\left\{\frac{1}{\psi'(y_{k+1})},\xi y_{k+1}\right\}\right) \label{eq:lenfor_recur_incre}
        \end{align}
    \end{subequations}
    where \(\iota\in[0,1)\) can be arbitrary. By letting $\iota=1-1/\sqrt{2}$ and multiplying by $4\sqrt{\alpha_k}$ on both sides, we obtain 
    \begin{equation*}
        \|x_{k+1}-x_k\| \leqslant 2\sqrt{2}(\psi(\tilde{z}_{k+1})-\psi(\tilde{z}_k))+2\max\left\{\frac{\alpha_k}{\psi'(y_{k+1})},\xi\alpha_ky_{k+1}\right\}.
    \end{equation*}
    Telescoping yields
    \begin{subequations}
        \begin{align}
        \sum_{k=0}^K\|x_{k+1}-x_k\| 
        \leqslant &~ 2\sqrt{2}(\psi(\tilde{z}_K)-\psi(\tilde{z}_0))+2\sum_{k=0}^{K-1}\max\left\{\frac{\alpha_k}{\psi'(y_{k+1})},\xi\alpha_ky_{k+1}\right\} + \|x_{K+1}-x_K\| \label{eq:telescope_a} \\
        \leqslant &~ 2\sqrt{2}\psi(\tilde{z}_K-\tilde{z}_0) +2\sum_{k=0}^K \max\left\{\frac{\alpha_k}{\psi'(y_k)},\xi\alpha_ky_k\right\} + 2(NL+L_g) \alpha_K \label{eq:telescope_b} \\
        \leqslant &~ 2\sqrt{2}\psi(z_0-z_K)+2\sum_{k=0}^K\max\left\{\frac{\alpha_k}{\psi'(y_k)},\xi\alpha_ky_k\right\} + 2(NL+L_g) \max_{k\in\llbracket0,K\rrbracket} \alpha_k
        \label{eq:telescope_c}
        \end{align}
    \end{subequations}
    where \eqref{eq:telescope_b} uses \citet[Lemma 3.5]{josz2023convergence}, $y_k \geqslant y_{k+1}$, and \cref{cor:Oalpha}. 

    In the second case, i.e., $\tilde{z}_0>\cdots>\tilde{z}_K$, \eqref{eq:increasing} becomes
    \begin{equation*}
        \psi(\tilde{z}_k)-\psi(\tilde{z}_{k+1}) \geqslant \frac{\left(\frac{\sqrt{\alpha_k}}{2\sqrt{2}}d(0,\partial\Phi(x_{k+1}))+\frac{1}{4\sqrt{\alpha_k}}\|x_{k+1}-x_k\|\right)^2}{d(0,\partial\Phi(x_k))+\max\left\{\frac{1}{\psi'(y_k)},\xi y_k\right\}}
    \end{equation*}
    and \eqref{eq:am-gm} becomes
    \begin{align*}
        &\frac{\sqrt{\alpha_k}}{2\sqrt{2}}d(0,\partial\Phi(x_{k+1})) +\frac{1}{4\sqrt{\alpha_k}}\|x_{k+1}-x_k\| \\
        &\quad \leqslant
        \frac{1}{2(1-\iota)\sqrt{\alpha_k}}(\psi(\tilde{z}_k)-\psi(\tilde{z}_{k+1})) + \frac{(1-\iota)\sqrt{\alpha_k}}{2}\left(d(0,\partial\Phi(x_k))+\max\left\{\frac{1}{\psi'(y_k)},\xi y_k\right\}\right).
    \end{align*}
    Setting $\iota=1-r^{-1}/\sqrt{2}$ and multiplying $4\sqrt{\alpha_k}$ on both sides yields 
    \begin{align*}
        &\sqrt{2}r^{-1}\alpha_{k+1}d(0,\partial\Phi(x_{k+1})) +\|x_{k+1}-x_k\| \\
        &\quad \leqslant 2\sqrt{2}r(\psi(\tilde{z}_k)-\psi(\tilde{z}_{k+1}))+\sqrt{2}r^{-1}\alpha_k\left(d(0,\partial\Phi(x_k))+\max\left\{\frac{1}{\psi'(y_k)},\xi y_k\right\}\right) 
    \end{align*}
    where we use the fact that $\alpha_{k+1}/\alpha_k \leqslant r$. This can be simplified to 
    \begin{align*}
        \|x_{k+1}-x_k\| \leqslant & 2\sqrt{2}r[\psi(\tilde{z}_k)-\psi(\tilde{z}_{k+1})] + \sqrt{2}r^{-1}[\alpha_kd(0,\partial\Phi(x_k))-\alpha_{k+1}d(0,\partial\Phi(x_{k+1}))] \\ &~ +\sqrt{2}\max\left\{\frac{\alpha_k}{\psi'(y_k)},\xi\alpha_ky_k\right\}. 
    \end{align*}
    Telescoping and $r\geqslant 1$ yield
    \begin{align*}
        \sum_{k=0}^K\|x_{k+1}-x_k\| 
        \leqslant & ~ 2\sqrt{2}r(\psi(\tilde{z}_0)-\psi(\tilde{z}_K)) + \sqrt{2}r^{-1}[\alpha_0 d(0,\partial\Phi(x_0))-\alpha_K d(0,\partial\Phi(x_K))] \\
        & ~ + \sqrt{2}\sum_{k=0}^K \max\left\{\frac{\alpha_k}{\psi'(y_k)},\xi\alpha_ky_k\right\} + \|x_{K+1} - x_K\| \\
        \leqslant &~ 2\sqrt{2}r(\psi(\tilde{z}_0)-\psi(\tilde{z}_K)) + \sqrt{2}r^{-1}\alpha_0d(0,\partial\Phi(x_0)) \\
        & ~ + \sqrt{2}\sum_{k=0}^K\max\left\{\frac{\alpha_k}{\psi'(y_k)},\xi\alpha_ky_k\right\} + \|x_{K+1} - x_K\|\\
        \leqslant &~ 2\sqrt{2}r\psi(\tilde{z}_0-\tilde{z}_K) + \sqrt{2}\alpha_0(\|\nabla f(x_0)\|+d(0,\partial g(x_0))) \\
        & ~ + \sqrt{2}\sum_{k=0}^K \max\left\{\frac{\alpha_k}{\psi'(y_k)},\xi\alpha_ky_k\right\} + 2(NL+L_g)\alpha_K \\
        \leqslant &~ 2\sqrt{2}r\psi(z_0-z_K) + \sqrt{2}(NL+L_g)\alpha_0 \\
        & ~ + \sqrt{2}\sum_{k=0}^K \max\left\{\frac{\alpha_k}{\psi'(y_k)},\xi\alpha_ky_k\right\} + 2(NL+L_g)\alpha_K \\
        \leqslant &~ 2\sqrt{2}r\psi(z_0-z_K) + 2\sum_{k=0}^K \max\left\{\frac{\alpha_k}{\psi'(y_k)},\xi\alpha_ky_k\right\} + 4(NL+L_g)\max_{k\in \llbracket 0,K\rrbracket} \alpha_k.
    \end{align*}
    Compared with the upper bound obtained in the first case \eqref{eq:telescope_c}, the above bound is larger. We can hence use it as a common bound for both cases.

    We next consider the case where $z_0>\cdots>z_K$ and there exist $0 \leqslant K_1 \leqslant \hdots \leqslant K_p \leqslant K-1$ such that
    \begin{equation*}
        [z_K,z_{K_p+1}) \cup \hdots \cup [z_{K_2},z_{K_1+1}) \cup [z_{K_1},z_0)
    \end{equation*}
    excludes the elements of $V$ and the averages of any two consecutive elements of $V$. For notational convenience, let $K_0 := -1$ and $K_{p+1} := K$. We have 
\begin{align*}
    \sum_{k=0}^K \|x_{k+1}-x_k\| 
    = & ~ \sum_{i=0}^p \sum_{k=K_i+1}^{K_{i+1}} \|x_{k+1}-x_k\| \\
    \leqslant & ~ \sum_{i=0}^p \left( 2\sqrt{2}r\psi(z_{K_i+1}-z_{K_{i+1}}) + 2\sum_{k=K_i+1}^{K_{i+1}} \max\left\{\frac{\alpha_k}{\psi'(y_k)},\xi\alpha_ky_k\right\} \right. \\
    & ~ +\left. 4(NL+L_g)\max_{k\in \llbracket K_i+1,K_{i+1}\rrbracket} \alpha_k \right) \\
    \leqslant & ~ 2\sqrt{2}r(p+1)\psi\left(\frac{1}{p+1}\sum_{i=0}^p z_{K_i+1}-z_{K_{i+1}}\right) + 2\sum_{k=0}^K \max\left\{\frac{\alpha_k}{\psi'(y_k)},\xi\alpha_ky_k\right\} \\
    & ~ +4(p+1)(NL+L_g)\max_{k\in \llbracket 0,K\rrbracket} \alpha_k  \\
    \leqslant & ~ 2\sqrt{2}r(p+1)\psi\left(\frac{z_0-z_K}{p+1}\right) + 2\sum_{k=0}^K \max\left\{\frac{\alpha_k}{\psi'(y_k)},\xi\alpha_ky_k\right\} \\
    & ~ + 4(p+1)(NL+L_g)\max_{k\in \llbracket 0,K\rrbracket} \alpha_k \\
    \leqslant & ~ 4\sqrt{2}mr~\psi\left(\frac{z_0-z_K}{2m}\right) + 2\sum_{k=0}^K \max\left\{\frac{\alpha_k}{\psi'(y_k)},\xi\alpha_ky_k\right\} \\
    & ~ + 8m(NL+L_g)\max_{k\in \llbracket 0,K\rrbracket} \alpha_k, 
\end{align*}
where the second inequality uses concavity of $\psi$ and the last inequality uses $p\leqslant 2m-1$ and the fact that $s\mapsto s\psi(a/s)$ is increasing over $(0,\infty)$ for any constant $a>0$ (\citet[Proposition 2.1]{marechal2005functional}). 

We now consider the general case where $z_0 \geqslant \cdots \geqslant z_K$. Observe that if $z_{k+1}=z_k$ at some iteration $k \in \llbracket 0,K-1\rrbracket$, then $x_{k+1}=x_k$ and $0 \in \partial \Phi(x_{k+1})$ by \eqref{eq:descent_property}. Hence such iterations do not contribute to the length, and happen only at critical points of $\Phi$. We can thus remove them and obtain the above length formula with the remaining indices, to which we can add back the discarded indices. Indeed, the added terms on the left hand side would be equal to zero, and the added terms on the right-hand side would only make the bound greater.

The desired formula now follows from simple observations. First,
\begin{equation*}
    z_0-z_K  = \Phi(x_0) + y_0 -\Phi(x_K) - y_K = \Phi(x_0) - \Phi(x_K) + \frac{c_1}{\xi}\sum_{k=0}^{K-1}\alpha_k^3.
\end{equation*}
Second, since $\psi'(t) \geqslant t^{-\theta}$ for all $t\geqslant 0$ and $y_k \leqslant \xi y_k = c_1 \sum_{i=k}^{K}\alpha_i^3$, we have
\begin{align*}
    \max\left\{\frac{\alpha_k}{\psi'(y_k)},\xi\alpha_ky_k\right\} &= \alpha_k \max\left\{\frac{1}{\psi'(y_k)},\xi y_k\right\}\\
    &\leqslant \alpha_k \max\left\{y_k^\theta,\xi y_k\right\}\\
    &\leqslant \alpha_k \max\left\{c_1^\theta \left(\sum_{i=k}^{K}\alpha_i^3\right)^\theta,c_1 \sum_{i=k}^{K}\alpha_i^3\right\}\\
    &\leqslant c_1 \alpha_k \max\left\{\left(\sum_{i=k}^{K}\alpha_i^3\right)^\theta,\sum_{i=k}^{K}\alpha_i^3\right\}.
\end{align*} 
The last inequality holds because either $c_1 = 0$ or $c_1 \geqslant 1$. Indeed, recall that $c_1:=\xi(N-1)^2N(2N-1)L^2M^2/12$. If $N=1$, then $c_1=0$ and both sides of the last inequality are zero, so it holds trivially; otherwise $N\geqslant 2$ and $c_1\geqslant \xi\cdot 1^2\cdot 2\cdot 3L^2M^2/12 = \xi L^2M^2/2$. Recall that we assume $\xi\geqslant 2$ and $M,L\geqslant 1$ in the statement of this Proposition, so $c_1\geqslant 2\cdot 1^2\cdot 1^2/2=1$ and $c_1^\theta\leqslant c_1$ for all $\theta\in(0,1)$ yields the last inequality. 
\end{proof}

The above length formula finds a particularly simple form if we borrow the step sizes used in \citet[Lemma 3.7]{li2023convergence}. It then agrees with the one we derived for the momentum method (\citet[Lemma 2.2]{josz2023convergence}). It is hence suitable for obtaining global convergence.

\begin{corollary}
    \label{cor:length}
    Suppose \cref{assumption_standing,assumption_smooth} hold and that $\Phi$ is semialgebraic. Let $\beta = 0$ if $N=1$, and $\beta \in (1/2,1)$ otherwise. Let $X\subset \dom \Phi$ be a bounded set. There exist $\bar{\alpha}>0$, $\eta,\kappa \geqslant 0$, and a desingularizing function $\psi$ of $\Phi$ over $X$ such that for all $K\in \mathbb{N}$, $\alpha\in(0,\bar{\alpha}]$, and $\gamma \in\mathbb{N}^*$, if $(x_0,\ldots,x_{K+1}) \in X \times \cdots \times X \times \mathbb{R}^n$ is generated by \cref{alg:prr} with $\alpha_k=\alpha/(k+\gamma)^\beta$, then 
\begin{equation*}
    \sum_{k=0}^K\|x_{k+1}-x_k\| \leqslant \psi\left(\Phi(x_0)-\Phi(x_K)+ \eta \alpha\right) + \kappa\alpha. 
\end{equation*}
\end{corollary}
\begin{proof}
    Since $\alpha_{k+1}/\alpha_k \leqslant \max\{1,[(k+\gamma)/(k+\gamma+1)]^\beta\} \leqslant 1$ for all $k\in\llbracket 0,K\rrbracket$, by \cref{prop:discrete_length_formula} the formula \eqref{ga:length_formula} holds with $r:=1$ and whatever $\theta \in (0,1)$ we desire, with some corresponding $\xi \geqslant 2$ and $m\in \mathbb{N}^*$. Replacing $\psi$ with $4\sqrt{2}m\psi(\cdot/(2m))$, it holds that
    \begin{align*}
    \label{ga:length_cor}
       \sum_{k=0}^K\|x_{k+1}-x_k\| \leqslant &\psi\left(\Phi(x_0)-\Phi(x_K)+c_1\sum_{k=0}^{K-1}\alpha_k^3\right) 
        + 2 c_1\sum_{k=0}^{K}\alpha_k\max\left\{\left(\sum_{i=k}^{K}\alpha_i^3\right)^\theta,\sum_{i=k}^{K}\alpha_i^3\right\} \\
        & ~ + 2m c_2\max_{k\in\llbracket 0,K\rrbracket}\alpha_k. 
    \end{align*}
    If $N=1$, then it suffices to take $\eta := c_1 = 0$ and $\kappa := 2m c_2$. If $N>1$, then let $\theta\in(\max\{(1-\beta)/(3\beta-1),1/2\},1)$. For any $i\in\mathbb{N}$, using classical series integral comparison arguments (see, e.g. \citet[Appendix B]{li2023convergence}), we have
    \begin{equation}
    \label{eq:l3-bound}
        \sum_{k=i}^\infty \alpha_k^3 = \sum_{k=i}^\infty\frac{\alpha^3}{(k+\gamma)^{3\beta}} \leqslant \alpha_i^3+\int_{i+\gamma}^{\infty}\frac{\alpha^3}{\nu^{3\beta}}\;d\nu \leqslant \alpha_i^3+\alpha^3\frac{(i+\gamma)^{1-3\beta}}{3\beta-1}. 
    \end{equation}
    Without loss of generality, we may assume that $\bar{\alpha} \leqslant 1/2$. Since $\beta\in(1/2,1)$, we obtain an upper bound on the argument of $\psi$ in \eqref{ga:length_formula}:
    \begin{equation*}
        \sum_{k=0}^K \alpha_k^3 \leqslant \sum_{k=0}^{\infty}\alpha_k^3 \leqslant \frac{\alpha^3}{\gamma^{3\beta}}+\alpha^3\frac{\gamma^{1-3\beta}}{3\beta-1} \leqslant \alpha^3+\alpha^3\frac{1}{3\beta-1} \leqslant 3\alpha^3 \leqslant \alpha,
    \end{equation*}
    where the second inequality is obtained by setting $i=0$ in \eqref{eq:l3-bound} and $k=0$ in $\alpha_k=\alpha/(k+\gamma)^\beta$. We next upper bound the second term on the right-hand side of \eqref{ga:length_formula}. Note that $\sum_{i=k}^K\alpha_i^3 < (\sum_{i=k}^K\alpha_i^3)^\theta$ as $\sum_{i=k}^K\alpha_i^3,\theta\in (0,1)$. Therefore,
    \begin{align*}
        \sum_{k=0}^{K}\alpha_k\max\left\{\left(\sum_{i=k}^{K}\alpha_i^3\right)^\theta,\sum_{i=k}^{K}\alpha_i^3\right\} &= \sum_{k=0}^{K}\alpha_k\left(\sum_{i=k}^{K}\alpha_i^3\right)^\theta \\
        &\leqslant \sum_{k=0}^K\frac{\alpha}{(k+\gamma)^\beta}\left(\alpha_k^3+\alpha^3\frac{(k+\gamma)^{1-3\beta}}{3\beta-1}\right)^\theta \\
        &\leqslant \sum_{k=0}^K\frac{\alpha}{(k+\gamma)^\beta}\left(\alpha_k^{3\theta}+\alpha^{3\theta}\frac{(k+\gamma)^{(1-3\beta)\theta}}{(3\beta-1)^\theta}\right) \\
        &\leqslant \sum_{k=0}^K\frac{\alpha^{3\theta+1}}{(k+\gamma)^{\beta(1+3\theta)}} + \frac{\alpha^{3\theta+1}}{(3\beta-1)^\theta}\sum_{k=0}^K(k+\gamma)^{(1-3\beta)\theta-\beta}. 
    \end{align*}
    Similar to the upper bound of $\sum_{i=0}^K\alpha_i^3$ in \eqref{eq:l3-bound}, we have
    \begin{equation*}
        \sum_{k=0}^K\frac{\alpha^{3\theta+1}}{(k+\gamma)^{\beta(1+3\theta)}} \leqslant \alpha^{3\theta+1}+\alpha^{3\theta+1}\frac{1}{(3\theta+1)\beta-1}
    \end{equation*}
    due to $\theta \in (1/2,1)$, and
    \begin{align*}
        \frac{\alpha^{3\theta+1}}{(3\beta-1)^\theta}\sum_{k=0}^K(k+\gamma)^{(1-3\beta)\theta-\beta} &\leqslant 2^\theta\alpha^{3\theta+1}\left(\frac{1}{\gamma^{\beta+(3\beta-1)\theta}}+\frac{\gamma^{1-\beta-(3\beta-1)\theta}}{\beta+(3\beta-1)\theta-1}\right) \\
        &\leqslant 2\alpha^{3\theta+1}\left(1+\frac{1}{\beta+(3\beta-1)\theta-1}\right). 
    \end{align*}
    Since $\alpha\in(0,1/2]$, we have $\alpha^{3\theta + 1} \leqslant \alpha$ and
    \begin{align*}
         \sum_{k=0}^{K}\alpha_k\max\left\{\left(\sum_{i=k}^{K}\alpha_i^3\right)^\theta,\sum_{i=k}^{K}\alpha_i^3\right\}  
         &\leqslant \alpha^{3\theta + 1}\left(3+\frac{1}{(3\theta+1)\beta-1}+\frac{2}{\beta+(3\beta-1)\theta-1}\right)\\
         &\leqslant \alpha\left(3+\frac{1}{(3\theta+1)\beta-1}+\frac{2}{\beta+(3\beta-1)\theta-1}\right) \\
        &\leqslant 3\alpha\left(1+\frac{1}{\beta+(3\beta-1)\theta-1}\right). 
    \end{align*}
The desired inequality then follows by posing
\begin{equation*}
    \eta := c_1 ~~~ \text{and} ~~~ \kappa :=  6c_1\left(1+\frac{1}{\beta+(3\beta-1)\theta-1}\right) + 2mc_2. 
\end{equation*}
\end{proof}

We are now ready to prove the main result of this subsection. Leveraging the tracking lemma (\cref{lemma:stoc_tracking}) and the length formula (\cref{cor:length}), we employ a framework similar to that used in \citet[Theorem 2.8]{josz2023convergence} to establish the desired result, \cref{thm:sgd_finite_len}. To provide a clear overview, we outline the proof below, with the full details available in \cref{app:thm4proof}. 

\begin{proof}[Sketch proof of \cref{thm:sgd_finite_len}.]
Let $X_0$ be a bounded subset of $\dom \Phi$. We can first use the tracking lemma (\cref{lemma:stoc_tracking}) and the length formula (\cref{cor:length}) to prove the existence of $\bar{\alpha}>0$ such that $\Gamma(X_0,\bar{\alpha})<\infty$ where 
\begin{subequations}\label{eq:sup_trajectory_dis_sgd_sk}
\begin{align}
    \Gamma(X_0,\bar{\alpha}) := & \sup\limits_{\substack{x\in (\mathbb{R}^n)^{\mathbb{N}\times \llbracket 0,N\rrbracket},\sigma\in\mathfrak{S}_{N}^{\mathbb{N}}\\ (\alpha,\gamma)\in(0,\bar{\alpha}]\times \mathbb{N}^*}} ~~ \sum_{k=0}^{\infty} \|x_{k+1,0}-x_{k,0}\| \\
  & ~~ \mathrm{s.t.} ~~~ 
\left\{ 
\begin{array}{l}
x_{k+1,0}=\prox_{\alpha_k g}(x_{k,N}),\\[1mm]
x_{k,i} = x_{k,i-1}  - \alpha_k \nabla f_{\sigma^k_i} (x_{k,i-1}), ~ \forall i \in \llbracket 1,N \rrbracket, \\
\alpha_k = \alpha/(k+\gamma)^\beta,~\forall k\in\mathbb{N}, ~ x_{0,0} \in X_0.
\end{array}
\right.
\end{align}
\end{subequations}
Above, $\mathfrak{S}_{N}$ denotes the symmetric group of degree $N$. 

The finiteness of $\Gamma(X_0,\bar{\alpha})$ naturally means that $\|x_k-x_0\| \leqslant \sum_{i=0}^{k-1}\|x_{i+1}-x_i\| \leqslant \Gamma(X_0,\bar{\alpha})$ for all $k\in\mathbb{N}$, that is, $x_k \in B(X_0,\Gamma(X_0,\bar{\alpha}))$ (where $x_k := x_{k,0}$). Let $L_0,M$ respectively be Lipschitz constants of $f_1,\hdots,f_N$ and $\nabla f_1,\hdots,\nabla f_N$ on $B(X_0,\Gamma(X_0,\bar{\alpha})+2)$. We then establish 
\begin{align*}
    d(0,\partial \Phi(x_{k+1})) \leqslant \left(\frac{1}{\alpha_k}+M\right)\|x_{k+1}-x_k\| + (N-1)\sqrt{\frac{N(2N-1)}{6}}ML_0 \alpha_k.
\end{align*}
Note that for all nonnegative decreasing sequence $u_0,u_1,u_2,\hdots$, we have $u_k \leqslant 2 (\sum_{i=\lfloor k/2\rfloor}^\infty u_i)/(k+2)$ as explained in \citet[Footnote 1]{josz2023global}. In particular, for $u_k := \min_{i\in \llbracket 0,k \rrbracket} \alpha_i d(0,\partial \Phi(x_{i+1}))$, we have
\begin{align*}
    \min_{i\in \llbracket 0,k \rrbracket} \alpha_i d(0,\partial \Phi(x_{i+1})) \leqslant \frac{2}{k+2} \sum_{i=\lfloor k/2\rfloor}^\infty (1+\alpha M)\|x_{i+1}-x_i\| +(N-1)\sqrt{\frac{N(2N-1)}{6}}ML_0\alpha_i^2.
\end{align*}
Thus, using $\alpha_k = \alpha/(k+1)^\beta$ and  series integral comparison arguments again, we can finally obtain 
\begin{small}
 \begin{align*}
    \min_{i\in \llbracket 0,k \rrbracket} d(0,\partial \Phi(x_{i+1})) 
    \leqslant \frac{2}{(k+1)^{1-\beta}} \left((\alpha^{-1}+M)\sum_{i=\lfloor k/2\rfloor}^\infty \|x_{i+1}-x_i\| + \frac{2(N-1)\sqrt{N(2N-1)}ML_0\alpha}{\sqrt{6}(2\beta-1)(\lfloor k/2\rfloor+1)^{2\beta-1}}\right). 
\end{align*}   
\end{small}
\end{proof}

It is worth noting that our proof of \cref{thm:sgd_finite_len} is based on a worst-case analysis, which holds for any permutation of indices used in the stochastic gradient updates. As a consequence, the resulting complexity bound may not be tight with respect to the sample size~$N$, and the step size upper bound $\bar{\alpha}$ requires a conservative choice to ensure the iterates track the subgradient trajectories accurately. In particular, while \cref{lem:approx_descent} sets $\bar{\alpha} := \min\{1/(2NL+2L_g), 1/(2NM)\}$, which yields the typical $\Theta(1/N)$ scaling seen in nonconvex settings (see, e.g., \citet[Eq.~(2.16)]{yu2023high}), the final value of $\bar{\alpha}$ used in the proof of \cref{thm:sgd_finite_len} may be significantly smaller. This tighter restriction is necessary to guarantee the recursive tracking of subgradient trajectories, particularly in the absence of additional assumptions.

One natural question is whether a high-probability analysis could improve the dependence on~$N$ or allow for a larger step size. Indeed, such improvements have been demonstrated in the literature. For example, applying the approach in \citet{yu2023high} to PRR under their global smoothness assumptions, one can use a constant step size and improve the last-iterate complexity from $O(N^{3/2} \epsilon^{-3})$ (deterministic) to $O(\max\{N \epsilon^{-2}, N^{1/2} \epsilon^{-3}\})$ (high probability). Correspondingly, the step size improves from $O(1/N^{3/2})$ to $O(1/N)$.

However, in our setting, we do not assume a global Lipschitz constant for~$\nabla f_i$ or a uniform bound on the iterates. To ensure boundedness and thereby control the local Lipschitz constants, we need to adopt step sizes that diminish faster than $1/\sqrt{k}$. Under such conditions, high-probability analysis does not lead to an explicit improvement in either the allowable step size or the last-iterate complexity. Indeed, we explored both worst-case and high-probability analyses and found that, when using step sizes close to $1/\sqrt{k}$, the last-iterate complexity remains approximately $O(\epsilon^{-6})$. This complexity is quite poor, reflecting significantly slower convergence of the last iterate compared to the convergence rate in terms of minimum gradient norm. Moreover, the dependence on $N$ remains unclear. For these reasons, we chose not to include any details in this work.

\section{Conclusion}\label{sec:conclusion}

We analyzed the proximal random reshuffling (PRR) algorithm for minimizing a finite sum of locally Lipschitz functions $f_i$ and a convex function $g$ that is locally Lipschitz over its closed domain. Our first main contribution is a new tracking lemma showing that PRR approximates a differential inclusion over any finite time interval, without requiring smoothness or boundedness of the iterates. This tool allows us to establish the reachability of $(\epsilon, \delta)$-near approximate stationary points (NAS) for PRR under minimal assumptions and any nonsummable step size rule. Notably, this result applies even to problems with unbounded iterates, such as certain matrix completion instances. When the $f_i$ are weakly convex, we prove convergence to $(\epsilon, \delta)$-NAS using well-tuned summable step sizes. Furthermore, assuming uniqueness and convergence of the solution to the associated differential inclusion, we obtain convergence to $(\epsilon, 0)$-NAS under the same step size schedule. In the case where the $f_i$ are locally smooth, semialgebraic, and admit bounded subgradient trajectories, we establish convergence of PRR to $(0,0)$-NAS: with a nearly $o(1/\sqrt{k})$ rate for $N > 1$ using diminishing step sizes, and a faster $o(1/k)$ rate for $N = 1$ with constant step sizes. These results combine our tracking lemma with a new length formula for PRR. As an application, we show bounded subgradient trajectories for nonnegative matrix factorization (NMF), leading to the first convergence guarantee to stationary points for solving NMF via the projected gradient method with constant step sizes.

This work has several limitations. For nonsmooth $f_i$, in order to avoid assuming bounded iterates, our results are relatively weak: either proving only reachability, or requiring summable step sizes, which offer limited practical guidance. For smooth $f_i$, the convergence result for $N > 1$ relies on diminishing step sizes, which yield slower convergence, whereas constant step sizes often perform well in practice, yet our current analysis cannot capture this behavior. Lastly, our analysis is worst-case in nature and does not reflect the potential acceleration offered by random reshuffling. Addressing these limitations leads to several promising directions for future research.

\noindent\textbf{Acknowledgements} The authors thank Tonghua Tian for fruitful discussions, and the area editor, the associate editor, and the anonymous referees for their careful reading and constructive comments, which improved the paper significantly. 

\phantomsection
\renewcommand{\bibsection}{\section*{References}}
\addcontentsline{toc}{section}{\protect References}

\bibliographystyle{abbrvnat}
\bibliography{references}

\appendix
\renewcommand{\theexample}{\thesection.\arabic{example}}
\renewcommand{\thefigure}{\thesection.\arabic{figure}}

\section{Omitted proof}\label{sec:omittedproof}
\setcounter{example}{0}
\setcounter{figure}{0}

\subsection{Proof of \cref{thm:sgd_finite_len}}\label{app:thm4proof}
Let $X_0$ be a bounded subset of $\dom \Phi$. Our goal is to show that there exist $\bar{\alpha}>0$ such that $\Gamma(X_0,\bar{\alpha})<\infty$, where 
    \begin{subequations}\label{eq:sup_trajectory_dis_sgd}
    \begin{align}
        \Gamma(X_0,\bar{\alpha}) := & \sup\limits_{\substack{x\in (\mathbb{R}^n)^{\mathbb{N}\times \llbracket 0,N\rrbracket},\sigma\in\mathfrak{S}_{N}^{\mathbb{N}}\\ (\alpha,\gamma)\in(0,\bar{\alpha}]\times \mathbb{N}^*}} ~~ \sum_{k=0}^{\infty} \|x_{k+1,0}-x_{k,0}\| \\
      & ~~ \mathrm{s.t.} ~~~ 
    \left\{ 
    \begin{array}{l}
    x_{k+1,0}=\prox_{\alpha_k g}(x_{k,N}),\\[1mm]
    x_{k,i} = x_{k,i-1}  - \alpha_k \nabla f_{\sigma^k_i} (x_{k,i-1}), ~ \forall i \in \llbracket 1,N \rrbracket, \\
    \alpha_k = \alpha/(k+\gamma)^\beta,~\forall k\in\mathbb{N}, ~ x_{0,0} \in X_0, 
    \end{array}
    \right.
    \end{align}
    \end{subequations}
    and $\mathfrak{S}_{N}$ denotes the symmetric group of degree $N$. 

    Let $\Sigma:\mathbb{R}_+\times\dom g\rightarrow\dom g$ be defined by $\Sigma(t,x_0):=x(t)$ where $x(\cdot)$ is the unique solution of
    \begin{equation*}
        x'(t) \in - \partial \Phi(x(t)) = -\nabla f(x(t)) - \partial g(x(t)) = -\sum_{i=1}^N\nabla f_i(x(t)) - \partial g(x(t)), \quad x(0)=x_0. 
    \end{equation*}
    Notice that the uniqueness can be easily derived from \citet[Theorem 4.3(b)]{cornet1983existence}. Let $C$ be the set of critical points of $\Phi$ in $\overline{\Sigma(\mathbb{R}_+,X_0)}$. Note that $C$ is bounded due to \citet[Proposition 5.3]{li2025bounded} (see also \citet[Lemma 1]{josz2023global}) and closed due to \citet[Proposition 2(a), p. 141]{aubin1984differential}. Thus, $C$ is compact, and either $X_0\subset C$ (in this case just take any $\epsilon\in(0,1]$) or there exists $\epsilon\in(0,1]$ such that $X_0\setminus\mathring{B}(C,\epsilon/6)\ne\emptyset$. 

    To prove $\Gamma(X_0,\bar{\alpha})<\infty$ for some $\bar{\alpha}>0$, we follow the two steps below. First, we will show that the iterates can go arbitrarily close to a critical point of $\Phi$. Second, we will show that either the iterates always stay within a neighborhood of the critical point, or they admit a sufficient decrease in function value when they step out of the neighborhood of the critical point. 

    For the first step, we show that for every $X_0$ and every $\delta\in(0,\epsilon/2)$, there exists $\hat{\alpha}$ such that for any $(x_k)_{k\in\mathbb{N}}$ generated by \cref{alg:prr} with $x_0\in X_0$ and $(\alpha_k)_{k\in\mathbb{N}}$ satisfying $\alpha_k\in(0,\hat{\alpha}]$ for all $k\in\mathbb{N}$, there exist $k^*\in\mathbb{N}$ and $x^*\in C$ such that $\|x_{k^*}-x^*\|\leqslant \delta$. If $X_0\subset C$, then the desired result trivially holds because we can simply take any $\hat{\alpha}>0$ with $k^*=0$ and $x^*=x_0$. If $X_0\setminus\mathring{B}(C,\epsilon/6)\ne\emptyset$, then the argument can be further divided into two parts. The first part is to show that the subgradient trajectory of $\Phi$ can go arbitrarily close to a critical point of $\Phi$, i.e., we show that for fixed $X_0$ and $\delta$, there exists $T>0$ such that for any solution $x(\cdot)$ of $x'(t)\in -\partial\Phi(x(t))$ with $x_0\in X_0$, we have that $\|x(t^*)-x^*\|\leqslant\delta/3$ for some $t^*\in [0,T]$ and $x^*\in C$. To prove it, define $\Sigma_0:=\Sigma(\mathbb{R}_+,X_0)$ and let 
    \begin{equation*}
        \nu:=\inf_{x\in \overline{\Sigma_0}\setminus \mathring{B}(C,\delta/3)}d(0,\partial\Phi(x)) = \inf \left\{\|y\|:y\in \partial \Phi\left(\overline{\Sigma_0}\setminus \mathring{B}(C,\delta/3)\right)\right\}. 
    \end{equation*}
    Recall that $\partial \Phi:= \nabla f+\partial g$ has closed graph by \cref{lemma:closed}. Since $\overline{\Sigma_0}\setminus \mathring{B}(C,\delta/3)$ is a nonempty compact subset of $\dom g$, the set $\partial\Phi(\overline{\Sigma_0}\setminus \mathring{B}(C,\delta/3))$ is nonempty and closed. Thus, $\nu$ is finite and attained. Combined with the definition of $C$, it holds that $\nu>0$. Let $T:=2\Gamma(X_0)/\nu$ where $\Gamma(X_0)$ is defined by
    \begin{align*}
        \Gamma(X_0) := & \sup\limits_{x \in \mathcal{A}(\mathbb{R}_+,\mathbb{R}^n)} ~~ \int_0^{\infty} \|x'(t)\|dt \\
      & ~~ \mathrm{subject~to} ~~~ 
    \left\{ 
    \begin{array}{l}
        x'(t) \in - \partial \Phi(x(t)), ~\mathrm{for~a.e.}~ t > 0,\\[3mm] x(0) \in X_0.
        \end{array}
        \right.
    \end{align*}
    and $\Gamma(X_0)<\infty$ by \citet[Proposition 5.3]{li2025bounded} (see also \citet[Lemma 1]{josz2023global}). Note that $T>0$ because $X_0\not\subset C$ and $\Gamma(X_0)>0$.  Furthermore, it follows that there exists $t^*\in[0,T]$ such that $d(0,\partial\Phi(x(t^*))) =\|x'(t^*)\|< 2\Gamma(X_0)/T$. We know from \cref{lem:chain} that $d(0,\partial\Phi(x(t))) =\|x'(t)\|$ for almost every $t\in [0,T]$.  Assume the contrary that $\|x'(t)\|\geqslant 2\Gamma(X_0)/T$ for all such $t$, then 
    \begin{equation*}
        \int_0^\infty\|x'(t)\|dt \geqslant \int_0^T\|x'(t)\|dt \geqslant 2\Gamma(X_0) > \Gamma(X_0)
    \end{equation*}
    contradicts the definition of $\Gamma(X_0)$. 
    According to the definition of $\nu$, we have $x(t^*)\in \mathring{B}(C,\delta/3)$. Hence, $\|x(t^*)-x^*\|\leqslant\delta/3$ for this $x^*\in C$ and the claim follows. 
    
    We next show that there exists $\hat{\alpha}>0$ such that for any $(x_k)_{k\in\mathbb{N}}$ generated by \cref{alg:prr} with $x_0\in X_0$ and $(\alpha_k)_{k\in\mathbb{N}}$ satisfying $\alpha_k\in(0,\hat{\alpha}]$ for all $k\in\mathbb{N}$, there exists $k^*\in\mathbb{N}$ such that $\|x_{k^*}-x(t^*)\|\leqslant 2\delta/3$. To see this, let 
    \begin{equation*}
        L_\Phi:=\sup_{x\in B(\overline{\Sigma_0},\epsilon)}\|\nabla f(x)\|+\sup_{x\in B(\overline{\Sigma_0},\epsilon)\cap\dom g}d(0,\partial g(x)),   
    \end{equation*}
  where $L_\Phi<\infty$ by \cref{lemma:min_grad}.  By \cref{lemma:stoc_tracking}, there exists $\hat{\alpha}\in(0,\delta/(3L_\Phi)]$ such that for any $(x_k)_{k\in\mathbb{N}}$ generated by \cref{alg:prr} with $x_0\in X_0$ and $(\alpha_k)_{k\in\mathbb{N}}$ satisfying $\alpha_k\in(0,\hat{\alpha}]$ for all $k\in\mathbb{N}$, we have that for all $k\in\mathbb{N}^*$, 
    \begin{equation*}
        \sum_{i=0}^{k-1}\alpha_i \leqslant T \implies \left\|x_k-x\left(\sum_{i=0}^{k-1}\alpha_i\right)\right\|\leqslant \frac{\delta}{3}, 
    \end{equation*}
  Since $\alpha_k\leqslant \hat{\alpha}\leqslant \delta/(3L_\Phi)$ and $\sum_{k = 0}^\infty \alpha_k = \infty$, there exists $k^*\in\mathbb{N}$ such that $t_{k^*}:=\sum_{k=0}^{k^*-1}\alpha_k\in[t^*-\delta/(3L_\Phi),t^*]$, where the convention $\sum_{k=0}^{-1}\alpha_k=0$ is used. Thus, 
    \begin{align*}
        \|x_{k^*}-x(t^*)\| &\leqslant \|x_{k^*}-x(t_{k^*})\| + \|x(t_{k^*})-x(t^*)\| \leqslant \frac{\delta}{3} + \int_{t_{k^*}}^{t^*}\|x'(\tau)\|\;d\tau \\
        &= \frac{\delta}{3} + \int_{t_{k^*}}^{t^*}d(0,\partial\Phi(x(\tau)))\;d\tau \leqslant \frac{\delta}{3}+L_\Phi(t^*-t_{k^*}) \leqslant 
        \frac{2\delta}{3}. 
    \end{align*}
    The claim of the first step of the proof follows easily from combining the above two parts, i.e., 
    \begin{equation*}
        \|x_{k^*}-x^*\|\leqslant \|x_{k^*}-x(t^*)\|+\|x(t^*)-x^*\|\leqslant \frac{2\delta}{3}+\frac{\delta}{3} = \delta. 
    \end{equation*}

    For the second step, define $K^*:=\inf\{k\geqslant k^*:x_k\not\in\mathring{B}(C,\epsilon)\}$. Clearly, we have that $K^*\geqslant k^*+1$. Furthermore, let $L_g:=\sup_{x\in B(\overline{\Sigma_0},\epsilon)\cap \dom g}d(0,\partial g(x))$ and $L:=\max\{L_1,\ldots,L_N\}$, where $L_i:=\sup_{x\in B(\overline{\Sigma_0},\epsilon+1)}\|\nabla f_i(x)\|$ for $i=1,\ldots,N$. Since $\Phi$ is continuous over $\dom g$ and $C$ is compact, there exists $\delta\in(0,\epsilon/2)$ such that 
    \begin{equation}\label{eq:sgd_main_continuity}
        \Phi(x) - \max_C \Phi \leqslant \frac{1}{4}\psi^{-1}\left(\frac{\epsilon}{3}\right), \quad \forall\,x\in B(C,\delta)\cap\dom g. 
    \end{equation}
    Applying the result in the first step to this $\delta$, we obtain the desired $\hat{\alpha}\in(0,\delta/(3L_\Phi)]$, $k^*$ and $x^*$. Since $B(\overline{\Sigma_0},\epsilon)$ is bounded, applying \cref{cor:length} to $B(\overline{\Sigma_0},\epsilon)\cap \dom \Phi$, there exist $\eta,\kappa>0$ and $\tilde{\alpha}\in(0,\hat{\alpha}]$ such that for every $\alpha\in (0,\tilde{\alpha}]$, 
    \begin{equation}\label{eq:sgd_main_discr_len}
        \sum_{k=0}^K\|x_{k+1}-x_k\| \leqslant \psi\left(\Phi(x_0)-\Phi(x_K)+ \eta \alpha\right) + \kappa\alpha. 
    \end{equation}
    This inequality holds for all $K=0,\ldots,K^*-1$ and for any $(x_0,\ldots,x_{K^*})\in(B(\overline{\Sigma_0},\epsilon))^{K^*}\times\mathbb{R}^n$ generated by \cref{alg:prr} with $\alpha
    _k=\alpha/(k+\gamma)^\beta$ where  $\alpha\in(0,\tilde{\alpha}]$. 
    
    If $K^*=\infty$, then $x_k\in \mathring{B}(C,\epsilon)$ for all $k\geqslant k^*$. Combined with the fact that $x_k\in B(\overline{\Sigma_0},\delta/3)$ for all $k=0,\ldots,k^*-1$, we have that $x_k\in B(\overline{\Sigma_0},\epsilon)$ for $k\in\mathbb{N}$. In this case, a limiting argument shows 
    \begin{align*}
        \sum_{k=0}^\infty\|x_{k+1}-x_k\| \leqslant \psi\left(\sup_{X_0}\Phi-\inf_{B(\overline{\Sigma_0},\epsilon)}\Phi+ \eta \alpha\right) + \kappa\alpha. 
    \end{align*}
    Thus, $\Gamma(X_0,\tilde{\alpha})<\infty$ immediately follows in this case. 
    
    If $K^*<\infty$, then consider any  $\alpha \in (0,\bar{\alpha}_0]$, where
    \begin{equation*}
        \bar{\alpha}_0 := \min\left\{\tilde{\alpha},\frac{\epsilon}{4(NL+L_g)},\frac{\epsilon}{6\kappa},\frac{1}{2\eta}\psi^{-1}\left(\frac{\epsilon}{3}\right)\right\}.
    \end{equation*}
    Furthermore, $K^*\geqslant k^*+2$ because \cref{cor:Oalpha} implies that 
    \begin{align*}
        \|x_{k^*+1}-x^*\| &\leqslant \|x_{k^*+1}-x_{k^*}\|+\|x_{k^*}-x^*\| \leqslant 2\alpha_{k^*}(NL+L_g)+\delta \\
        &\leqslant 2\frac{\epsilon}{4(NL+L_g)}(NL+L_g)+\delta < \frac{\epsilon}{2}+\frac{\epsilon}{2} = \epsilon. 
    \end{align*}
    We claim that at iteration $K^*-1$, the objective value has admitted a sufficient decrease from the largest critical value in $C$, i.e., 
    \begin{equation}\label{eq:sgd_main_suffic_decre}
        \Phi(x_{K^*-1}) \leqslant \max_C \Phi - \frac{1}{4}\psi^{-1}\left(\frac{\epsilon}{3}\right). 
    \end{equation}
    To verify the claim, first notice that  
    \begin{equation*}
        \sum_{k=k^*}^{K^*-1}\|x_{k+1}-x_k\| \geqslant \|x_{K^*}-x_{k^*}\| \geqslant \|x_{K^*}-x^*\|-\|x_{k^*}-x^*\| \geqslant \epsilon-\delta \geqslant \frac{\epsilon}{2}.  
    \end{equation*}
    Furthermore, using the range of $\alpha$ and the fact that $\psi^{-1}$ is increasing, we have 
    \begin{equation*}
        \psi^{-1}\left(\sum_{k=k^*}^{K^*-1}\|x_{k+1}-x_k\| - \kappa\alpha\right) \geqslant  \psi^{-1}\left(\frac{\epsilon}{2}-\frac{\epsilon}{6}\right) = \psi^{-1}\left(\frac{\epsilon}{3}\right). 
    \end{equation*}
    On the other hand, \eqref{eq:sgd_main_discr_len} implies that 
    \begin{equation*}
        \Phi(x_{k^*})-\Phi(x_{K^*-1}) \geqslant \psi^{-1}\left(\sum_{k=k^*}^{K^*-1}\|x_{k+1}-x_k\| - \kappa\alpha\right) - \eta\alpha 
    \end{equation*}
    Therefore, it yields that 
    \begin{equation}\label{eq:lenfor_decrease}
        \Phi(x_{k^*})-\Phi(x_{K^*-1}) \geqslant \psi^{-1}\left(\frac{\epsilon}{3}\right) - \eta\alpha \geqslant \frac{1}{2}\psi^{-1}\left(\frac{\epsilon}{3}\right). 
    \end{equation}
    Finally, since $x_{k^*}\in B(C,\delta)\cap\dom g$, combining \eqref{eq:sgd_main_continuity} and \eqref{eq:lenfor_decrease}, we obtain the desired claim \eqref{eq:sgd_main_suffic_decre}. 

    Define the successive initial set of $X_0$ by 
    \begin{equation}\label{eq:X1}
        X_1 := B(C,\epsilon) \bigcap \left\{ x \in\dom g : \Phi(x) \leqslant \max_C \Phi - \frac{1}{4}\psi^{-1}\left(\frac{\epsilon}{3}\right) \right\}.
    \end{equation}
    From \eqref{eq:sgd_main_suffic_decre}, we know that $x_{K^*-1}\in X_1$. By \eqref{eq:sgd_main_discr_len}, for any feasible points $(x,\sigma,\alpha,\gamma)$ of \eqref{eq:sup_trajectory_dis_sgd}, we have 
    \begin{align*}
        \sum_{k=0}^\infty \|x_{k+1}-x_k\| &= \sum_{k=0}^{K^*-2} \|x_{k+1}-x_k\| + \sum_{k=K^*-1}^\infty \|x_{k+1}-x_k\| \\
        &\leqslant \psi\left(\Phi(x_0)-\Phi(x_{K^*-2})+ \eta \alpha\right) + \kappa\alpha + \sum_{k=K^*-1}^\infty \|x_{k+1}-x_k\| \\
        &\leqslant \psi\left(\sup_{X_0}\Phi-\inf_{B(\overline{\Sigma_0},\epsilon)}\Phi+\eta\bar{\alpha}_0\right) + \kappa\bar{\alpha}_0 + \Gamma(X_1,\bar{\alpha}_0), 
    \end{align*}
    where in the last inequality we use the fact that 
    if $(x,\sigma,\alpha,\gamma)$ is feasible for \eqref{eq:sup_trajectory_dis_sgd}, then for any $\bar{k}\in \mathbb{N}$, $((x_{k,\cdot})_{k\geqslant \bar{k}},(\sigma_k)_{k\geqslant \bar{k}},\alpha,\bar{k}+\gamma)$ is again feasible. Using the convention $\Gamma(\emptyset,\cdot)=-\infty$, by taking supremum over all feasible points $(x,\sigma,\alpha,\gamma)$ of \eqref{eq:sup_trajectory_dis_sgd}, we can combine the case when $K^*=\infty$ and $K^*<\infty$ into 
    \begin{equation*}
        \Gamma(X_0,\bar{\alpha}_0) \leqslant \psi\left(\sup_{X_0}\Phi-\inf \Phi+\eta\bar{\alpha}_0\right) + \kappa\bar{\alpha}_0+\max\{\Gamma(X_1,\bar{\alpha}_0),0\}. 
    \end{equation*}
    Note that the above inequality is also valid if $\bar{\alpha}_0$ is replaced by any $\alpha\in (0,\bar{\alpha}_0]$. 
    
    If $X_1\neq \emptyset$, we may repeat the above arguments by replacing $X_0$ with $X_1$. Recursively, we obtain sequences of $X_\ell, \bar{\alpha}_\ell, \psi_\ell, \eta_\ell,\kappa_\ell$ such that
    \begin{equation}\label{eq:main_recur}
        \Gamma(X_\ell,\bar{\alpha}_\ell) \leqslant \psi_\ell\left(\sup_{X_\ell}\Phi-\inf \Phi+\eta_\ell\bar{\alpha}_\ell\right) + \kappa_\ell\bar{\alpha}_\ell+\max\{\Gamma(X_{\ell+1},\bar{\alpha}_\ell),0\},
    \end{equation}
    for any $\ell \in \mathbb{N}$ as long as $X_\ell \neq \emptyset$, where  $\psi_0 := \psi, \eta_0:=\eta, \kappa_0 = \kappa$. Again, the above inequality remains valid when reducing $\bar{\alpha}_\ell$. 
    
    We finally show that $X_\ell = \emptyset$ for some $\ell\in \mathbb{N}$. Let $v_\ell$ be the maximum critical value of $\Phi$ in $(-\infty, \max_{X_\ell} \Phi]$. By the construction of the sequence of initial sets \eqref{eq:X1}, it is evident that if neither $X_\ell$ nor $X_{\ell + 1}$ is empty, then $v_{\ell} > \max_{X_{\ell+1}} \Phi \geqslant v_{\ell+1}$. As $\Phi$ has finitely many critical values by \citet[Corollary 9]{bolte2007clarke}, eventually $X_\ell = \emptyset$. Let $\bar{\ell} \geqslant 1$ be the smallest index such that $X_{\bar{\ell}} = \emptyset$. Telescoping \eqref{eq:main_recur} with $\bar{\alpha}_\ell$ replaced by $\bar{\alpha}:= \min\{\bar{\alpha}_\ell:\ell = 0,\ldots,\bar{\ell}-1\}$ yields
    \begin{equation*}
        \Gamma(X_0,\bar{\alpha}) \leqslant \sum_{\ell = 0}^{\bar{\ell} - 1} \psi_\ell\left(\sup_{X_\ell}\Phi-\inf \Phi+\eta_\ell\bar{\alpha}_\ell\right)+ \kappa_\ell\bar{\alpha}_\ell<\infty.
    \end{equation*}

The finiteness of $\Gamma(X_0,\bar{\alpha})$ naturally means that $\|x_k-x_0\| \leqslant \sum_{i=0}^{k-1}\|x_{i+1}-x_i\| \leqslant \Gamma(X_0,\bar{\alpha})$ for all $k\in\mathbb{N}$, that is, $x_k \in B(X_0,\Gamma(X_0,\bar{\alpha}))$ (where $x_k := x_{k,0}$). Let $L_0,M$ respectively be Lipschitz constants of $f_1,\hdots,f_N$ and $\nabla f_1,\hdots,\nabla f_N$ on $B(X_0,\Gamma(X_0,\bar{\alpha})+2)$. After possibly reducing $\bar{\alpha}$, we have $\bar{\alpha} \leqslant 1/(NL_0)$. Since $x_{k+1} = \prox_{\alpha_k g} (x_{k,N})$, Fermat's rule (\citet[Theorem 10.1]{rockafellar2009variational}) implies that $0 \in \alpha_k \partial g(x_{k+1}) + x_{k+1} - x_{k,N}$ and hence
\begin{equation*}
    \partial \Phi(x_{k+1}) = \nabla f(x_{k+1}) + \partial g(x_{k+1}) \ni \nabla f(x_{k+1}) - \frac{x_{k+1}-x_{k,N}}{\alpha_k}.
\end{equation*}
Since $x_{k,N} = x_k - \alpha_k \sum_{i=1}^N\nabla f_{\sigma_i^k}(x_{k,i-1})$, by \cref{cor:Oalpha} it follows that
    \begin{align*}
        d(0,\partial \Phi(x_{k+1})) & \leqslant \left\|\nabla f(x_{k+1}) - \frac{x_{k+1}-x_{k,N}}{\alpha_k}\right\| \\
        & \leqslant \frac{\|x_{k+1}-x_k\|}{\alpha_k} + \left\|\nabla f(x_{k+1}) - \sum_{i=1}^N\nabla f_{\sigma_i^k}(x_{k,i-1})\right\| \\
        & \leqslant \frac{\|x_{k+1}-x_k\|}{\alpha_k} + \|\nabla f(x_{k+1}) - \nabla f(x_k) \| + \left\| \nabla f(x_k) - \sum_{i=1}^N\nabla f_{\sigma_i^k}(x_{k,i-1})\right\| \\
        & \leqslant \left(\frac{1}{\alpha_k}+M\right)\|x_{k+1}-x_k\| + (N-1)\sqrt{\frac{N(2N-1)}{6}}ML_0 \alpha_k.
    \end{align*}
    The last inequality admits a similar derivation as in the proof of \cref{lem:approx_descent}. Note that for every nonnegative decreasing sequence $u_0,u_1,u_2,\hdots$, we have $u_k \leqslant 2 (\sum_{i=\lfloor k/2\rfloor}^\infty u_i)/(k+2)$ as explained in \citet[Footnote 1]{josz2023global}. In particular, for $u_k := \min_{i\in \llbracket 0,k \rrbracket} \alpha_i d(0,\partial \Phi(x_{i+1}))$, we have
    \begin{align*}
        \min_{i\in \llbracket 0,k \rrbracket} \alpha_i d(0,\partial \Phi(x_{i+1})) & \leqslant \frac{2}{k+2} \sum_{i=\lfloor k/2\rfloor}^\infty \min_{j \in \llbracket 0,i\rrbracket} \alpha_j d(0,\partial \Phi(x_{j+1})) \\
        & \leqslant \frac{2}{k+2} \sum_{i=\lfloor k/2\rfloor}^\infty \alpha_i d(0,\partial \Phi(x_{i+1})) \\
        & \leqslant \frac{2}{k+2} \sum_{i=\lfloor k/2\rfloor}^\infty (1+\alpha M)\|x_{i+1}-x_i\| +(N-1)\sqrt{\frac{N(2N-1)}{6}}ML_0\alpha_i^2.
    \end{align*}
    Thus, using $\alpha_k = \alpha/(k+1)^\beta$, we have
    \begin{small}
    \begin{align*}
        \min_{i\in \llbracket 0,k \rrbracket} d(0,\partial \Phi(x_{i+1})) 
        \leqslant & ~  \frac{2}{\alpha_k(k+2)} \sum_{i=\lfloor k/2\rfloor}^\infty (1+\alpha M) \|x_{i+1}-x_i\| 
        + (N-1)\sqrt{\frac{N(2N-1)}{6}}ML_0 \alpha_i^2 \\
        \leqslant & ~  \frac{2\alpha^{-1}}{(k+1)^{1-\beta}} \sum_{i=\lfloor k/2\rfloor}^\infty (1+\alpha M) \|x_{i+1}-x_i\| 
        + (N-1)\sqrt{\frac{N(2N-1)}{6}}ML_0 \frac{\alpha^2}{(i+1)^{2\beta}} \\
        \leqslant &  ~ \frac{2\alpha^{-1}}{(k+1)^{1-\beta}} \left((1+\alpha M)\sum_{i=\lfloor k/2\rfloor}^\infty \|x_{i+1}-x_i\|\right. \\
        &\left.+ (N-1)\sqrt{\frac{N(2N-1)}{6}}ML_0\left(\frac{\alpha^2}{(\lfloor k/2\rfloor+1)^{2\beta}}+\sum_{i=\lfloor k/2\rfloor+1}^\infty\int_i^{i+1}\frac{\alpha^2}{v^{2\beta}}\;dv\right)\right) \\
        \leqslant &  ~ \frac{2\alpha^{-1}}{(k+1)^{1-\beta}} \left((1+\alpha M)\sum_{i=\lfloor k/2\rfloor}^\infty \|x_{i+1}-x_i\|\right. \\
        &\left.+ (N-1)\sqrt{\frac{N(2N-1)}{6}}ML_0\left(\frac{\alpha^2}{(\lfloor k/2\rfloor+1)^{2\beta}}+\int_{\lfloor k/2\rfloor+1}^{\infty}\frac{\alpha^2}{v^{2\beta}}\;dv\right)\right) \\
        \leqslant &  ~ \frac{1}{(k+1)^{1-\beta}} \left(2(\alpha^{-1}+M)\sum_{i=\lfloor k/2\rfloor}^\infty \|x_{i+1}-x_i\|\right. \\
        &\left.+ (N-1)\sqrt{\frac{N(2N-1)}{6}}ML_0\left(\frac{2\alpha}{(\lfloor k/2\rfloor+1)^{2\beta}}+\frac{2\alpha}{(2\beta-1)(\lfloor k/2\rfloor+1)^{2\beta-1}}\right)\right) \\
        \leqslant &  ~ \frac{1}{(k+1)^{1-\beta}} \left(2(\alpha^{-1}+M)\sum_{i=\lfloor k/2\rfloor}^\infty \|x_{i+1}-x_i\| + \frac{4(N-1)\sqrt{N(2N-1)}ML_0\alpha}{\sqrt{6}(2\beta-1)(\lfloor k/2\rfloor+1)^{2\beta-1}}\right). 
    \end{align*}
    \end{small}
\hfill $\square$

\subsection{Proof of \cref{prop:conv}}\label{app:propconvproof}

For simplicity, we can consider the setting $U\subset B(0,1)$. By assumption, there exist $r_x,r_y>0$ such that for every $k\in\mathbb{N}$, $x_k(t)\in B(0,r_x)$ and $y_k(t)\in B(0,r_y)$ for all $t\in I\setminus Z_k$, where $Z_k$ has zero measure. Let $Z:=\bigcup_{k\in\mathbb{N}}Z_k$; then $Z$ also has zero measure. Thus, $x_k(t)\in B(0,r_x)$ and $y_k(t)\in B(0,r_y)$ for all $k\in\mathbb{N}$ and $t\in I\setminus Z$. In addition, there exists $r_F>0$ such that $F(x_k(t))\subset F(B(0,r_x))\subset B(0,r_F)$ for all $k\in\mathbb{N}$ and $t\in I\setminus Z$ as $F$ is locally bounded. 

Furthermore, there is a set $Z'$ with zero measure such that for all $t\in I\setminus Z'$, for any neighborhood $U$ of $0$ such that $U\subset B(0,1)$, there exists $k_0\in\mathbb{N}$ satisfying 
\begin{equation*}
    (x_k(t),y_k(t)) \in \gra(F+G)+U, \quad \forall\,k\geqslant k_0.
\end{equation*}
Since $y_k(t) \in F(x_k(t))+G(x_k(t)) + B(0,1)$, there exist $\alpha_k(t)\in F(x_k(t))$, $\beta_k(t)\in G(x_k(t))$ and $\gamma_k(t)\in B(0,1)$ such that $y_k(t)=\alpha_k(t)+\beta_k(t)+\gamma_k(t)$ for all $k\geqslant k_0$ and $t\in I\setminus (Z\cup Z')$. Notice that $\|\beta_k(t)\|\leqslant r_y+r_F+1$. Set $r:=r_y+r_F+1$, and define the cut set-valued mapping $G_r$ as $G_r(\cdot) := G(\cdot)\cap B(0,r)$. Note that $F+G_r$ is proper, locally bounded, and has convex values. In addition, $F+G_r$ has closed graph by \cref{lemma:closed}. Using similar arguments as in \citet[2.1.5(d) Proposition]{clarke1990}, it follows that $F+G_r$ is upper semicontinuous (\citet[p. 59]{aubin1984differential}), and thus upper hemicontinuous (\citet[p. 60]{aubin1984differential}). Since $\beta_k(t)\in G_r(x_k(t))$, for all $t\in I\setminus (Z\cup Z')$, we have 
\begin{equation*}
    (x_k(t),y_k(t)) \in \gra(F+G_r)+U, \quad \forall\,k\geqslant k_0.
\end{equation*}
It is easy to see the above also holds for neighborhoods not necessarily subsets of $B(0,1)$ by taking the same $k_0$ as the one for neighborhood $U\cap \mathring{B}(0,1)$. Therefore, applying \citet[Theorem 1, p. 60]{aubin1984differential} to $F+G_r$, we can conclude that for a.e. $t\in I$, 
\begin{equation*}
    (x(t),y(t))\in \gra(F+G_r) \subset \gra(F+G). 
\end{equation*}
\hfill $\square$

\section{Supporting figures}\label{sec:suppfigs}
\label{app:examples}
\setcounter{example}{0}
\setcounter{figure}{0}

This section contains supporting figures for numerical experiments mentioned in \cref{sec:Introduction}. 

\begin{figure}[!ht]
\centering 
\includegraphics[width=0.48\textwidth]{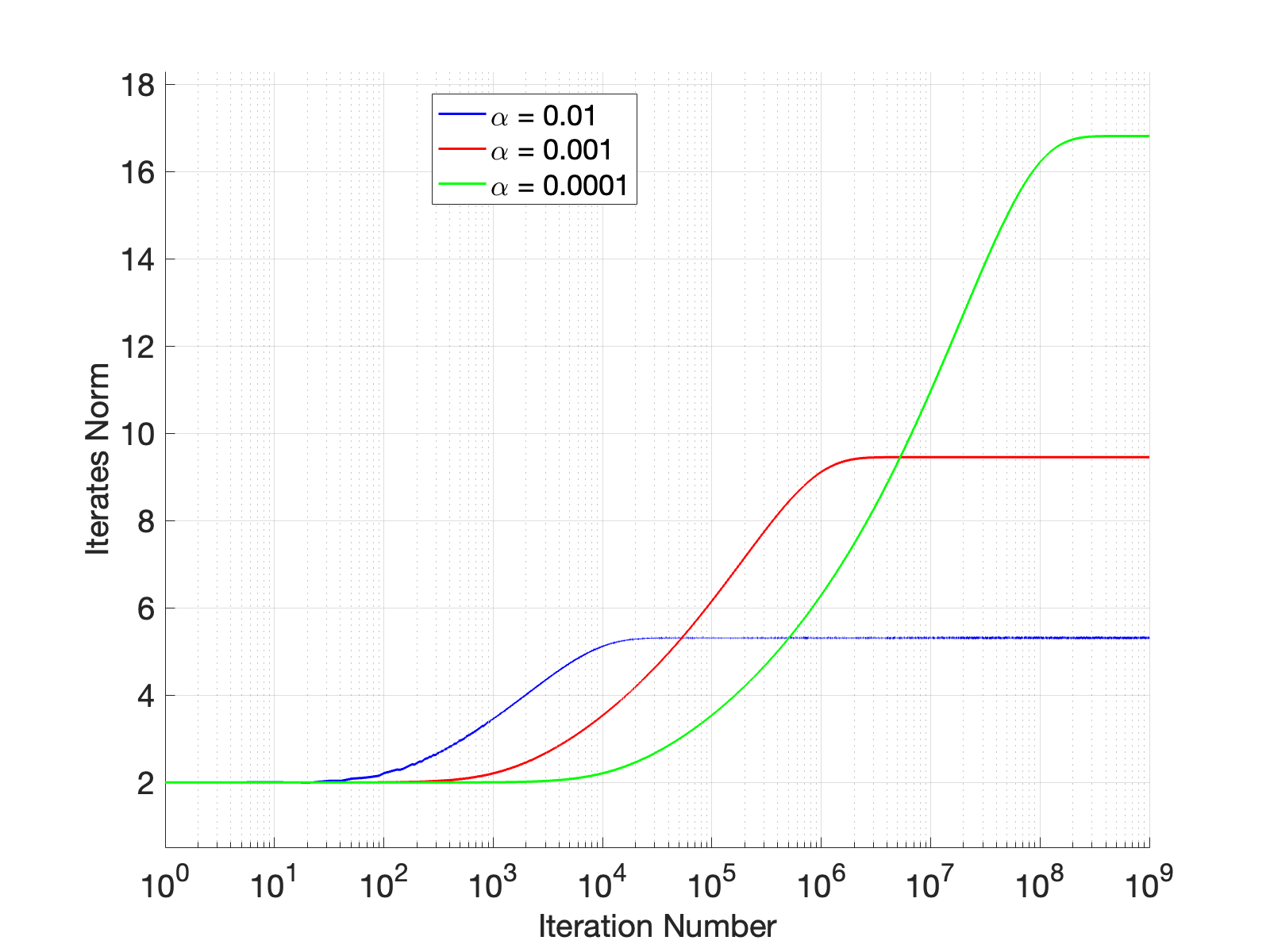}\hspace{-1ex}
\includegraphics[width=0.48\textwidth]{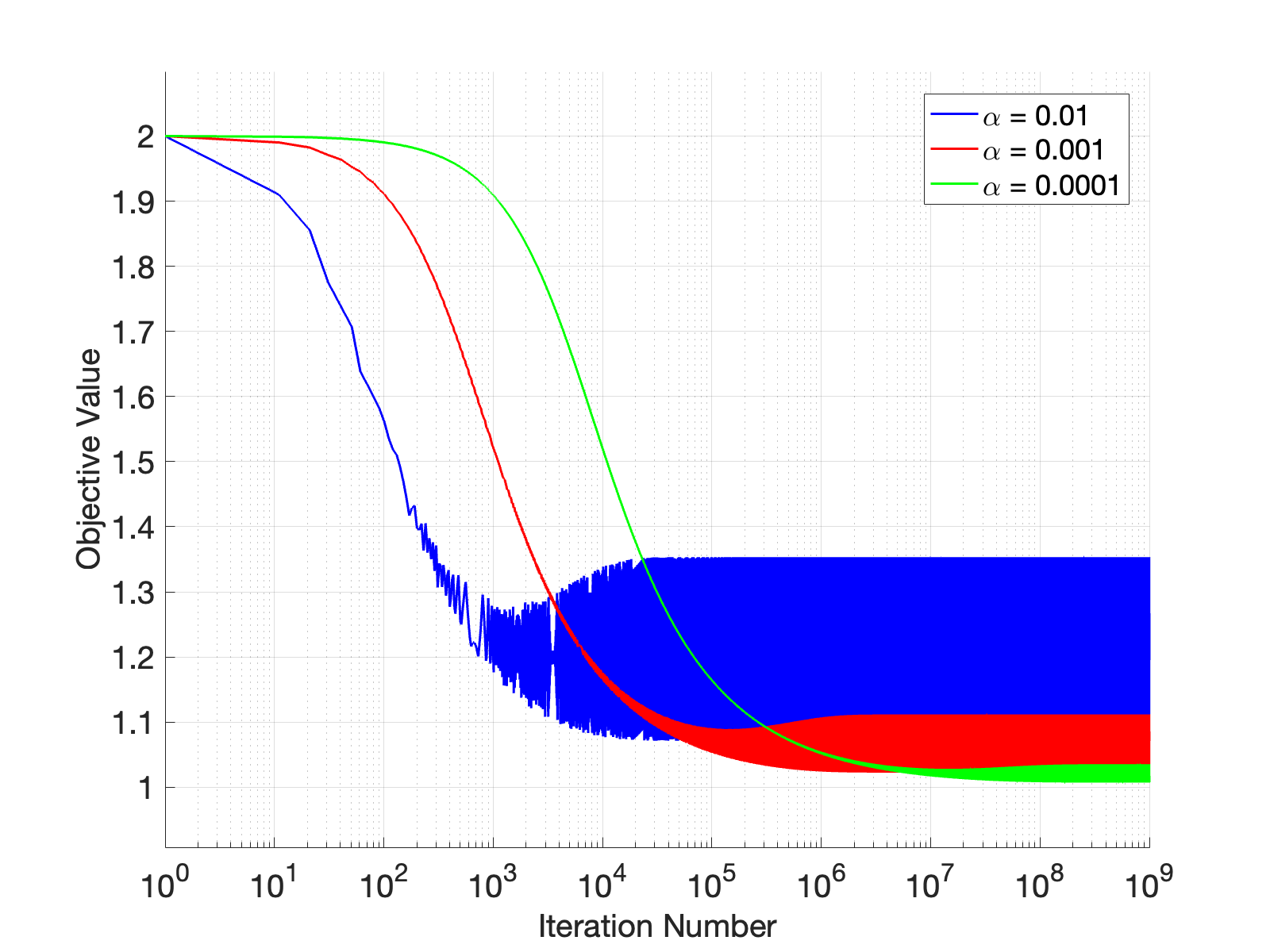}\hspace{-1ex}
\caption{Iterate norms and function values of the subgradient method with constant step sizes $\alpha_k=\alpha$ applied to $\ell_1$ matrix completion instance $f(x_1,x_2,y_1,y_2):=|1-x_1y_1|+|1-x_2y_1|+|1-x_2y_2|$.}
\label{fig:constant}
\end{figure}

\begin{figure}[!ht]
\centering 
\includegraphics[width=0.48\textwidth]{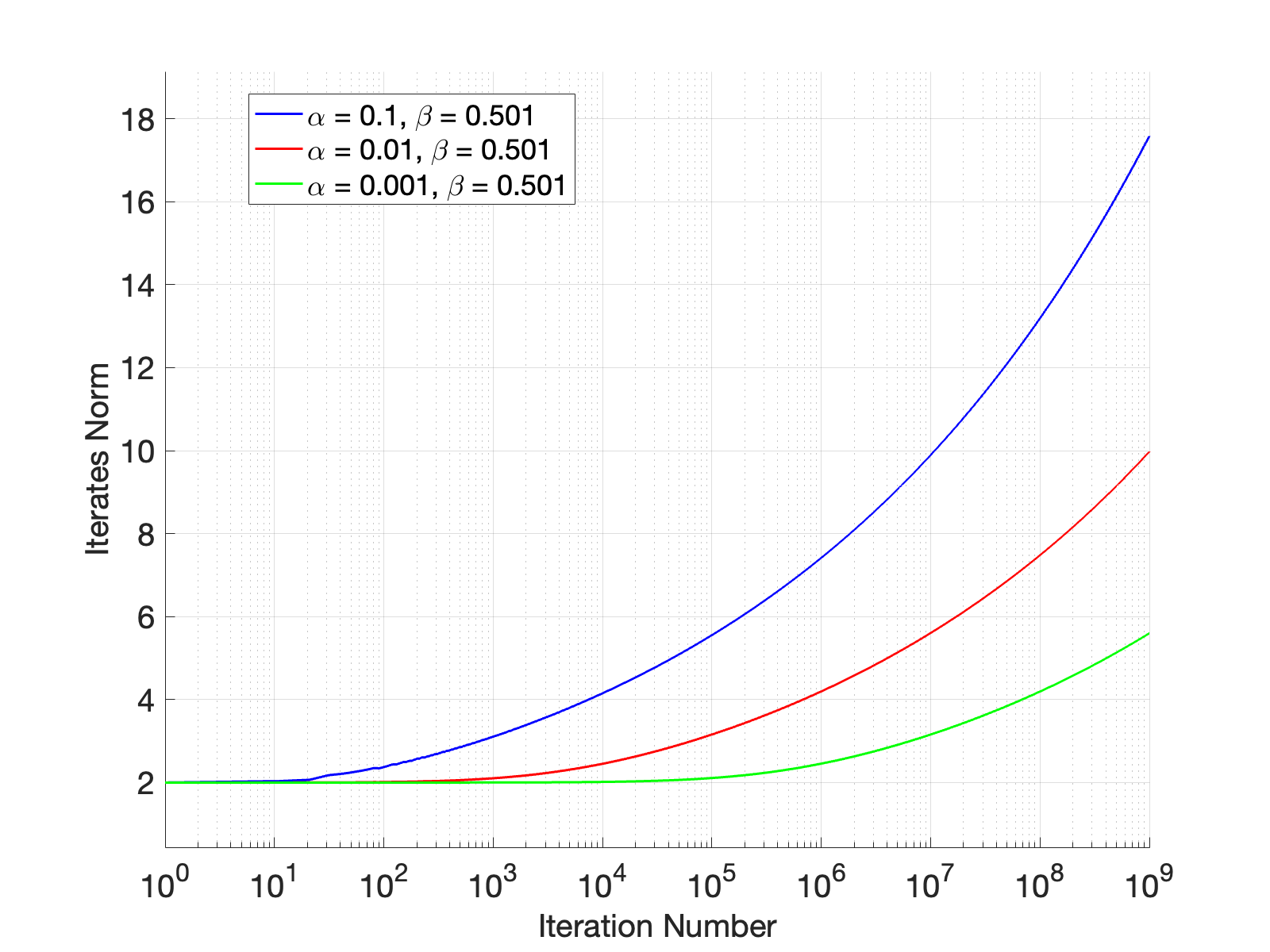}\hspace{-1ex}
\includegraphics[width=0.48\textwidth]{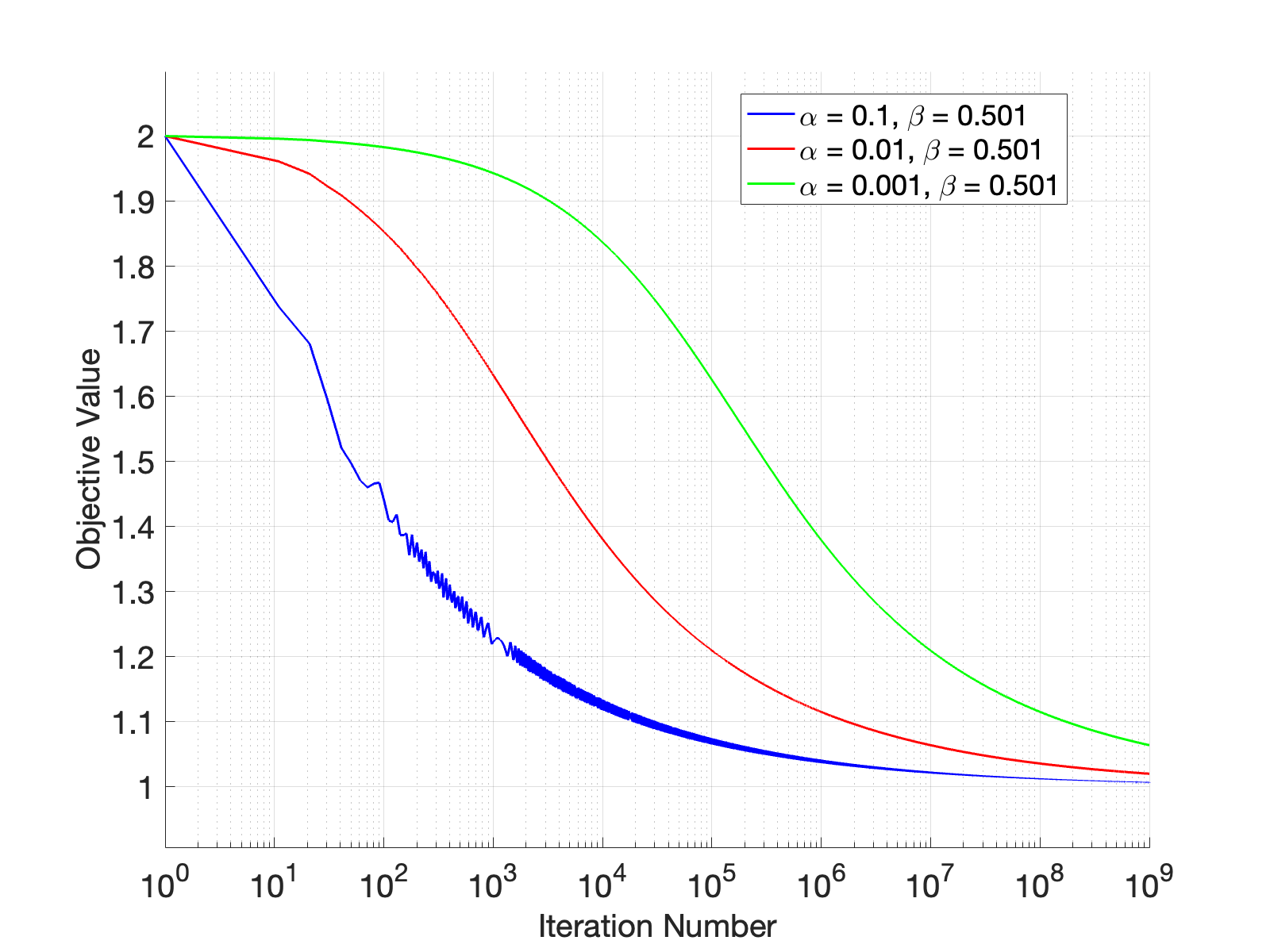}\hspace{-1ex}
\caption{Iterate norms and function values of the subgradient method with non-summable diminishing step sizes $\alpha_k=\alpha/(k+1)^\beta$ applied to $\ell_1$ matrix completion instance.}
\label{fig:unbdd_iter}
\end{figure}

\begin{figure}[!ht]
\centering 
\includegraphics[width=0.48\textwidth]{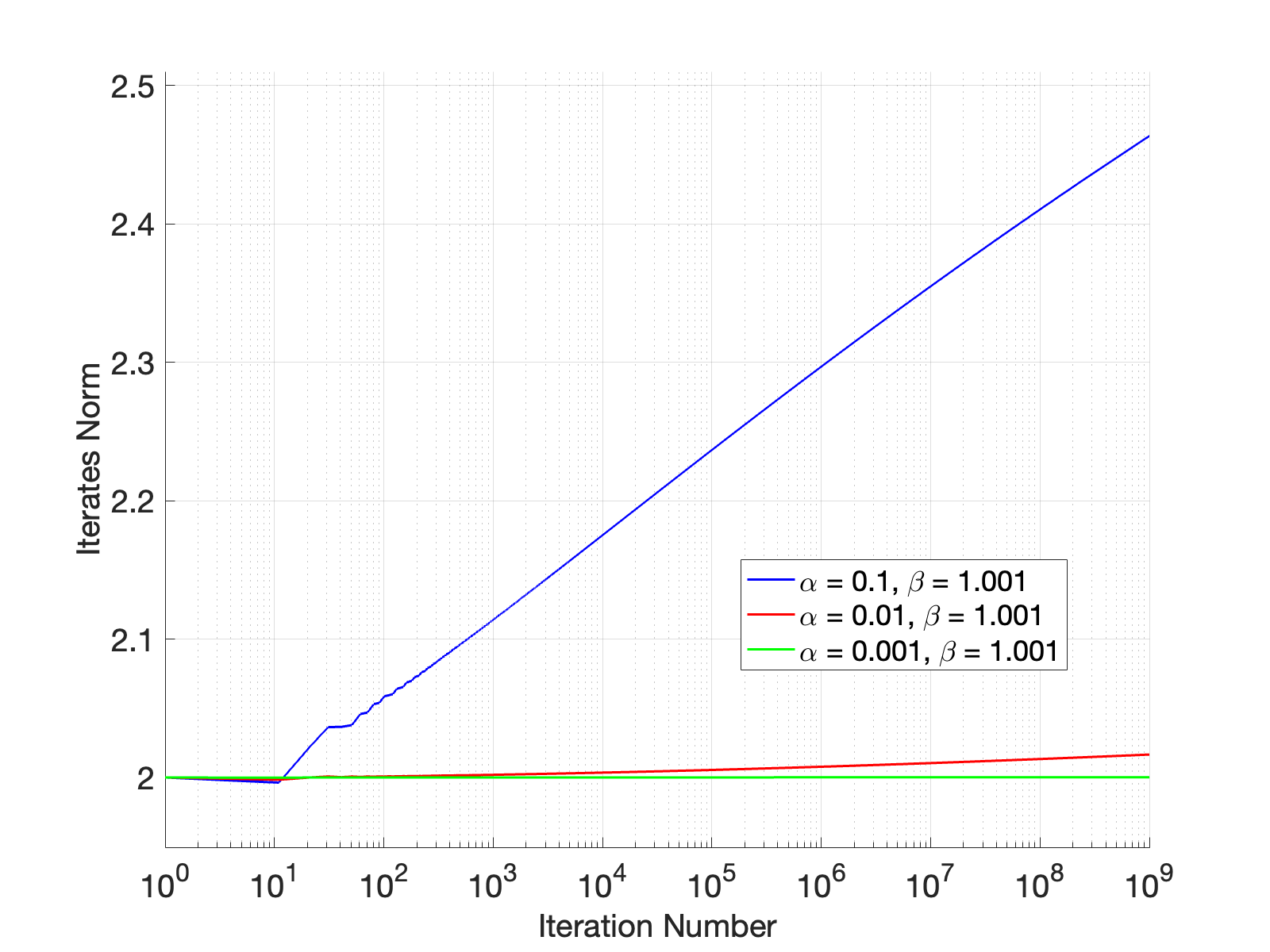}\hspace{-1ex}
\includegraphics[width=0.48\textwidth]{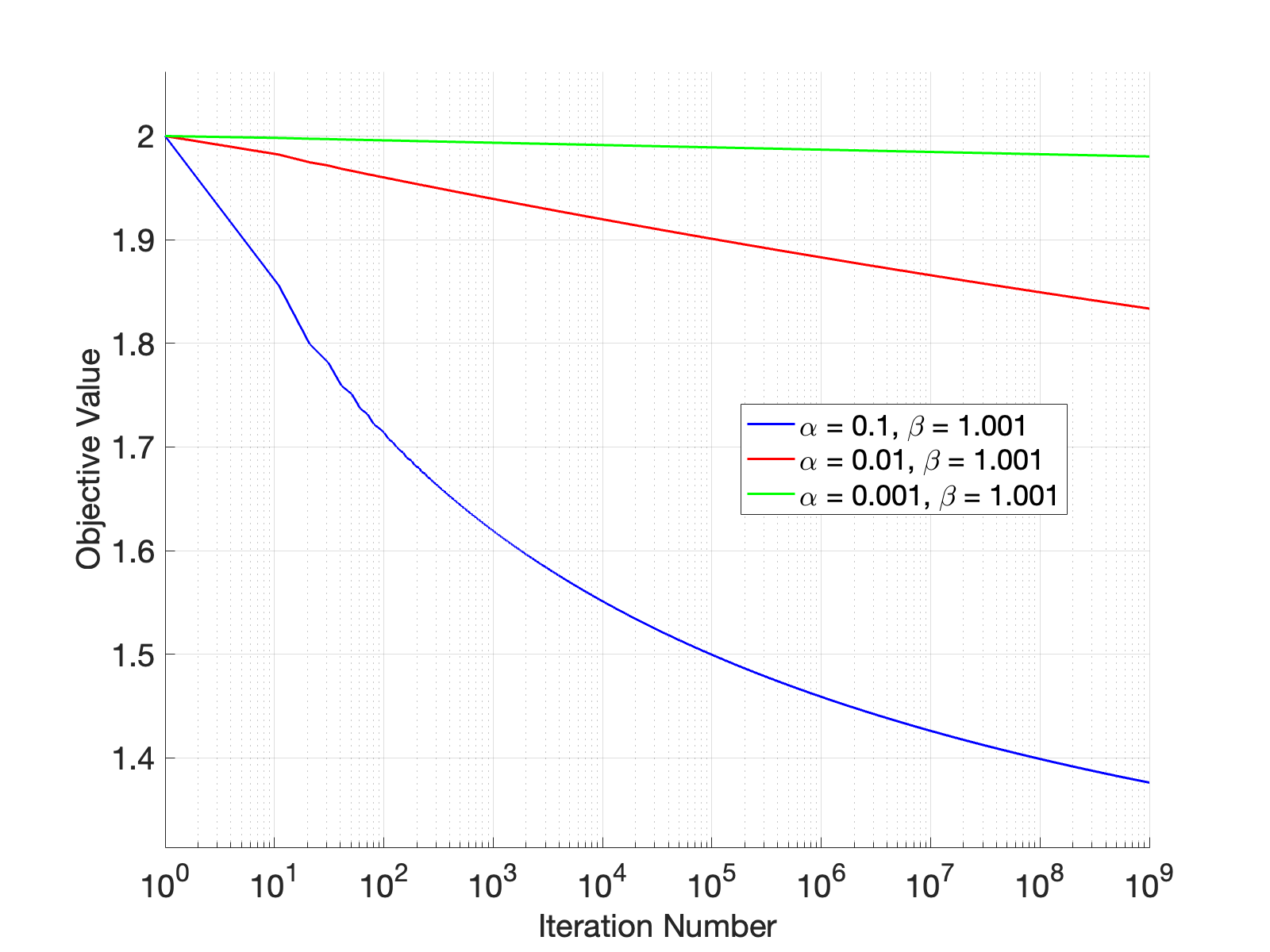}\hspace{-1ex}
\caption{Iterate norms and function values of the subgradient method with summable diminishing step sizes $\alpha_k=\alpha/(k+1)^\beta$ applied to $\ell_1$ matrix completion instance.}
\label{fig:unbdd_summable}
\end{figure}

\end{document}